%% file: main.tex
\documentclass[a4paper,12pt]{article}
\pdfoutput=1
\usepackage[utf8]{inputenc}

\usepackage[margin=0.9in]{geometry}

\usepackage[utf8]{inputenc}
\usepackage{amssymb}
\usepackage{amsmath}
\usepackage{amsfonts}
\usepackage{amsthm}
\usepackage{multirow}
\usepackage{mathrsfs}
\usepackage{authblk}

\DeclareUnicodeCharacter{00A0}{ }
\usepackage{algorithm}
\usepackage{algpseudocode}

\usepackage{graphicx}
\usepackage{epstopdf}
\usepackage{etex}
\usepackage[percent]{overpic}
\usepackage{transparent}
\usepackage{tikz}
 \usepackage{booktabs}
\usepackage{multimedia}
\usetikzlibrary{arrows,shapes}

\usepackage{pst-node}

\usepackage[skins,theorems]{tcolorbox}
\tcbset{highlight math style={enhanced,
  colframe=red,colback=white,arc=0pt,boxrule=1pt}}

\usepackage[customcolors]{hf-tikz}

\tikzstyle{every picture}+=[remember picture]

\usepackage{pgf}
\usepackage{pgfplots}
\pgfplotsset{compat=1.5}
\usepackage{pgfplotstable}

\usepackage{subfig}

\usepackage{refdef}
\usepackage{ifontdef}

\newcommand{\citep}[1]{\cite{#1}}

\usepackage{pgfplots} 
\pgfplotsset{
  tick label style = {font=\color{white!8!black}},
  every axis label = {font=\color{white!8!black}},
  legend style = {font=\color{white!8!black},font=\small},
  label style = {font=\color{white!8!black}}
}
\usepackage{epstopdf}
\usepackage{pdfpages}
\usepackage{tikz}
\usetikzlibrary{shapes,arrows}
\usetikzlibrary{plotmarks}

\newtheorem{prob}{Problem}

 \newcommand{\vek}[1]{\mathchoice{\displaystyle\boldsymbol{#1}}
 {\textstyle\boldsymbol{#1}}{\scriptstyle\boldsymbol{#1}}
 {\scriptscriptstyle\boldsymbol{#1}}}

 \newcommand{\mat}[1]{\mathchoice{\displaystyle\mathbf{#1}}
 {\textstyle\mathbf{#1}}{\scriptstyle\mathbf{#1}}
 {\scriptscriptstyle\mathbf{#1}}}

\title{Bayesian stochastic multi-scale analysis via energy considerations}
\author{M.S. Sarfaraz$^{1}$, B. Rosi\'{c}$^{2}$, H.G. Matthies$^{1}$, A. Ibrahimbegovi\'{c}$^{3}$}
\affil[1]{Institute of Scientific Computing, Muehlenpfordtstrasse 23, TU Braunschweig, 38106 Braunschweig, Germany}
\affil[2]{Applied Mechanics and Data Analysis, Faculty of Engineering Technology, P.O. Box 217, University of Twente, 7500 Enschede, Netherlands}
\affil[3]{Centre de Recherches, Lab. Roberval Mecanique, U.C.T., 60203 Compiegne, France}

\begin{document}

\maketitle

\begin{abstract}
In this paper physical multi-scale processes governed by their own principles for
evolution or equilibrium on each scale are coupled by matching the stored and
dissipated energy, in line with the Hill-Mandel principle.  In our view the
correct representations of stored and dissipated energy is essential to the
representation irreversible material behaviour, and this matching
is also used for upscaling.  The small scales, here the meso-scale,
is assumed to be described probabilistically,
and so on the macroscale also a probabilistic model is identified in a
Bayesian setting, reflecting the randomness of the meso-scale, the loss of
resolution due to upscaling, and the uncertainty involved in the Bayesian process.
In this way multi-scale processes become hierarchical systems in which the information is transferred across the scales by Bayesian identification on coarser levels.  The
quantities to be matched on the coarse-scale model are the stored and dissipated
energies. In this way probability distributions of macro-scale material parameters
are determined, and not only in the elastic region, but also for the irreversible and
highly nonlinear elasto-damage regimes, refelcting the aleatory uncetainty at the
meso-scale level. For this purpose high dimensional meso-scale stochastic simulations in a
non-intrusive functional approximation forms are mapped to the macro-scale models in an
approximative manner by employing a generalised version of the Kalman filter. To reduce the overall computational cost, a model reduction of the meso-scale simulation is achieved by combining the unsupervised learning techniques based on the Bayesian copula variartional inference with the classical functional approximation forms from the field of uncertainty quantification.

\end{abstract}

\input{introduction}

\input{section1}

\input{section2}

\input{section3}

\input{section4}

\input{example1}

\input{example2}

\input{conclusion}

\bibliographystyle{plain}
\bibliography{paper_bib}

\end{document}

%% file: introduction.tex
\section{Introduction}

The predictive modelling of concrete as a 
heterogeneous material requires more 
realistic mathematical models. This is especially true when describing the nonlinear material behaviour, which 
is not fully resolved unless observed on multiple scales 
going from nano-to macro-scale descriptions. 
As the detailed description on the macroscopic level is not computationally
feasable for large scale structures, the multi-scale approaches are often utilised in the numerical practice. In this paper only macro- and meso-scale descriptions of concrete will be considered, however, lower scales can be introduced as well in order to explicitely describe heterogeneities 
characterising the material structure of aggregates, the mortar matrix or the interfacial zone without any further modifications. 

Conceptually, the prevalent computational methods to tackle the multiscale problems can be classified into concurrent and non-concurrent approaches. Concurrent schemes consider both the coarse and fine-scales during the course of the simulation e.g.~FE-squared method \cite{FEYEL2000309,Feyel2003}, whereas the non-concurrent ones are based on a scale separation idea by which the desired quantity of interest (QoI), e.g.~average stresses or strains, are estimated given numerical experiments on a representative volume element (RVE), see \cite{geers2017homogenization} for a recent overview on related techniques. Although, the non-concurrent method has proved to work
very well for elastic properties, it does not explicitly include complex loading paths induced
by structural effect, which are of crucial concern when dealing with material non-linearites
such as plasticity and damage. In addition, the so-called size-effect problem, e.g. a 
problem of determining the size of the representative 
element, appears and has to be resolved. In (martan's thesis) this is achieved by considering the mesh in element approach (MIEL)
in which the meso-scale sturcture is embedded in a macro-scale finite element. Multi-scale here means an
explicit distinction made between microscopic or micro variables – describing the fine
scale representation – and macroscopic or macro variables – describing the coarser-scale.

The main problem in a multi-scale simulation is the process of bridging the scales, especially those of different nature, e.g. discrete versus continuum 
finite element descriptions on the meso-scale and macro-scale, respectivelly. In previous work, see \cite{sadiq}, the infromation transfer is achieved in a
Bayesian probabilistic manner. The meso-scale response is taken 
as measurement, and the properties of the macro-element are assumed to be unknown, and hence uncertain. 
As they cannot be directly measured, their estimation is obtained indirectly from the measurement data. Such an inversion is not 
well posed in a sense of Hadamard, and hence the regularisation is needed. In a Bayesian setting this matches 
with introducing the prior information, i.e.~expert's knowledge or epistemic uncertainty, onto the parameter set. 
This further results in an well-posed problem, the solution of which is stochastic and not deterministic any more. 
The resulting probability distributions are however only epistemic uncertainties and represent our confidence in 
obtained estimates. 
In this paper we go one step further and extend the previous approach, further reffered as the problem of stochastic upscaling, by including a proper 
treatment of aleatory uncertainties used to describe the meso-scale models. These are characterised by variations reflecting in uncertainties describing the geometry, 
the spatial distribution and the material properties of the individual material
constituents and their mutual interaction. 

 Stochastic homogenization is usually performed on an ensemble of RVEs in order to extract the relevant statistical QoI \cite{PW1,PW2,PW3,PW4,PW5,SAVVAS2016340,SAVVAS201847,STEFANOU2015384,
Stefanou2015,STEFANOU2017319}. For example, in \citep{GRAHAMBRADY20102398,liu2016perturbation}, the moving-window approach  is used to characterize the probabilistic uni-axial compressive strength of random micro-structures. Another active direction of research is to develop stochastic surrogate models for strain energy of random micro-structures as in \cite{clement1,clement2,clement3,STABER20171}. The main goal is to mitigate the effect of curse of dimensionality due to large number of stochastic dimensions. With the rapid expansion of machine learning and data driven techniques, the current trend is to train neural network based approximate models \citep{Lu2019,articleBaAnh,2019DeepCM,ungerNnet,UNGER20081994} as a cheaper computational alternative in multi-scale methods. Furthermore, to obtain a probabilistic description of macro-scale characteristic  by incorporating micro-scale measurements, Bayesian methods have been applied to such problems with promising success, please see \cite{franck2017multimodal,
franck2016sparse,
franck2017constitutive,
koutsourelakis2016variational,bruder2018beyond} for recent application in formulation of high dimensional probabilistic inverse problem generally, and estimating distribution of material parameter specifically. Moreover \citep{koutsourelakis2007stochastic,schoberl2018predictive,felsberger2018physics,sadiq} demonstrate the application of Bayesian framework to multi-scale problems.        

This paper serves as an extension of the previously mentioned approaches.  The idea is to design an appropriate macro-scale model, as well as its corresponding parameters such that the meso-scale energy is preserved. For this purpose the non-dissipative and dissipative parts of meso-scale strain energies are captured, and further used as measurements in the Bayesian estimation of the macro-scale properties. Here, three different kinds of algorithms 
are considered. The first addresses direct simulation of the probabilistic meso-scale model, which further results in a high-dimensional strain energy proxy model. The latter one is further mapped to the macro-scale parameter by using the generalised Gauss-Markov-Kalman filter as previously described in \cite{Matthies2016}. Even though the direct simulation allows for the deterministic estimates of the stochastic quantities when treated in a functional approximation setting, this approach is not attractive in real situations when only the meso-scale measurement data are available, or when the parameteric dimension of the model exponentially grows. Therefore, we present another upscaling version in which the non-dissipative and dissipative parts of meso-scale energy are sampled, and further mapped to its lower dimensional causin by an unsupervised data driven learning technique. The idea is to map the non-Guassian energy measurement to lower dimensional Gaussian space 
by a nonlinear mapping. Here this is achieved by emloying the copula based variational inference on a generalised mixture model. Due to high nonlinearity and history dependence such a mapping is not easy to construct, and therefore we employ additional Bayesian correction of the estimate in a sparse form. This is then contrasted to a more practical solution in which the upscaled parameters on the macro-scale level are sampled instead of energies, and further mapped to the lower dimensional Gaussian space. As further discussed this approach is computationally more convenient as the correlation structure 
between the material parameters is easier to learn.

The paper is organised as follows: in section 1 the generalised model problem is presented, and the research questions are defined. Section 2 describes 
Bayesian framework for upscaling random meso-scale information with the particular focus on the approximate posterior estimation. From computational point of view this is further discussed in Section 3. The process of model reduction of the meso-scale computations is further described in Section 4, whereas the energy conservation principle is discussed in Section 5. The algorithm performance is analysed in Section 6 on two different numerical examples, which are further concluded in Section 7.

%% file: section1.tex
\section{Abstract model problem} \label{SSS:set-up}

Let $(\mathcal{Z}_M,\mathcal{E}_M,\mathcal{D}_M)$ be an abstract structure of the general 
rate-independent small-strain homogeneous macro-model in which $\mathcal{Z}_M$ 
denotes the state space, $\mathcal{E}_M: [0,T]\times \mathcal{Z}_M\rightarrow \mathbb{R}$
is a time-dependent energy functional, and $\mathcal{D}_M: \mathcal{Z}_M\times \mathcal{Z}_M \rightarrow [0,\infty]$ is a convex and 
lower-semicontinuous dissipation potential satisfying $\mathcal{D}_M(z,0)=0$ and 
the homogeneity property $\forall z,y \in \mathcal{Z}_M:\textrm{ }\mathcal{D}_M(z,\lambda y)=\lambda \mathcal{D}_M(z,y)$ for all $\lambda>0$.
Then, in an abstract manner
the macro-mechanical system can be described mathematically by the sub-differential inclusion 
\begin{equation}\label{main_eq}
 z:[0,T]\rightarrow \mathcal{Z}:\quad\quad \partial_{\dot{z}} \mathcal{D}_M(\kappa,z,\dot{z})+\textrm{D}_z\mathcal{E}_M(t,\kappa,z) \ni  0
\end{equation}
in which $\textrm{D}_z$ stands for the G\^{a}teux partial differential with respect to the state variable $z$, and the derivative of $\mathcal{D}_M$ 
is given in terms of the set-valued sub-differential $\partial \mathcal{D}_M$ in the sense of the convex analysis, see \cite{mielke}. 
Furthermore, we assume that the rate-independent system in \refeq{main_eq} is parameterised by 
a vector $\kappa$ representing homogeneous material characteristics only. This parameterisation can be extended by including 
boundary conditions, loadings, etc. into the set. However, this situation will not be considered here.

Given mathematical description in \refeq{main_eq}, the goal is to find the unknown set $\kappa$ such that the structural 
response of the macro-scale model matches the 
response of the meso-scale model as well as possible. The latter one is taken as a more detailed 
description of the macro-scale counterpart--the one that accounts for material and geometrical heterogenieties at a lower-scale level, and hence 
 represents the system that
we can evaluate at possibly high cost.
The meso-scale mathematical model reads
\begin{equation}  \label{eq:model-2}
 \partial_{\dot{z}_m} \mathcal{D}_m(z_m,\kappa_m,\dot{z}_m)+\textrm{D}_z\mathcal{E}_m(t,\kappa_m,z_m) \ni  0,
\end{equation}
and is subjected to the same boundary conditions and forcings as the one in \refeq{main_eq} . 
Here, we assume that the state space $\mathcal{Z}_m$, the energy functional $\mathcal{E}_m$ and the dissipation potentials $\mathcal{D}_m$ have same 
form as before with the only difference 
that the material parameters $\kappa$ are given more realistic description $\kappa_m$, and are known. 
However, any 
other kind of model which allows a ``measurement'' resp.\
computation of stored and dissipated energies could be also used.

As $\C{Z}_M \neq \C{Z}_m$, the states $z_M$ and $z_m$ cannot be directly
compared, and the two models are to be compared by some
observables or measurements $y \in \C{Y}$, where $\C{Y}$ is typically some
vector space like $\D{R}^m$.  In other words, let 
\begin{equation}  \label{eq:meas-2}
y_m = Y_m(z_m(\kappa_m,f_m))+\hat{\epsilon},
\end{equation}
be the meso-scale observable (e.g.~energy, stress or strain etc.) 
in which $Y_m$ describes the measurement operator, $f_m$ is the external excitiation and $\hat{\epsilon}$ 
respresents the measurement noise. On the other side, let 
\begin{equation}  \label{eq:meas-11}
y_M = Y_M(\kappa,z_M(\kappa,f_M))
\end{equation}
be the prediction of the same observation on the macro-scale level this time described by the measurement operator $Y_M$
and the external excitation $f_M$ of the same type as $f_m$. To incorporate the prediction of discrepancy $\hat{\epsilon}$
as in \refeq{eq:meas-2}, one has to model 
$y_M$ as a noisy variant. For this purpose we introduce the probability space $(\varOmega_\epsilon,\mathfrak{B}_\epsilon,\mathbb{P}_\epsilon)$ 
and add to $Y_M(\kappa,z_M(\kappa,f_M))$ 
the random variable $\epsilon(\omega_\epsilon)\in L_2(\varOmega_\epsilon,\mathfrak{B}_\epsilon,\mathbb{P}_\epsilon;\mathbb{R}^d)$ 
that best describes our knowledge about $\hat{\epsilon}$. 
Hence, \refeq{eq:meas-1} becomes stochastic and reads
\begin{equation}  \label{eq:meas-12}
y_M(\omega_\epsilon) = Y_M(\kappa,z_M(\kappa,f_M))+\epsilon(\omega_\epsilon).
\end{equation}
Typically, $\epsilon(\omega_\epsilon)$ is modelled as a zero-mean Gaussian random variable
$\epsilon\sim \mathcal{N}(0,C_\epsilon)$ with covariance $C_\epsilon$.
However, other models for $\epsilon(\omega_\epsilon)$ can also be introduced without modifying the general setting presented in this paper.

In variaty of literature the measurement in \refeq{eq:meas-2} is considered as deteriministic, e.g. 
the measurement is a function 
of the state $z_m$ characterised by one 
realisation of the heterogeneous material. This further defines the first problem of our consideration:

\begin{prob}\label{Prob1}
 Find deterministic $\kappa$ in \refeq{main_eq} such that the predictions of \refeq{eq:meas-12} match those of \refeq{eq:meas-2} in a measurement sense.
\end{prob}

 However, in a real multi-scale analysis the parameters $\kappa_m$ in \refeq{eq:model-2} vary, or are not fully known. 
Hence, the meso-scale model is characterised by aleatory and epistemic uncertainties, which have to be encountered into the modelling process.
Taking the probabilistic view on uncertainty, we model $\kappa_m$ as a random variable/field in $L_2(\varOmega_\kappa,\mathfrak{B}_\kappa,
\mathbb{P}_\kappa;\mathcal{K}_m)$ defined by mapping 
\begin{equation}
 \kappa_m(\omega_\kappa):=\kappa_m(x,\omega_\kappa): \mathcal{G}\times \varOmega_\kappa \mapsto \mathcal{K}_m.
\end{equation}
Here, $\mathcal{K}_m$ is the parameter space which depends on the application. 
As a consequence, the evolution problem described by $(\mathcal{Z}_m,\mathcal{E}_m,\mathcal{D}_m)$ also becomes uncertain, and hence 
\refeq{eq:model-2} rewrites to 
\begin{equation}  \label{eq:model-2s}
 \partial_{\dot{z}_m} \mathcal{D}_{m}(z_m(\omega_\kappa),\kappa_m(\omega_\kappa),\dot{z}_m(\omega_\kappa))+
 \textrm{D}_z\mathcal{E}_{m}(t,\omega,\kappa_m(x,\omega_\kappa),z_m(\omega_\kappa)) \ni  0 \textrm{ a.s.}
\end{equation}
in which 
\begin{eqnarray}
 \mathcal{E}_{m}(t,z_m(\omega_\kappa))=\int_{\varOmega_\kappa} \mathcal{E}_m(\kappa_m(x,\omega_\kappa),z_m(\omega_\kappa)) \mathbb{P}(d\omega_\kappa),&&\nonumber \\
 \mathcal{D}_{m}(z_m(\omega_\kappa),\dot{z}_m(\omega_\kappa))=\int_{\varOmega_\kappa} \mathcal{D}_m(\kappa_m(\omega_\kappa),z_m(\omega_\kappa),\dot{z}_m(\omega_\kappa))
 \mathbb{P}(d\omega_k),&&
\end{eqnarray}
respectively. 

Once the uncertainty is present in the meso-scale model, the observation in \refeq{eq:meas-2} rewrites to 
\begin{equation}  \label{eq:meas-21}
y_m(\omega_y) = Y_m(z_m(\kappa_m(\omega_\kappa),f_m(\omega_\kappa)))+\hat{\epsilon}(\omega_{e}), \quad \omega_y:=(\omega_\kappa,\omega_e)
\end{equation}
in which both the model and the error $\hat{\epsilon}(\omega_e) \in L_2(\varOmega_e,\mathfrak{B}_e,\mathbb{P}_e;\mathbb{R}^d)$ are stochastic. 
This in turn modifes Problem \ref{Prob1} into

\begin{prob}\label{Prob2}
 Find stochastic $\kappa_M$ in \refeq{main_eq} such that the predictions of \refeq{eq:meas-12} match those of \refeq{eq:meas-21} in a measurement sense.
\end{prob}

The upscaling process that is related to Problem \ref{Prob1} is already considered in \cite{sadiq}, and hence will not be repeated here. However, 
as we show later this problem 
is a special case of Problem \ref{Prob2} that is the main topic of this paper.

%% file: section2.tex
\section{Bayesian upscaling of random mesostructures}

Let the meso-scale parameter set $\kappa_m(\omega_\kappa)$ define the observations in \refeq{eq:meas-21}, here assumed to be continous. The goal is to use information in
\refeq{eq:meas-21} in order to calibrate (upscale) the set of material parameters $\kappa$
of the coarse-scale model in \refeq{main_eq}. 
To achieve this, $\kappa$ is assumed to be uncertain (unknown) and further modelled a priori as a random variable $\kappa(\omega)$
---prior---belonging to $L_2(\Omega,\mathcal{F},\mathbb{P};\mathcal{K})$. Hence, the model in \refeq{main_eq} rewrites to 
\begin{equation}\label{main_eqs}
 \partial_{\dot{z}} \mathcal{D}_M(\kappa(\omega),z(\kappa(\omega)),\dot{z}(\kappa(\omega)))+
 \textrm{D}_z\mathcal{E}_M(t,\kappa(\omega),z(\kappa(\omega))) \ni  0,
\end{equation}
and subsequentually a priori prediction of the macro-scale measurement becomes
\begin{equation}  \label{eq:meas-1}
y_M(\kappa(\omega),\epsilon(\omega_\epsilon)) = Y_M(\kappa(\omega),z_M(\kappa(\omega)),f_M))+\epsilon(\omega_\epsilon)
\end{equation}
with $\epsilon(\omega_\epsilon)\in L_2(\varOmega_\epsilon,\mathfrak{B}_\epsilon,\mathbb{P}_\epsilon)$. 
The goal is to identify the vector $\kappa(\omega)$ given $y_m$ using Bayes's rule such that
\begin{equation}\label{bayes_rule}
 p(\kappa|y_m)=\frac{p(y_m|\kappa)p(\kappa)}{p(y_m)}
\end{equation}
holds. 
As $\kappa$ is positive definite, this constraint has to be taken into consideration. Therefore, instead of Problem \ref{Prob1} and 
Problem \ref{Prob2} we consider their modified versions:

\begin{prob}\label{Prob4}
 Find deterministic $q:=\textrm{log }\kappa$ in \refeq{main_eq} such that the predictions of \refeq{eq:meas-12} match those of \refeq{eq:meas-2} in a measurement sense.
\end{prob}
and
\begin{prob}\label{Prob3}
 Find stochastic $q:=\textrm{log }\kappa$ in \refeq{main_eq} such that the predictions of \refeq{eq:meas-12} match those of \refeq{eq:meas-21} in a measurement sense.
\end{prob}
Instead of calibrating the vector $\kappa$ directly, we calibrate its logartithm as in this way computationally whatever approximations or
linear operations are
performed on the numerical representation of $q(x,\omega)$, in the end $\textrm{exp}(q(x, \omega))$
is always going to be positive.  Additionally, the multiplicative group of positive
real numbers—a (commutative) one-dimensional Lie group—is thereby put
into correspondence with the additive group of reals, which also represents
the (one-dimensional) tangent vector space at the group unit, the number
one. This is the corresponding Lie algebra. A positive quadratic form on the
Lie algebra—in one dimension necessarily proportional to Euclidean distance
squared—can thereby be carried to a Riemannian metric on the Lie group. Therefore,
\begin{equation}\label{bayes_rule_log}
 p(q|y_m)=\frac{p(y_m|q)p(q)}{p(y_m)}.
\end{equation}
However, in practice the estimation of the full posterior $p(q|y_m)$ is not analytically tractable. 
Numerical estimation on the other hand is expensive either due to evidence estimation or 
due to slow convergence of the random walk algorithms. 
As in engineering practice one is often not interested in estimating the full posterior measure,
in this paper we investigate the estimation of the posterior functional given in a form of conditional expectation
\begin{equation}\label{cond_exp_integration}
 \mathbb{E}(q|y_m)=\int_\varOmega q p(q|y_m) d q
\end{equation}
as it is computationally simpler. Instead of direct intergation over the posterior measure, the conditional expectation can be straightforwardly estimated by projecting the random variable $q$ onto the subspace generated by the sub-sigma algebra $\mathfrak{B}:=\sigma(Y)$ of measurement. 
To achieve this, one has to compute the minimal distance of $q$ to the point $q^*$ which can be achieved in different ways. As shown by \cite{banerjee,Ro19}, the notion of distance can be generalised given a strictly convex, differentiable function
$\varphi: \mathbb{R}^d \mapsto \mathbb{R}$ with the 
hyperplane tangent $\mathcal{H}_{\hat{q}}(q)=\varphi(\hat{q})+\langle q-\hat{q},\nabla \varphi(\hat{q})\rangle$ to $\varphi$ at point $\hat{q}$, such that
\begin{equation}
\label{eq:optimality_blf}
 q^*:=\mathbb{E}(q|\mathfrak{B})=
 \underset{\hat{q}\in L_2(\varOmega,\mathcal{B},\mathbb{P};\mathcal{Q})}{\arg \min}
 \textrm{} \mathbb{E}(\mathcal{D}_\varphi(q||\hat{q}))
\end{equation}
holds. Here, $\mathcal{D}_\varphi(q||\hat{q})=\mathcal{H}_{q}(q)-\mathcal{H}_{\hat{q}}(q)$ denotes the distance term that is also known as the Bregman's loss function (BLF) or divergence. In general the projection in \refeq{eq:optimality_blf} is of non-orthogonal kind and reflects the generalised inequality
\begin{equation}
\label{bregman_var_ineq}
 \mathbb{E}(D_\varphi(q||\check{q}))\geq \mathbb{E}(D_\varphi(q||q^*))+\mathbb{E}(D_\varphi(q^*||\check{q}))
\end{equation}
that also holds for any arbitrary $\mathfrak{F}$-measurable random variable $\check{q}$. In case when $q^*=\mathbb{E}(q)$ the previous relation rewrites to
\begin{equation}
  \mathbb{E}(\mathcal{D}_\varphi(q||\mathbb{E}(q)))=\mathbb{E}(\mathcal{D}_\varphi(q||\check{q}))
  -\mathbb{E}(\mathcal{D}_\varphi(\mathbb{E}(q)||\check{q}))\geq 0,
  \end{equation}
in which the distance $\mathbb{E}(\mathcal{D}_\varphi(q||\mathbb{E}(q)))$ has notion of variance, further refered as Bregman's variance. 
In terms of convex function $\varphi$ the Bregman's variance obtains the following form
  \begin{eqnarray}\label{jansen}
 \textrm{var}_\varphi(q):&=&\mathbb{E}(\mathcal{D}_\varphi(q||\mathbb{E}(q)))=\mathbb{E}(\varphi(q)-\varphi(\mathbb{E}(q))
 +\langle q-\mathbb{E}(q),\nabla \varphi(\mathbb{E}(q))\rangle)\nonumber\\
 &=&\mathbb{E}(\varphi(q))-\varphi(\mathbb{E}(q))\geq 0,
\end{eqnarray}
which is also matching the definition of Jensen's inequality in the context of the probability theory. By minimising the Bregmann's variance in the last expression,
one obtains the mean as the corresponding minimumm. Subsequentually, one may conclude that the mean is the same minimum point for any expected Bregman divergence, e.g.
  \begin{equation}
   \mathbb{E}(q)=\underset{\hat{q}\in L_2(\varOmega,\mathcal{F},\mathbb{P};\mathcal{Q})}{\arg \min} 
   \mathbb{E}(\mathcal{D}_\varphi(q||\hat{q})).
  \end{equation}
Notice that similar can be derived for $\mathfrak{B}$-measurable random variables. The last relation then matches \refeq{eq:optimality_blf}.

For computational purposes, the Bregman's distance in \refeq{eq:optimality_blf} is further taken as the squared Euclidean distance by assuming $\varphi(q)=\frac{1}{2}q^2$, such that 
$\mathcal{D}_\varphi(q||\hat{q})=\frac{1}{2}\|q-\hat{q}\|^2$ and
\begin{equation}
\label{eq:optimality_blf1}
 q^*:=\mathbb{E}(q|\mathfrak{B})=
 \underset{\hat{q}\in L_2(\varOmega,\mathcal{B},\mathbb{P};\mathcal{Q})}{\arg \min} \textrm{} \mathbb{E}(\|q-\hat{q}\|^2).
\end{equation}
In this manner \refeq{bregman_var_ineq} reduces to the Pythagorean theorem
\begin{equation}
 \|q-\hat{q}\|^2=\|q-q^*\|^2+\|q^*-\hat{q}\|^2
\end{equation}
%
in which by taking  $q^*=\mathbb{E}(q|\mathcal{B})$ one obtains
\begin{equation}
 \|q-\hat{q}\|^2=\|q-\mathbb{E}(q|\mathcal{B})\|^2+\|\mathbb{E}(q|\mathcal{B})-\hat{q}\|^2.
\end{equation}
Integrating over the probability measure of $q$, the previous equality can be rewritten as inequality:
\begin{equation}
 \mathbb{E}(\|q-\hat{q}\|^2)\geq \mathbb{E}(\|q-\mathbb{E}(q|\mathcal{B})\|^2)
\end{equation}
which shows that the conditional expectation is the minimiser of the mean squared error. In addition, if $\hat{q}=\mathbb{E}(q)$ then 
it is also the minimum variance unbiased estimator. 

\subsection{Upscaling based on the conditional expectation approximation}

Following \refeq{eq:optimality_blf1}, one may 
 decompose the random variable $q$ belonging to  $(\varOmega,\mathfrak{F},\mathbb{P})$ into projected $q_p$ and residual
 $q_r$
components such that
\begin{equation}\label{rv_decomposition}
 q=q_p+q_o=P_\mathfrak{B} q+(I-P_\mathfrak{B})q
\end{equation}
holds. Here, $q_p=P_\mathfrak{B} q=\mathbb{E}(q|\mathfrak{B})$ is the orthogonal projection of the random variable $q$
onto the
space $(\varOmega,\mathfrak{B},\mathbb{P})$ of all distributions consistent with the data,
whereas $q_o:=(I-P_\mathfrak{B})q$
is its orthogonal residual.
%
%
%
%
%
To give \refeq{rv_decomposition} more practical form,  the projection term $P_\mathfrak{B} q$  is further described by a measurable mapping $\phi: y \mapsto q $
according to the Doob-Dynkin lemma which states:
\begin{equation}\label{Doobmap}
 \mathbb{E}(q|y)=\phi(y(q))=\phi\circ y\circ q.
\end{equation} 
As a result, \refeq{rv_decomposition} rewrites to:
\begin{equation}\label{rv_decomp1}
 q=\phi(y(q))+(q-\phi(y(q)))
\end{equation} 
in which 
the first term in the sum, i.e.~$\phi(y(q))$, is taken to be the projection of $q$ onto the meso-scale
data set $y_m$, whereas $(q-\phi(y(q)))$ 
is the residual component defined by a priori knowledge $q_f$.
Following this, \refeq{rv_decomp1} recasts to the update equation for the random variable $q$ as
\begin{equation}\label{general_map}
 q_a=q_f+\phi(y_m)-\phi(y_M)
\end{equation}
in which $y_M$ (see \refeq{eq:meas-12}) is the random variable representing our prior prediction/forecast of the measurement data, and $q_a$ is 
the assimilated random variable. 
%
%
Therefore, to estimate $q_a$ one requires only information on 
the map $\phi$.

 For the sake of computational simplicity, the map in \refeq{Doobmap} is further approximated in a Galerkin manner by 
\begin{equation}\label{poly_poly}
 \mathcal{Q}_n=\{ \phi(y) \in \mathcal{Q} \textrm{ }|\textrm{ } \phi_n: y \rightarrow q
 \textrm{ a n-th degree polynomial}\}
\end{equation}
such that the filter in \refeq{general_map} rewrites to
\begin{equation}\label{update_eq}
 q_a(\omega)=q_f(\omega)+\phi_n(y_m)-\phi_n(y_M).
\end{equation}
As the map in \refeq{poly_poly} is parametrised by a set of coefficients $\beta$, i.e. $\phi_n(y;\beta)$, these further can be found
by minimising the residual component (the optimality condition in \refeq{eq:optimality_blf1}):
%
%
\begin{equation}
\label{optim_probv}
 \beta^*=\underset{\beta}{\arg \min }\textrm{ } \mathbb{E}(\|q_f-\phi_n(y_M;\beta)\|^2).
\end{equation}
In an affine case, when $n=1,\phi_1(y_M;\beta)=Ky+b$, the previous optimisation outcomes in a formula defining 
the well-known Kalman gain:
\begin{equation}
 K=\textrm{cov}_{q_f,y_M}\textrm{cov}_{y_M}^{-1}
\end{equation}
such that the update formula in \refeq{update_eq} reduces to the generalisation of the Kalman filter
\begin{equation}
\label{lin_update}
 q_a=q_f+K(y_m-y_M),
\end{equation}
here refered as Gauss-Markov-Kalman filter, 
for more details please see \cite{Matthies2016}. 

However, in most of practical situations one has to deal with the nonlinearity of the measurement-parameter map, and in such a case
the formula given in \refeq{lin_update} is not optimal. One example is the process of upscaling in which the energy observation is used to predict the material characteristics.
To overcome this issue, one may introduce higher order terms in the approximation in \refeq{poly_poly} either in
a monomial form
\begin{equation}
 \phi(y;\beta)=\sum_{i=0}^p K_i(y^{\vee i}), \quad \beta:=(K_i)
\end{equation}
in which $y^{\vee i}:=\textrm{Sym}(y^{\otimes i})$ 
is the symmetric tensor product of $y$ taken $i$ times with itself, or in a polynomial one
\begin{equation}
 \phi(y;\beta)=\sum_{\alpha \in \mathcal{I}} K^{(\alpha)}V_\alpha( y)
\end{equation}
in which $\mathcal{I}$ is the multi-index set with elements $\alpha:=(\alpha_1,...,\alpha_n)$, and $V_\alpha$ 
are multivariate polynomials possibly of orthogonal kind. Following this, one may distingusih the monomial filtering formula
\begin{equation}
\label{high_update_mono}
 q_a=q_f+\sum_{i=0}^p K_i(y_m^{\vee i}-y_M^{\vee i})
\end{equation}
from the polynomial one
\begin{equation}
\label{high_update_poly}
 q_a=q_f+\sum_{\alpha \in \mathcal{I}} K^{(\alpha)}(V_\alpha( y_{m})- V_\alpha(y_{M})).
\end{equation}
 
Note that the form of \refeq{high_update_poly}, given in monomials, is
numerically not a good form—except for very low orders $n$
—and hence straightforward
use in computations is not recommended.

Finally, in a similar manner as before one may estimate the unknown coefficients $\beta$
given the stationarity or orthogonality condition in \refeq{optim_probv}. The uniquness of the solution as well as different forms of 
approximations are studied in 
\cite{Matthies2016}.

\subsubsection{Bayesian Gauss-Markov-Kalman filter}\label{BGMK}

The filter presented in the previous section is known to be computationally expensive when the dimension of $q_f$ and $y_M$ is high. Therefore, one may substitute $y_M$ by a surrogate model
\begin{equation}\label{forward_model_map}
 \hat{y}_M=\phi_f(q_f;\beta)
\end{equation}
in which $\phi_f$ is usually taken to be of polynomial type in arguments $q_f$ with the coefficients $\beta$. 
To estimate $\beta$ one may use the same type of estimator as described in the previous section, however, in a slightly different form, as this time 
one estimates the map
$q_f \mapsto y_M$, i.e.
\begin{equation}
 \label{est_back_map}
 \beta^*=\underset{\beta}{\arg \min }\textrm{ } \mathbb{E}(\|y_M-\phi_f(q_f;\beta)\|^2)
\end{equation}
as further described in details in \cite{Ro19} on an example of ordinary differential equation. However, 
the estimate provided by 
\refeq{est_back_map} requires the knowledge of the random variable $y_M$. Although we may estimate $y_M$ by propagating forward the uncertainty through the meso-scale model, this is not very efficient due to the high-dimensionality of the problem. Therefore, we assume that $y_M$ is in general only 
known as a set of samples, and our goal is 
to match
$\hat{y}_M$ with $y_M$ in a distribution sense. 
To achieve this, we consider another 
type of Bregman's 
divergence in which $\varphi(y)=\mathbb{E}(\pi_y \textrm{log } \pi_y )$ 
is the continous entropy 
with $\pi_y$ being the density distribution 
of the random variable $y$. In such a case \refeq{est_back_map} rewrites to
\begin{equation}\label{parameterbeta}
 \beta^*=\underset{\hat{y}_M:=\phi_f(q,\beta)}{\arg \min }\textrm{ } D_{KL}(y_M||\hat{y}_M)
\end{equation}
with 
\begin{equation}
\label{kl_in_q}
 D_{KL}(y_M||\hat{y}_M):=\int \pi_{y_M} \textrm{ log }\frac{\pi_{y_M}}{\pi_{\hat{y}_M}} dy_M
\end{equation}
being the Kullback-Leibler (KL) divergence,
i.e. the distance betwen the prior variable $y_M$ and its approximation $\hat{y}_M$ obtained
by propagation of $q_f$ to $y_M$ via nonlinear map 
$\phi_f(q_f,\beta)$ parameterised by $\beta$. 
However, as the backward map $\phi_f$ is not explicitely known,
one may approximate $\pi_{\hat{y}_M}$ by
conditional density $\pi_{\hat{y}_M|x}$ coming from Bayes's rule, see \cite{Ro19}
\begin{equation}
\label{bayes_f}
\pi_{\hat{y}_M|x}=\frac{\pi_{\hat{y}_M,x}}{p_x}
\end{equation}
in which $x:=(q_f(\omega_i),y_M(\omega_i))_{i=1}^N$ is the full set of data (samples)
describing the forward propagation of $q_f$ to $y_M$. In order to specify $\hat{y}_M$ 
the only thing we need to find are the coefficients of the map $\phi_f$. Therefore, we infer $\beta$ given data $x$ using Bayes's rule
\begin{equation}
\label{bayes_coeff}
 p(\beta|x)=\frac{p(x,\beta)}{\int p(x,\beta) d\beta}
\end{equation}
instead of \refeq{bayes_f}. 
In general case the marginalisation in \refeq{bayes_coeff} can be expensive,
and therefore in this paper we use the variational Bayesian
inference instead.  The idea is to introduce a family $\mathcal{D}:=\{g(\beta):=g(\beta|\lambda,w)\}$
over $\beta$ indexed by a set of free parameters $(w,\lambda)$
such that $\hat{y}_M \sim y_M$. Thus, the idea is to optimise the parameter values by 
minimising the Kullback-Leibler divergence
\begin{equation}
 g^*(\beta)=\underset{g(\beta)\in \mathcal{D}}{\arg \min \textrm{ }} D_{KL}(g(\beta)||p(\beta|x))=
 \underset{g(\beta)\in \mathcal{D}}{\arg \min} \int g(\beta) \textrm{log }\frac{g(\beta)}
 {p(\beta|x) }d\beta
\end{equation}
instead of \refeq{parameterbeta}. After few derivation steps as depicted in \cite{Hoffman}, the previous minimisation problem reduces to 
\begin{equation}\label{min_elbo}
 \beta^*=\arg \max \mathcal{L}(g(\beta)):= \mathbb{E}_g(\textrm{log }p(x,\beta))-\mathbb{E}_g(\textrm{log }g(\beta))
\end{equation}
in which $\mathcal{L}(g)$ is the evidence lower 
bound (ELBO), or variational free energy. 
To obtain closed-form solution for $\beta^*$, the usual practice is to assume that both
posterior $p(\beta|x)$ as well as its approximation $g(\beta)$ 
can be factorised in a 
mean sense, i.e.
\begin{equation}
 p(\beta|x)=\prod p(\beta_i|x), \quad 
 g(\beta)=\prod g(\beta_i)
\end{equation}
in which each factor $p(\beta_i|x)$, $g(\beta_i)$ is independent and belongs to the exponential family.
Similarly, their complete conditionals given all other variables and observations are also assumed to belong 
to the exponential families and are assumed to be independent. Obviously these assumptions 
lead to conjugacy relationships, and closed form solution of \refeq{min_elbo} as further discussed in more detail in \cite{Hoffman}.

\subsection{Upscaling based on the polynomial chaos approximation}

To discretise RVs in \refeq{lin_update}, \refeq{high_update_mono} and \refeq{high_update_poly} we use {functional approximations}.  This means that all RVs,
say $y_M(\omega)$, are described as functions of {known} RVs
$\{\theta_1(\omega),\dots,\theta_n(\omega),\dots\}$.  Often, when for example
stochastic processes or random fields are involved, one has to deal here
with {infinitely} many RVs, which for an actual computation have to be
truncated to a finte vector $\vek{\theta}(\omega)=[\theta_1(\omega),\dots,\theta_L(\omega)]\in
\varTheta \cong \D{R}^L$ of significant RVs.
We shall assume that these have been chosen such as to be
independent, and often even normalised Gaussian and independent.  The reason to not use $q_f$ directly 
is that in the process of identification of
$q$  they may turn out to be correlated, whereas $\vek{\theta}$ can stay
independent as they are.  

To actually describe the functions $y_M(\vek{\theta}), q_f(\vek{\theta})$,
one further chooses a finite set of linearly independent functions 
$\{\Psi_\alpha\}_{\alpha\in\C{J}_Z}$ of the variables $\vek{\theta}(\omega)$, where
the index $\alpha=(\dots,\alpha_k,\dots)$ often is a {multi-index},
and the set of multi-indices for approximation $\C{J}_Z$ is a finite set with
cardinality (size) $Z$.  Many different systems of functions can be used, classical
choices are \citep{wiener38, ghanemSpanos91,
xiu2010} multivariate polynomials ---
leading to the {polynomial chaos expansion} (PCE) or generalised PCE (gPCE)
\citep{xiu2010},
kernel functions, radial basis functions, or functions derived
from fuzzy sets.  The particular choice is immaterial for the further development.
But to obtain results which match the above theory as regards $\E{L}$-invariant
subspaces, we shall assume that the set $\{\Psi_\alpha\}_{\alpha\in\C{J}_Z}$ includes
all the {linear} functions of $\vek{\theta}$.  This is easy to achieve with
polynomials, and w.r.t\ kriging it corresponds to {universal} kriging.
All other functions systems can also be augmented by a linear trend.

Thus a RV $\vek{y}_M(\vek{\theta})$ would be replaced by a {functional approximation}
--- this gives these methods its name, sometimes also termed {spectral} approximation ---
\begin{equation}  \label{eq:FA}
   \vek{y}_M(\vek{\theta}) = \sum_{\alpha\in\C{J}_Z} \vek{y}_M^{(\alpha)} 
     \Psi_\alpha(\vek{\theta}) = \vek{y}_M(\vek{\theta}), 
\end{equation}
 and analogously $\vek{q}_f$ by
\begin{equation}  \label{eq:FA-x}
    \vek{q}_f(\vek{\theta}) = 
   \sum_{\alpha\in\C{J}_Z} \vek{q}_f^{(\alpha)} \Psi_\alpha(\vek{\theta}).
\end{equation}

Like always, there are several alternatives to determine the coefficients
$\vek{q}_f^{(\alpha)},\vek{y}_M^{(\alpha)}$ in the previous equations. In this paper, they are estimated in a data-driven way by 
using the variational method presented in Section \ref{BGMK}. Namely, \refeq{eq:FA} is taken as an example of \refeq{forward_model_map}
in which the coefficients $\vek{v}:=\{\vek{y}_M^{(\alpha)}\}_{\alpha \in \mathcal{J}_Z}$ are estimated by minimising ELBO analogous to the one in \refeq{min_elbo}
by using the variational relevance vector machine method \cite{bishopvar}. Namely, the measurement forecast $\vek{y}_s:=\{\vek{y}_M(\vek{\theta}_j)\}_{j=1}^N $ can be rewritten in a vector form as
\begin{equation}\label{regression}
\vek{y}_s=\vek{v}\mat{\Psi}
\end{equation}
in which $\mat{\Psi}$ is the matrix of collection of basis functions $\Psi_\alpha(\vek{\theta})$ evaluated at the set of sample points $\{\vek{\theta}_j\}_{i=1}^N$. However, the expression in the previous equation is not complete as the PCE  in \refeq{eq:FA} is truncated. This implies presence of the modelling error. Under Gaussian assumption, the data then can be modelled as
\begin{equation}
p(\vek{y}_s)\sim \mathcal{N}(\vek{v}\mat{\Psi},\varsigma^{-1}\vek{I})
\end{equation}
in which $\varsigma \sim \Gamma(a_\varsigma,b_\varsigma)$ denotes the imprecision parameter, here assumed to follow Gamma distribution. 
The coefficients $\vek{v}$ are given a normal distribution under the independency assumption:
\begin{equation}
p(\vek{v}|\vek{a})\sim \prod_{i=0}^Z \mathcal{N}(0,\zeta_i^{-1})
\end{equation}
in which $Z$ denotes the cardinality of the PCE, and $\vek{\zeta}:=\{\zeta_i\}$ is a vector of hyperparameters. To promote for sparsity, 
the vector of hyperparameters is further assumed to follow Gamma distribution
\begin{equation}
p(\zeta_i) \sim \Gamma(a_{i},b_i)
\end{equation}
under the independency assumption. In this manner the posterior for $\vek{\beta}:=\{\vek{v},\vek{\zeta},\vek{\varsigma})$, i.e. $p(\vek{\beta}|\vek{y}_s)$, can be approximated by a variational mean field form 
\begin{equation}
g(\vek{\beta})=g(\vek{v})g(\vek{\zeta})g(\vek{\varsigma}),
\end{equation}
 the factors of which are chosen to take same distribution type as the corresponding prior due to the conjugacy reasons. Once this assumption is made, one may maximise the corresponding ELBO in order to estimate the parameter set.

Following this the filtering equation as presented in \refeq{general_map} obtains its purely deterministic flavor. 
Once the random variables of consideraton are substituted by their discrete versions one obtains the functional approximation based filter 
\begin{equation}\label{hit_eq}
\sum_{\alpha \in \mathcal{J}} \vek{q}_a^{(\alpha)}\varPsi_\alpha(\vek{\theta}(\omega))=\sum_{\alpha \in \mathcal{J}}
\vek{q}_f^{(\alpha)}\varPsi_\alpha(\vek{\theta}(\omega))+
\varphi_n(\vek{y}_m)-\varphi_n(\sum_{\alpha \in \mathcal{J}} \vek{y}_M^{(\alpha)}\varPsi_\alpha(\vek{\theta}(\omega)))
\end{equation}
which is purely deterministic. 
If one choses $\varPsi_\alpha(\vek{\theta}(\omega))$ as orthogonal polynomial chaos basis functions, one may project the previous equation onto
$\varPsi_\beta(\vek{\theta}(\omega))$ to obtain
the index-wise formula
\begin{equation}\label{eqn:expansoin_update}
\forall \beta \in \mathcal{J}:\quad \vek{q}_{a}^{(\beta)}=\vek{q}_{f}^{(\beta)}+\mathbb{E}\left(\varphi_n(\vek{y}_m)-
\varphi_n\left(\sum_{\alpha \in \mathcal{J}} \vek{y}_M^{(\alpha)}\varPsi_\alpha)\right)\Psi_\beta\right).
\end{equation}
In affine case the previous equation reduces to
\begin{equation}\label{eqn:expansoin_updatea}
\forall \alpha \in \mathcal{J}:\quad \vek{q}_{a}^{(\alpha)}=\vek{q}_{f}^{(\alpha)}+\vek{K}_g(\vek{y}_m^{(\alpha)}
-\vek{y}_{{f}}^{(\alpha)}),
\end{equation}
with $\vek{K}_g$ being the well-known Kalman gain defined as
\begin{equation}\label{eqn:kalman_gain}
\vek{K}_g=\vek{C}_{q_{f}y_{M}}(\vek{C}_{y_{M}}+\vek{C}_{\epsilon})^{-1},
\end{equation}
and evaluated directly from the polynomial chaos coefficients. In particular, the formula for
covariance matrix $\vek{C}_{y_{M}}$ reads
\begin{equation}
 \vek{C}_{y_{M}}=\sum_{\alpha>0} \vek{y}_M^{(\alpha)}\otimes \vek{y}_M^{(\alpha)}\alpha!,
\end{equation}
and analogously can be defined for other covariance matrices appearing in \refeq{eqn:kalman_gain}.

%% file: section3.tex
\section{Upscaling aleatory uncertainty}

The main issue in \refeq{bayes_rule_log} is that $q$ is not deterministic as 
in most of cases found in literature. In contrast, the vector $q$ is a
random variable of unknown probability distribution, and represents the aleatory uncertainty 
reflected in the measurement data $y_m$. The formula as given in \refeq{eqn:expansoin_update} is in general used for 
the estimation of unknown $q$ given deterministic measurement $y_m$, and not the random variable, see \cite{Ro19}. 
In such a case $q_a$ is a random variable representing the state of our knowledge about the mean 
of posterior of $q$. On the other hand, 
if the measurement data $y_m$ are stochastic, then the estimate in
 \refeq{eqn:expansoin_update} is stochastic as well, and represents
 the mapped aleatoric uncertainty from $y_m$ to $q$ plus the remainings of a priori knowledge. 
 If $y_m$ and its polynomial chaos expansion are known, then one may use \refeq{eqn:expansoin_update}
 as a computationally cheap trick to map $y_m$ to $q$. In a mutli-scale analysis this is often the case as
 $y_m$ are mostly obtained 
 by virtual computational simulations. However, in general (e.g.~when performing the real experiment) we do not have the continuous data $y_m$ (in a form of a random variable), but its discrete version 
given as a set of random meso-scale realisations 
observed via 
\begin{equation}\label{eq_coll_data}
 y_m(\omega_i) = Y_m(z_m(\kappa_m(\xi_\kappa(\omega_i)))+\hat{\epsilon}(\xi_\epsilon(\omega_i)), \quad i=1,...,n
\end{equation}
such that the measurement data set reads $y_d:=(y_m(\omega_1),...,y_m(\omega_n))^T$. Furthermore, we assume that this is also 
the case
when $y_m$ is obtained by virtual simulation coming from a high-dimensional problem 
as further described in numerical examples. Namely, in such a case the straightforward 
 estimation of the polynomial chaos expansion of $y_m$ can be computationally expensive. 
 Therefore, we distinguish two type of problems in the following text:
 \begin{prob}
 find stochastic $q$ in a form of polynomial chaos expansion such that the macro-scale predictions $y_M$ match the meso-scale ones $y_m$ given in a form of a polynomial chaos expansion 
 \begin{equation}\label{meso_pce}
 \vek{y}_m (\vek{\xi}(\omega)) \approx \sum_{\alpha \in \mathcal{J}_m} \vek{y}_m^{(\alpha)}\Gamma_\alpha(\vek{\xi}(\omega))
 \end{equation} 
 in which $\mathcal{J}_m$ is the multi-index set, $\Gamma_\alpha(\vek{\xi}(\omega))$ is the set of orthogonal polynomials with random variables 
 $\vek{\xi}(\omega):=(\vek{\xi}_\kappa,\vek{\xi}_\epsilon)$ as arguments. Note that $\vek{\xi}(\omega)$ describe next to the natural variablity 
the modelling (e.g. discretisation) error. 
 \end{prob}
 In such a case the posterior $q$ following \refeq{eqn:expansoin_update} reads:
\begin{equation}\label{hit_eq1}
\sum_{\alpha \in \mathcal{J}_a} \vek{q}_a^{(\alpha)}H_\alpha(\vek{\theta}(\omega),\vek{\xi}(\omega))=\sum_{\alpha \in \mathcal{J}}
\vek{q}_f^{(\alpha)}\varPsi_\alpha(\vek{\theta}(\omega))+
\varphi_n( \sum_{\alpha \in \mathcal{J}_m} \vek{y}_m^{(\alpha)}\Gamma_\alpha(\vek{\xi}(\omega)))-\varphi_n(\sum_{\alpha \in \mathcal{J}} \vek{y}_M^{(\alpha)}\varPsi_\alpha(\vek{\theta}(\omega)))
\end{equation}
in which $H_\alpha$ is the generalised polynomial chaos expansion with random variables $(\vek{\theta}(\omega),\vek{\xi}(\omega))$ as arguments. As $\vek{\theta}(\omega)$ describe the a priori uncertainty, one may take mathematical expectation of the previous equation w.r.t. $\vek{\theta}(\omega)$ to obtain the natural variability:
\begin{eqnarray}
\sum_{\alpha \in \mathcal{J}_m} \vek{q}_a^{(\alpha)}\Gamma_\alpha(\vek{\xi}(\omega)) &=&\mathbb{E}_{\vek{\theta}}\left(\sum_{\alpha \in \mathcal{J}_a} \vek{q}_a^{(\alpha)}H_\alpha(\vek{\theta}(\omega),\vek{\xi}(\omega))\right)=\varphi_n( \sum_{\alpha \in \mathcal{J}_m} \vek{y}_m^{(\alpha)}\Gamma_\alpha(\vek{\xi}(\omega)))+\nonumber\\
&& \mathbb{E}_{\vek{\theta}}\left(\sum_{\alpha \in \mathcal{J}}
\vek{q}_f^{(\alpha)}\varPsi_\alpha(\vek{\theta}(\omega))-\varphi_n(\sum_{\alpha \in \mathcal{J}} \vek{y}_M^{(\alpha)}\varPsi_\alpha(\vek{\theta}(\omega)))\right).
\end{eqnarray}
However, the previous estimate is only possible in virtual simulations. Otherwise, the problem generalises to
  \begin{prob}
 find stochastic $q$ in a form of polynomial chaos expansion such that the macro-scale predictions $y_M$ match the meso-scale ones $y_m$ given
 as a set of random meso-scale realisations 
$y_m(\omega_i) = Y_m(z_m(\kappa_m(\xi_\kappa(\omega_i)))+\hat{\epsilon}(\xi_\epsilon(\omega_i)), \quad i=1,...,n$.
 \end{prob}
 The previous problem can be considered from two different perspectives: using \refeq{eqn:expansoin_update} one may estimate $q(\omega_i)$
 for each given $y_m(\omega_i) $, and then estimate its polynomial chaos approximation, or one may first estimate approximation of $y_m$ and then use 
 \refeq{hit_eq1} to upscale the material parameters. Unfortunately, both of these approaches are unsupervised, and hence hard to solve. 
 
 In the first case scenario one may estimate $q_a(\omega_i), \quad i=1,...,n$ by repeating update formula in \refeq{hit_eq} $n$ times:
\begin{equation}\label{hit_eqre}
\forall \omega_i:\vek{q}_a(\omega_i):= \sum_{\alpha \in \mathcal{J}} \vek{q}_i^{(\alpha)}\varPsi_\alpha(\vek{\theta})=\sum_{\alpha \in \mathcal{J}}
\vek{q}_f^{(\alpha)}\varPsi_\alpha(\vek{\theta})+
\varphi_n(\vek{y}_m(\omega_i))-\varphi_n(\sum_{\alpha \in \mathcal{J}} \vek{y}_M^{(\alpha)}\varPsi_\alpha(\vek{\theta})).
\end{equation}
By avaraging over a priori uncertainty, we may further define a set of samples:
 \begin{equation}\label{hit_eqreav}
\forall \omega_i: \bar{\vek{q}}_i=\mathbb{E}_{\vek{\theta}}(\vek{q}_a(\omega_i)),\quad i=1,...,n,
\end{equation}
i.e. the data which are to be used for the estimation of the conditional distribution of $q$ in a parameteric form. To achieve this, we search for an approximation of $q$ given the incomplete data set $q_d:=(\bar{\vek{q}}_i)_{i=1}^n$, here embodied as a nonlinear mapping of some basic standard random variable $\vek{\eta}$ such as Gaussian or uniform, i.e.
\begin{equation}
 \varphi_q:\vek{ \eta} (\omega_i) \rightarrow \bar{\vek{q}}_i
\end{equation}
such that
\begin{equation}\label{measurement}
 \bar{\vek{q}}_i(\omega_i)=\varphi_q(\vek{w}_q,\vek{\eta}(\omega_i))
\end{equation}
holds. Note that then both the mapping $\varphi_q$ and the random variable $\vek{\eta}(\omega)$ (e.g. Gaussian) are unknown, and 
have to be found. 

Similar discussution holds for a given set of measurements $y_d:=(\vek{y}_m(\omega_i))_{i=1}^n$.
  Generally speaking, one may model $y_m$
as a nonlinear mapping of some basic standard random variable $\vek{\zeta}$ 
\begin{equation}
 \varphi_y: \vek{\zeta}(\omega_i)  \rightarrow y_m \vek{y}_m(\omega_i)
\end{equation}
such that
\begin{equation}\label{measurement1}
 \vek{y}_m(\omega_i)=\varphi_y(\vek{w}_y,\vek{\zeta}(\omega_i))
\end{equation}
holds similarly to \refeq{measurement}. As the mapping in both \refeq{measurement} and \refeq{measurement1} is not unique, one has that $\vek{\xi}$ is not same as $\vek{\eta}$ or $\vek{\zeta}$. Furthermore,  having $\vek{\eta},\vek{\zeta}$ as well as parameters $\vek{w}_q$ and $\vek{w}_y$ as unknowns, both \refeq{measurement} and \refeq{measurement1} outcome in an undetermined set of equations as the number of unknowns grows with the 
data set, and its always larger than its size. Therefore, the appropriate regularisation has to be applied. This is further explained on the example of the measurement set, but similar holds if $q_d$ is to be approximated.

\subsection{Identification of the meso-scale observation}

Let the measurement be approximated by
%
\begin{equation}
 y_m=\varphi_y(w,\eta)
\end{equation}
in which $\varphi_y$ is an analytical function (e.g. Gaussian mixture model, neural network, etc.) parameterised by global variables/parameters $w$
describing the whole data set,
and the latent local/hidden variables $\eta$
that describe each data point. An example is the generalised mixture model in which parameters $w$ include statistics of 
individual components, and the mixture weights, whereas the hidden variable $\eta$ stands for the indicator variable 
that describes the membership of data points to the mixture components.
The goal is to estimate
the pair ${\beta}:=(w,\eta)$ given data $y_d$ with the help of Bayes's rule, e.g.
\begin{equation}\label{y_m_via_bayes}
 p(\beta|y_d)=\frac{p(y_d,\beta)}{\int p(y_d,\beta) d\beta}.
\end{equation}
Following theory in Section \ref{BGMK}, \refeq{y_m_via_bayes} is reformulated to the computationally 
simpler variational inference problem. In other words, we 
introduce a family of density functions $\mathcal{D}:=\{g(\beta):=g(\beta|\lambda,\varpi)\}$
over $\beta$ indexed by a set of free parameters $(\varpi,\lambda)$ 
that approximate the posterior density $p(\beta|y_d)$, and further optimise the variational parameter
values by minimising the Kullback-Leibler divergence between the approximation $g(\beta)$ 
and the exact posterior $p(\beta|y_d)$. 
In other words, following \refeq{min_elbo} we maximise the ELBO 
\begin{eqnarray}\label{first_fun}
 \mathcal{L}(g)=\mathbb{E}_{g(\beta)}(\textrm{log }p(y_d,\beta))
 -\mathbb{E}_{g(\beta)}(\textrm{log }g(\beta))))
 \end{eqnarray}
 by using the mean-field factorisation assumption, and conjugacy relationships. The optimisation problem 
attiains closed form solution in which the lower bound is iteratively
optimised with respect to the global parameters keeping the local parameters fixed, and in
the second step the local parameters are updated and the global parameters are hold fixed. 
The algorithm can be imporved by 
considering the stochastic optimisation in which the noisy estimate of the gradient is used instead of the natural one. 
In a special case the previous algorithm reduces to the expectation-maximisation (EM) one. Assuming that the variational density 
matches the posteror one, the first term in \refeq{first_fun} vanishes such that the log-evidence is equal to the ELBO. 
Then, fixing the parameter $\lambda$, i.e. the global variables,
and taking $p_\lambda(\lambda|y_d)=\delta(\lambda-\lambda^*|y_d)$, $g_\lambda(\lambda)=\delta(\lambda-\lambda^*)$ one may rewrite \refeq{first_fun} into 
\begin{eqnarray}
 \textrm{log }p(y_d)=\mathcal{L}(g)=\mathbb{E}_{g(\beta)}(\textrm{log }p(y_d,\beta))
 -\mathbb{E}_{g(\beta)}(\textrm{log }g(\beta))))
\end{eqnarray}
 i.e.
\begin{eqnarray}
 \textrm{log }p(y_d|\lambda^*)&=&\mathcal{L}(g_\varpi(\varpi),\lambda^*)\nonumber\\
 &=&\mathbb{E}_{g_\eta}(\textrm{log }p(y_d,\lambda^*,\varpi)) 
 -\mathbb{E}_{g_\varpi}(\textrm{log }g_\varpi(\varpi)))\nonumber
 \end{eqnarray}
 Here, $\mathbb{E}_{g_\varpi}(\textrm{log }g_\varpi(\varpi)))$ is a constant and represents the entropy of $\varpi$ given $y_d$.
Hence, it is not 
 taken into consideration when maximizing the ELBO.  
 To maximise the previous functional, one may use the iterative scheme consisting of expectation and maximisation  steps which are altered
 as further described in \cite{Hoffman}.

Note that the presented approach is similar to the one 
described in Section \ref{BGMK}. The only
difference
is that this time we do not have full set of data $(\xi(\omega_i),y_m(\omega_i))$ but its incomplete
version generated only by $y_m(\omega_i)$.
Hence, the estimation is more general
and includes the problem described in Section \ref{BGMK} as a special case. 

\subsection{Dependence estimation}

The mean field factorisation as presented in the previous section is computationally simple, but not accurate. For example, one cannot assume independency
between the stored and dissipation energies coming from the same experiment. 
In other words, the correlation among the 
latent variables is not explored. As a result, 
the covariance of the measurement will be underestimated. To introduce the dependency 
into the factorisation, one may extend the mean-field approach
via copula factorisations \cite{tran,Tran2015CopulaVI}:
\begin{equation}
\label{copula_factorisation}
g(\beta)=c(\varPhi(\beta_1),...,\varPhi(\beta_m),\chi)\prod_{i=1}^m g(\beta_i)
\end{equation}
in which $c(\varPhi(\beta_1),...,\varPhi(\beta_m),\chi)$ is the representative of the copula family, 
$\varPhi(\beta_i)$ is the marginal cumulative distribution function of the random variable $\beta_i$ 
 and $\chi$ is the 
set of parameters describing the copula family. Similarly, $g(\beta_i)$
represent the independent marginal densities. In this manner any distribution type can be 
represented by a formulation as given in \refeq{copula_factorisation} according to the Sklar's theorem\cite{sklar1959}. 

Following \refeq{copula_factorisation}, the goal is to find $g(\beta)$ such that the Kullback-Leibler divergence 
to the exact posterior distribution is minimised. Note that if the true posterior is described by 
\begin{equation}
p(\beta|y_d)=c_t(\varPhi(\beta_1),...,\varPhi(\beta_m),\chi_t) \prod_{i=1}^m f(\beta_i)
\end{equation}
then the Kullback-Leibler divergence reads: 
\begin{equation}
D_{KL}(g(\beta)|| p(\beta|y_m)=D_{KL}(c||c_t)+\sum_{i=1}^m D_{KL}(g(\beta_i)||f(\beta_i))
\end{equation}
and contains one additional term than the one characterising the mean field approximation.
When the copula is uniform, the previous equation reduces to the mean field one, and hence only the second term is minimised. On the other hand, if the mean field factorisation is not good assumption and the dependence 
relations are neglected, then the total approximation error will dominantly be specified by the first term. To avoid this, the ELBO derived in \refeq{first_fun}
modifies to 
\begin{equation}
\mathcal{L}(g)=\mathbb{E}_{g(\beta)}(\textrm{log } p(y_m,\beta))- \textrm{log } g(\beta,\chi)
\end{equation}
and is a function of parameters of the latent variables $\beta$, as well as of the copula parameters $\chi$. Therefore, the algorithm applied here 
consists of iteratively finding the parameters of the mean field approximation, as well as those of the copula. The algorithm 
is adopted from \cite{Tran2015CopulaVI} and is a black-box algorithm as it only depends on
the likelihood $p(y_m,\beta)$ and copula descripton in a vine form. Note that when the copula is equal to identity, i.e. uniform, 
the previous factorisation collapses to the mean field one.

%% file: section4.tex
\section{Energy based Bayesian upscaling}

An example of abstract model in Section 1 represents the coupled elasto-damage model 
introduced in \cite{ibrahimbegovic} in which
the state variable $z\in \mathcal{Z}_M=\mathcal{U}\times \mathcal{W}$ consists of the displacement
$u\in \mathcal{U}:=H_\Gamma^1(\mathcal{G};\mathbb{R}^d)=
\{u \in H^1(\mathcal{G};\mathbb{R}^d) \textrm{ }|\textrm{ } u_{|\Gamma}=0\}$, and the generalised damage strain tensor 
$\vek{E}_d:=(\varepsilon_d, \vek{\nu})\in \mathcal{W}$. Here, the generalised damage strain $\vek{E}_d\in \mathcal{W}=
L_2(\mathcal{G};\mathbb{R}^{d\times d}_{\textrm{vol}})\times L_2(\mathcal{G},
L_2(\mathcal{G};\mathbb{R}^{d\times d}_{\textrm{vol}}\times \mathbb{R})$ accounts for the damage strain 
$\varepsilon_d\in L_2(\mathcal{G};\mathbb{R}_{\textrm{vol}}^{d\times d})$ with
$\mathbb{R}^{d\times d}_{\textrm{vol}}=\{ \varepsilon_d \in \mathbb{R}^{d\times d}\textrm{ }|\textrm{ }
\varepsilon_d=\varepsilon_d^T, \textrm{ dev}(\varepsilon_d)=0\}$, 
and internal variables $(\varepsilon_d,\varsigma)\in L_2(\mathcal{G};\mathbb{R}^{d\times d}_{\textrm{vol}}\times \mathbb{R})$ . 
The total energy $\mathcal{E}_M$
then reads:
\begin{equation}
 \mathcal{E}(t,u,\varepsilon_p,\varepsilon_d)=\int_\mathcal{G}\left(\frac{1}{2}\vek{C}\varepsilon_e(u):\varepsilon_e(u) 
+\frac{1}{2} \varepsilon_d D^{-1}\varepsilon_d+\frac{1}{2}\vek{H}_d \vek{E}_d:\vek{E}_d 
 -f_{ext}(t) u \right) dx
\end{equation}
in which $\vek{C}$, $\vek{D}$, $\vek{H}_d$ denote positive-definite elasticity and damage compliance and hardening tensors, respectively. 
In further numerical experiments the previous model is used on both meso- and macro-scales, specified by the spatial avarage of the
complementary total energy
\begin{equation}
\psi_M = \int_{\mathcal{G}}\frac{1}{2}(\sigma\varepsilon_e  + \sigma\varepsilon_d + \chi^{d}\varsigma) dx
\end{equation}
 and dissipated potential
\begin{equation}
\varphi_M = \int_{\mathcal{G}}(\frac{1}{2}\sigma\dot{D}\sigma+\chi^{d}\dot{\varsigma})dx 
\end{equation}
 in the domain --- one quadrilateral element --- of
the coarse-scale model. Here, $\chi_d$ is a conjugate hardening force, whereas $\kappa$ is a collection of parameters specifing the
detailed character of the functions $\psi_M$ and $\varphi_M$, i.e. $\kappa:=\{\vek{C},\vek{D},\vek{H}_d\}$. Note that 
the failure criterion is specified as
\begin{equation}
\label{const3}
f_{d}(\sigma,\chi_{d})=  \langle -\text{tr}(\sigma) \rangle - (\sigma_{f}-\chi_{d}) ,
\end{equation}
and depends on the parameter $\sigma_f$, the limit at which the damage occurs. 
 As we only focus 
on the isotropic case, one may recognise that $\kappa=\{K,G,\sigma_f,K_d\}$ with $K$ and $G$ being the bulk and shear moduli, respectively, 
and $K_d$ being the isotropic damage hardening.

Finally, the upscaling is considered for the energy-type of measurement $y_m:=( \mathcal{E}_e,\mathcal{E}_d, \mathcal{E}_h)$
in which 
\begin{equation}
\label{measur2}
\mathcal{E}_e:=\int_{\mathcal{G}}\frac{1}{2}\sigma\varepsilon_e dx, ,\quad \mathcal{E}_d:=\int_{\mathcal{G}} \frac{1}{2}\sigma\varepsilon_d dx, \quad \mathcal{E}_h:=\int_{\mathcal{G}} \frac{1}{2} \chi^{d}\varsigma dx,
\end{equation}
respectively. In this manner,  the non-dissipated as well as dissipated portion of energies are conserved when moving
 from the meso- to the macro-scale model. 
 

To represent the measurement data $y_m$ one may use the generalised mixture models. In our praticular application the measurement is positive definite. Therefore,
we use samples $x:=(\textrm{log }y_m(\omega_i))_{i=1}^N$ to approximate
$\textrm{log }y_m$ as a Gaussian mixture model \cite{Bishop}
\begin{equation}
 p(x)=\sum_{k=1}^K \pi_K \mathcal{N}(x|\nu_k), \quad \sum_{k=1}^K \pi_k=1, \quad 0 \leq \pi_k \leq 1
\end{equation}
described by parameter set $\vek{\nu}=(\mu_k,\varSigma_k)_{k=1}^K$ with $\nu_k:=(\mu_k,\varSigma_k)$
being the statistics parameters of Gaussian components, and $\pi_k$ are the mixing coefficients. These generate the 
parameter vector $w$.
The hidden variable $\eta$ is the indicator vector $z_k$ of dimension $N$ that describes the membership of each data point 
to the Gaussian cluster. Following this, the joint distribution is given as
\begin{equation}
 p(x,Z,\mu,\Sigma,\pi)=p(x|z,\mu,\Sigma,\pi)p(z|\pi)p(\pi)p(\mu|\Sigma)p(\Sigma)
\end{equation}
in which
\begin{equation}
 p(z|\pi)=\prod_{n=1}^N \prod_{k=1}^K \pi_k^{z_{nk}}, \quad p(z_k=1)=\pi_k
\end{equation}
and
\begin{equation}
 p(x|z,\mu,\Sigma,\pi)=\prod_{n=1}^N \prod_{k=1}^K \mathcal{N}(x_n|\mu_K,\varSigma_K^{-1})^{z_{nk}}.
\end{equation}
The priors are chosen such that $p(\pi)$ is the Dirichlet prior, whereas $p(\mu,\Sigma)$ follows an independent Gaussian-Wishart 
prior governing the
mean and precision of each Gaussian component. Hence, our parameter set $\beta$ is described by a set of global parameters
$w:=(\mu,\Sigma,\pi)$ and the hidden variable $z$.

To incorporate correlations, the copula dependence structure of Gaussian mixture as in \refeq{copula_factorisation} is found, and the measurement data
are represented in a functional 
form, which is different than the polynomial chaos representation that is required in \refeq{meso_pce}. 
In other words, the measurement is given in terms of dependent random variables, and not independent ones.
Therefore, the dependence structure has to be mapped to the independent one. In a Gaussian copula case, 
the Nataf transformation can be used. Otherwise, the Rosenblat transformation is applied. 
For high-dimensional copulas such as regular vine copula \cite{Aas}  provided  algorithms 
to  compute  the  Rosenblatt  transform and  its
inverse. The result of the transformation are mutually independent and marginally 
uniformly distributed random variables, which further can be mapped to 
Gaussian ones or other types of standard random variables via marginals, \cite{sudret}.

%% file: example1.tex
\subsection{Bayesian upscaling of linear elastic material model}
In this section the proposed upscaling scheme is applied on a random heterogeneous material modelled as linear elastic. The example specimen consists 
of a 2D block described by 64 circular inclusions of equal size randomly distributed in the domain. In the first case scenario only one meso-scale realisation 
is observed for the validation purpose. Furthermore, the stochastic upslacing is considered.

\begin{figure}
\begin{center}
\includegraphics[width=0.7\textwidth]{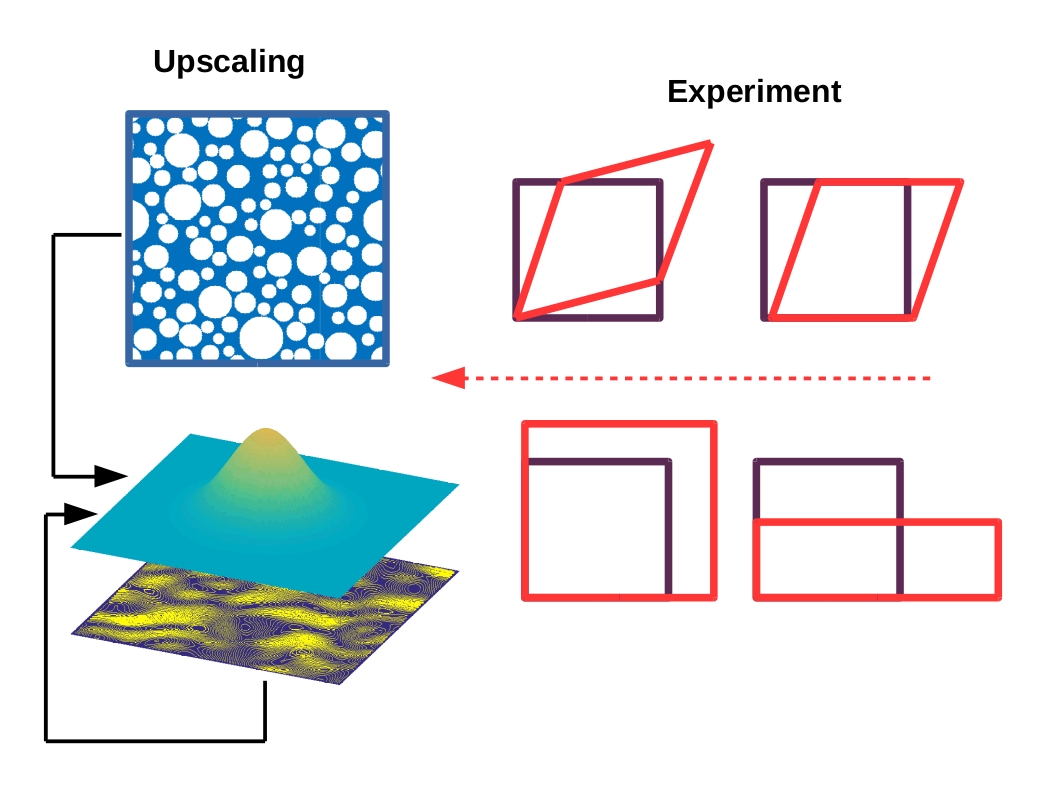}
\caption{Experimental setup}
\label{Figexpseup2d}
\end{center}
\end{figure}

 The meso-scale characteristics are upscaled 
in a Bayesian manner to the coarse scale homogeneous isotropic finite element described by material properties taking the form of a posteriori random variable as schematically shown in Fig.~\ref{Figexpseup2d}. To gather as much as information as possible in observation data,  we consider different types of
loading conditions including shearing and compression solely or their combination as shown in Fig.~\ref{Figexpseup2d}. 
These are further enumerated from the left top to the right bottom as  Exp 1 - Exp 4. Note that in a similar manner one may also use the heterogeneous description 
of the meso-scale model based on the random field realisations of the corrresponding material quantities.

\subsubsection{Validation}

\begin{figure}
\begin{center}
\includegraphics[width=0.39\textwidth]{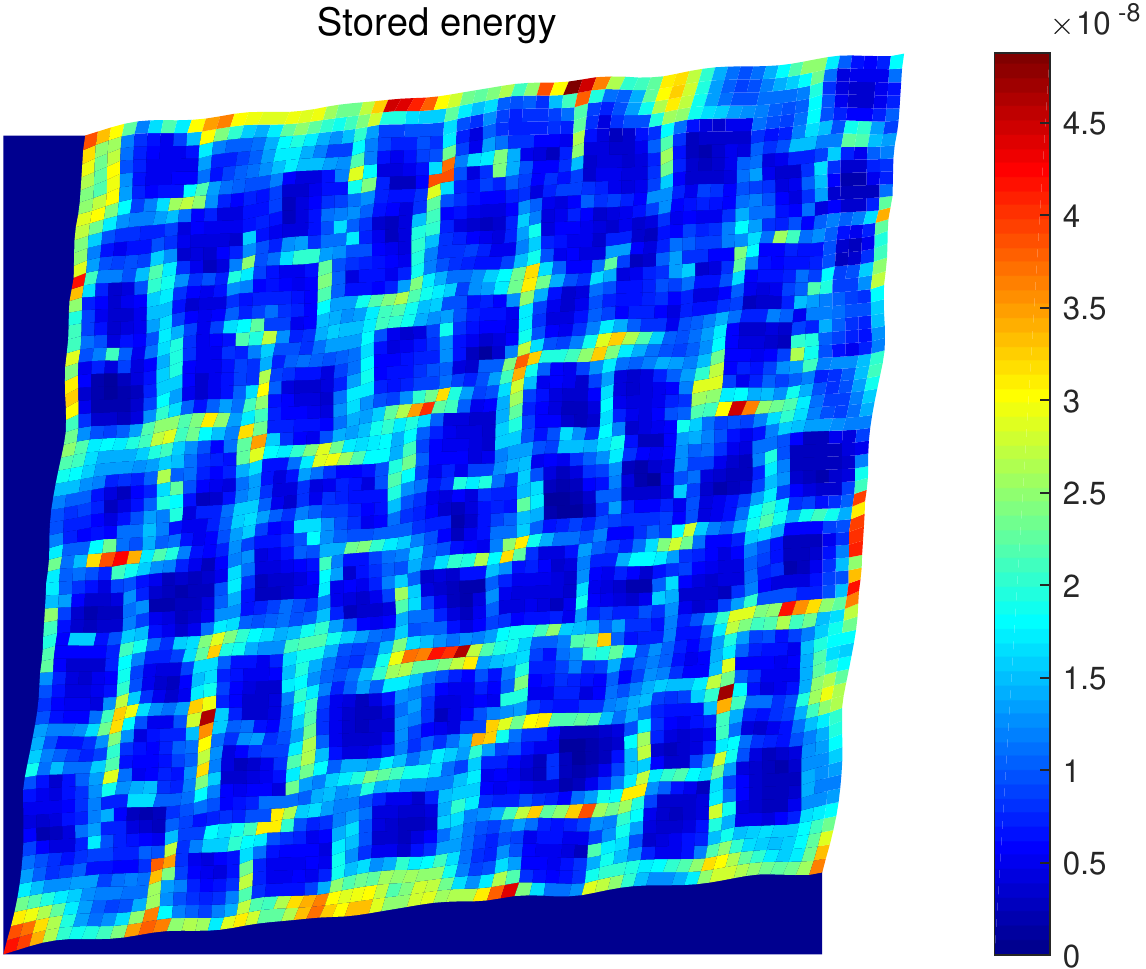}
\end{center}
\caption{Fine scale response under shearing with periodic boundary conditions}
\label{fig:shearing_stored_ener}
\end{figure}

To validate our method, we compare Bayesian upscaling procedure to the deterministic homogenisation approach as presented in \cite{dethom}.
Therefore, we initally observe only one realisation of the random mesostructure and apply periodic boundary conditions.
 An example 
of the fine scale response in terms of stored energy function
is shown in
Fig.~\ref{fig:shearing_stored_ener} for Exp 1. 


\begin{figure}
\includegraphics{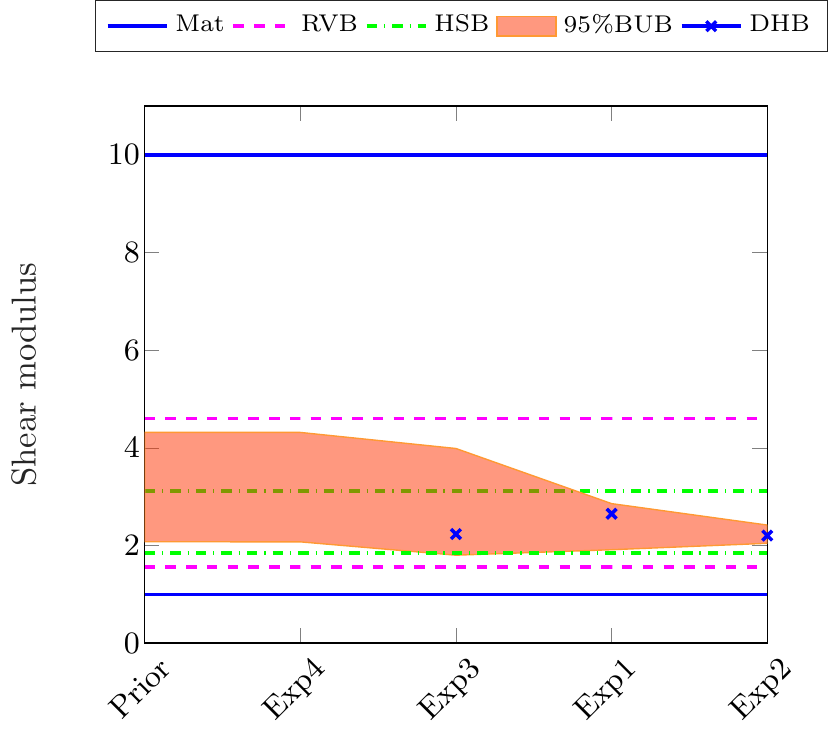}
\includegraphics{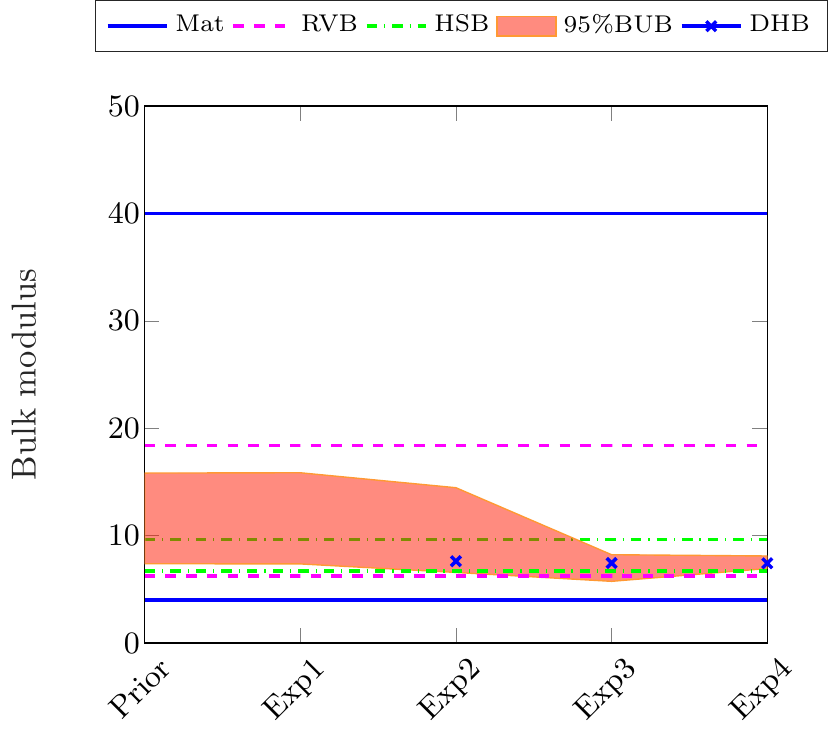}
\begin{center}
{\tiny{RVB$=$Reuss-Voight bounds, HSB$=$Hashin-Strikman bounds, BUB$=$Bayesian updating bounds, DHB$=$classical homogenisation}}
\end{center}
\caption{Upscaling of deterministic material properties: a) shear mnoduli b) bulk moduli}
\label{fig2:shear_bulk}
\end{figure}

As depicted in Fig.~\ref{fig2:shear_bulk}, the deterministic 
homogenisation (DHB) fails to predict the shear moduli in purely compression state (Exp 4). Similar holds for the bulk moduli 
in case when only shear loading conditions are applied (Exp 1). The reason is that in the absence of measurement information on one of the material 
characteristics, the deterministic homogenisation becomes ill-posed.
On the other hand, the Bayesian updating approach (BUB) in a form of \refeq{hit_eq} regularises the problem by introducing the prior information.
This means that the posterior estimate matches the prior when the data are not informative about the parameter set. 
Otherwise, the posterior mode and the deterministic homogenised value are identical. To conclude, the Bayesian upscaling procedure is more robust than the classical one. In addition, one may show that 
the Bayesian upscaling procedure is more appropriate as one may sequentially introduce the measurement data into the upscaling process. For example, 
one may first use the measurement information coming from the fourth experiment to obtain the upscaled material properties. 
These further can be used as a new prior for the third experiment, 
and so on, see \ref{fig2:shear_bulk}.

\subsubsection{Upscaling of random elastic material}

To quantify randomness on the meso-scale level, the previously described experiment is repeated several times, and the avaraged stored elastic energies per experiment are collected. In particular, we observe realisations of
the meso-scale elastic material described by 
randomly placed inclusions depicting the volume fraction of  40\%.  
Initially,  the stored energy is identified given observed data by using 
the variational Bayesian inference method as described in Section 4. The logarithm of energy is modelled by a
copula Gaussian mixture model, and the individual components are identified. 
The optimal number of mixture components is further decoupled by inverse
transform. The resulting uncorrelated random variables are  then further used to obtain the polynomial chaos surrogate of the measurement data. 

\begin{figure}
\begin{center}
\includegraphics{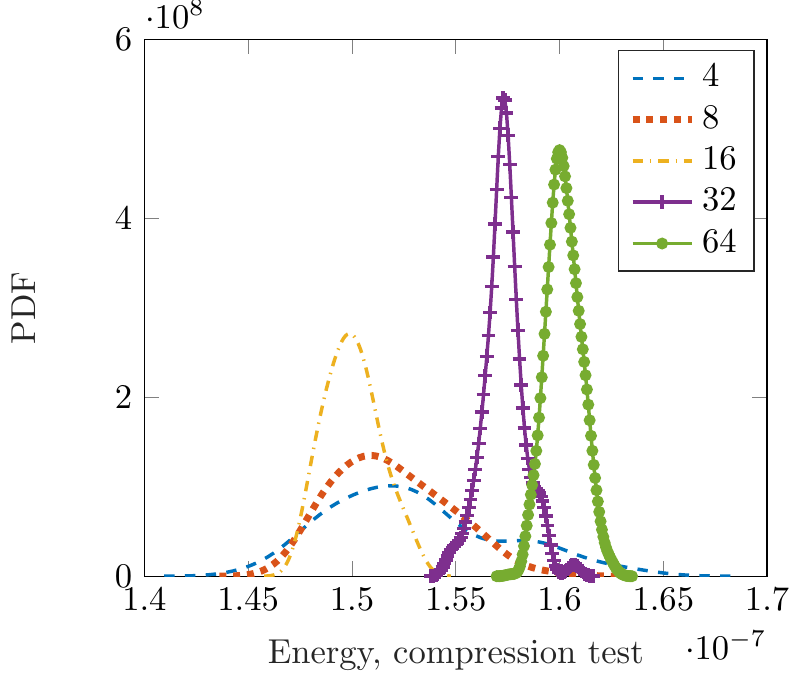}
\includegraphics{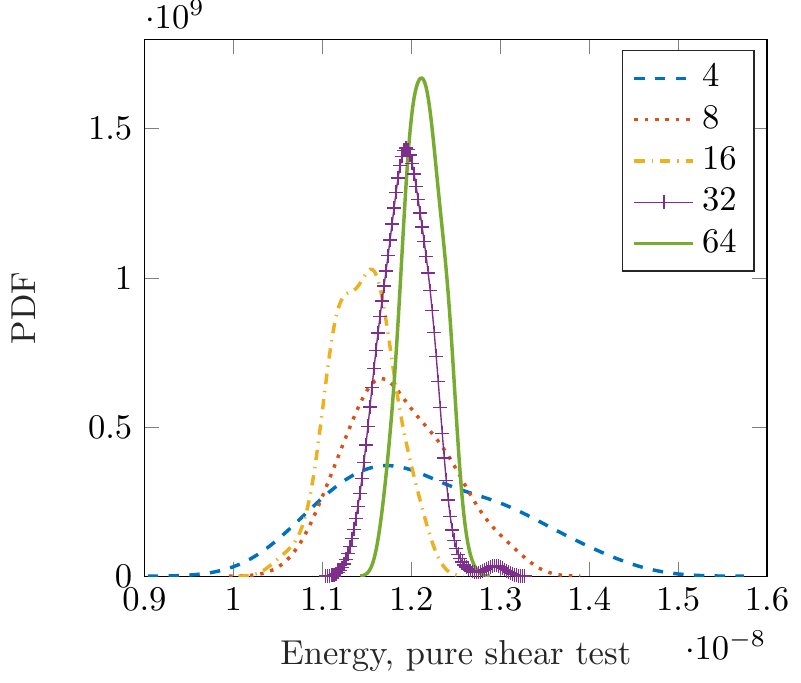}
\caption{Measured energy on 100 mesostructure realisations for different number of inclusions and random position only. The boundary condition is linear displacement.}
\label{Fig:energy_100_micro}
\end{center}
\end{figure}

\begin{figure}
\begin{center}
\includegraphics[width=0.45\textwidth]{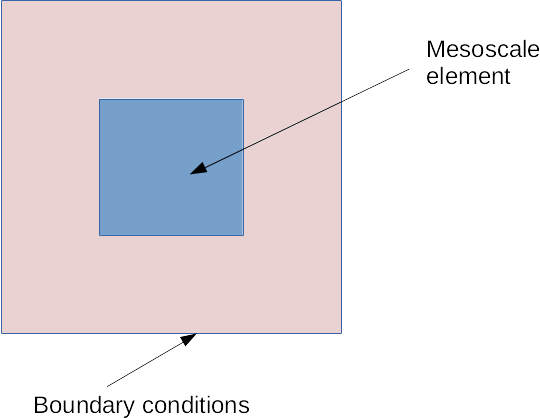}
\caption{Element definition for the upscaling procedure}
\label{figelemupscale}
\end{center}
\end{figure}

The fine-scale simulation is performed on the 2D micro-structure with increasing number of particles
and linear displacement based boundary condition. The material properties of the matrix phase
are taken to be: the bulk moduli is $K_m=4$ and shear moduli is $G_m=1$, 
whereas the inclusions characteristics are prescribed to be ten times higher. 
For a given number of particles embedded in the matrix phase, 
an ensemble of 100 realizations of stored energy is considered to gather corresponding measurement set. 
In Fig.~\ref{Fig:energy_100_micro} the PDFs 
for the identified elastic energies are shown for pure shear and 
bi-axial compression test, respectively. As expected, the variation of stored energy reduces with the increase of the number of particles in the matrix phase. However, it is interesting to note that in compression case the mean responses 
of stored energies vary more than in a pure shear test. This is closely related to the way how boundary conditions are imposed. 
Namely, we take into consideration directly the element on which the loadings are imposed, and thus one may recognise the strong influence of boundary conditions 
on the 
obtained results. In a more proper analysis one would take into consideration only internal element, which is away from the boundary as shown in 
Fig.~\ref{figelemupscale}.

\begin{figure}[htpb]
\begin{center}
\includegraphics{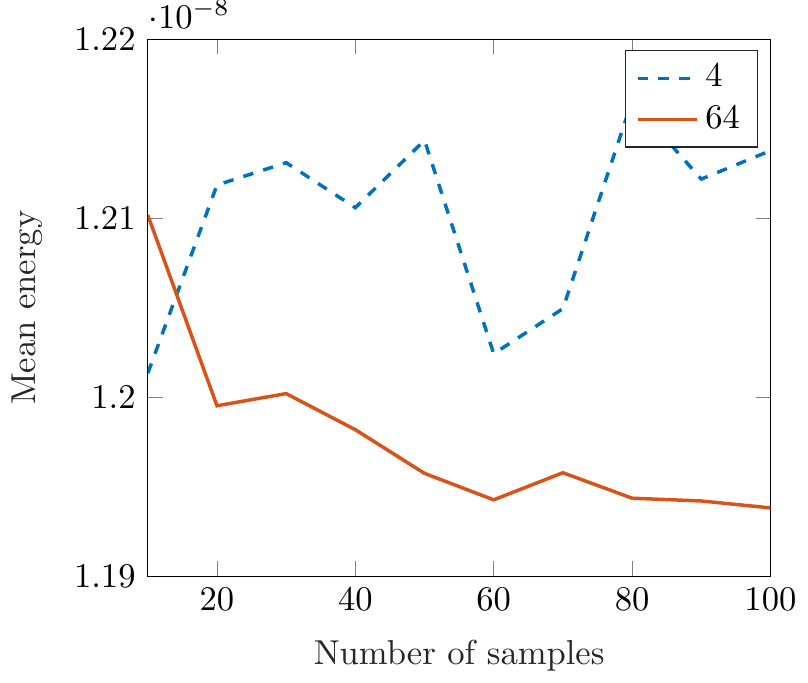}
\includegraphics{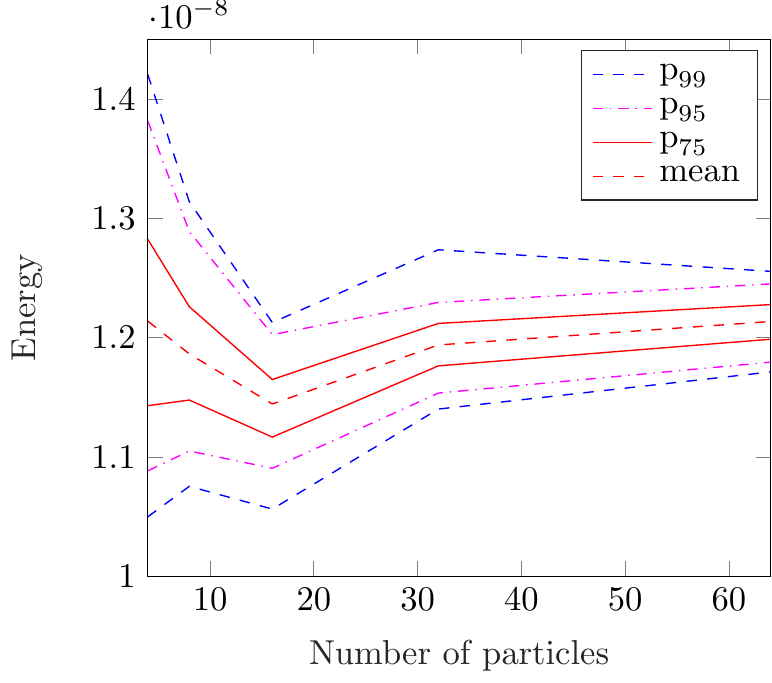}
\caption{Measured energy on 100 mesostructure realisations for different number of inclusions and random position only. The boundary condition is linear displacement. Left: the mean convergence w.r.t. number of Monte Carlo samples, right: Quantiles of energy w.r.t. number of particles.}
\label{fig_pure_shear_2d_elast_rand_pos_energ_no_inclus_g}
\end{center}
\end{figure}

In case of 64 embedded particles, the mean meso-scale energy tends to converge as one increases the size of Monte Carlo  ensemble set. This is in contrast with 
 its 4 particle counterpart, the energy of which keeps on fluctuating, see 
Fig.~\ref{fig_pure_shear_2d_elast_rand_pos_energ_no_inclus_g}a). 
This behaviour is more pronounced in Fig.~\ref{fig_pure_shear_2d_elast_rand_pos_energ_no_inclus_g}b), which depicts the dependence of 
the $p_{50}$, $p_{75}$, $p_{95}$ and $p_{99}$ quantiles on the number of embbeded particles in a matrix phase. The mean energy tends to stabilize to 
a fixed value, whereas the corresponding bounds tend to shrink towards the mean value when the number of particles increases. Hence, the uncertainty becomes less pronounced. 

\begin{figure}
\begin{center}
\includegraphics{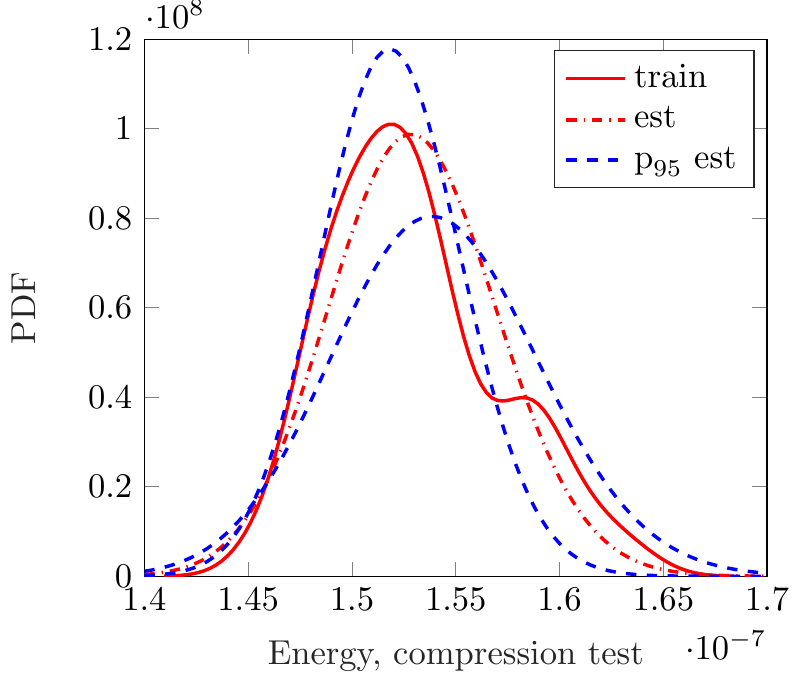}
\includegraphics{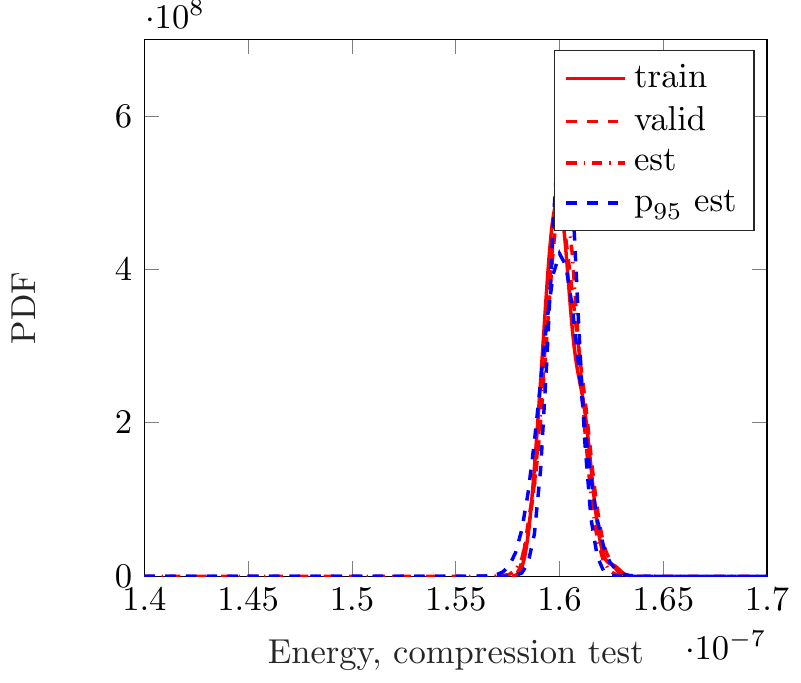}
\caption{ Confidence bounds of estimate a) 4 particles b) 64 particles for the compression experiment}
\label{fig_pure_shear_2d_elast_rand_pos_energ_no_inclus_bb}
\end{center}
\end{figure}

\begin{figure}
\begin{center}
\includegraphics{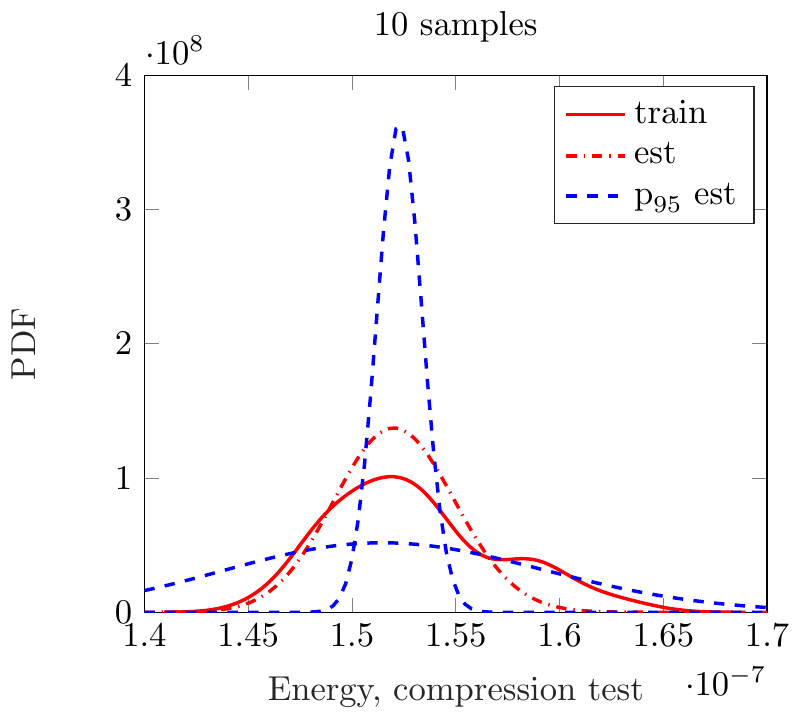}
\includegraphics{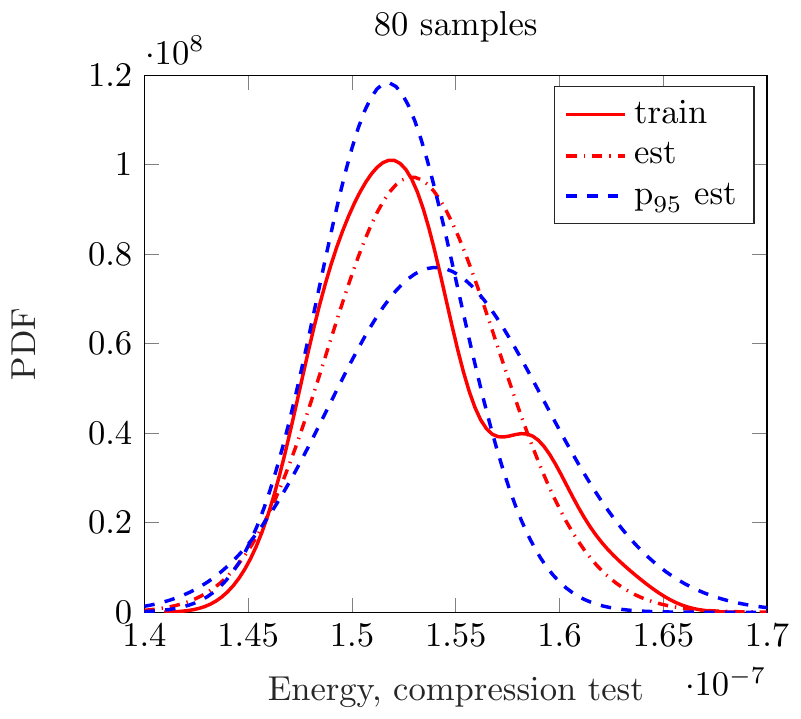}
\caption{Dependence of confidence bounds on the number of samples used in estimation. THe number of particles is 4. The loading is compression}
\label{fig_pure_shear_2d_elast_rand_pos_energ_no_inclus_ba}
\end{center}
\end{figure}

The previously discussed results are conditional estimates of stored energy given its samples. These consist 
of two kinds of uncertainties: aleatoric (meso-scale randomness) and epistemic (prior information). Estimating 
the confidence intervals 
w.r.t. to the epistemic uncertainties one obtains the corresponding PDFs of upscaled material properties: the mean 
PDF which 
represents purely aleatory uncertainty and $p_{95}$ upper and lower PDF's that describe $95\%$ quantiles
on the mean PDF, see Fig.~\ref{fig_pure_shear_2d_elast_rand_pos_energ_no_inclus_bb}. The result indicates that the 
PDF in case of 4 particles 
is slightly underestimated, whereas the PDF in case of 64 particles is overestimated. This phenomenon is often
observed when variational inference 
is used instead of full Bayes's rule. Naturaly the quantile intervals strongly depend on the size of the
measurement set. In 
Fig.~\ref{fig_pure_shear_2d_elast_rand_pos_energ_no_inclus_ba} one may notice that with the smaller 
measurement set (e.g. 10 samples) 
our confidence about the 
upscaled PDF is lower than in case of higher number of measurements (e.g. 80 samples). 

\begin{figure}
\begin{center}
\includegraphics{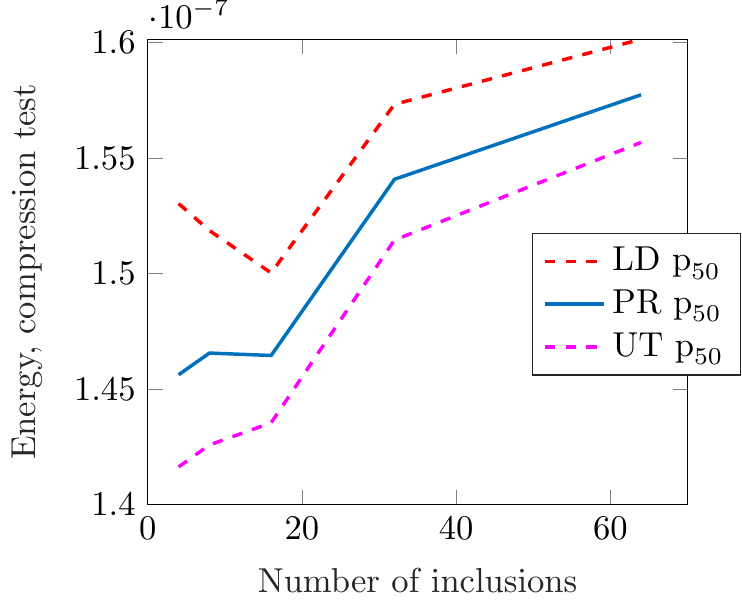}
\includegraphics{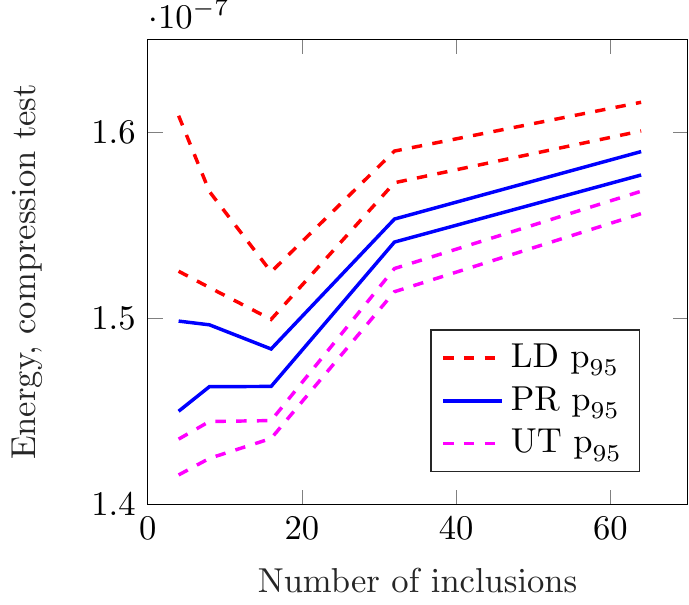}
\caption{Energy wrt boundary conditions and number of inclusions}
\label{fig_compr_diff_loading}
\end{center}
\end{figure}

\begin{figure}
\begin{center}
\includegraphics{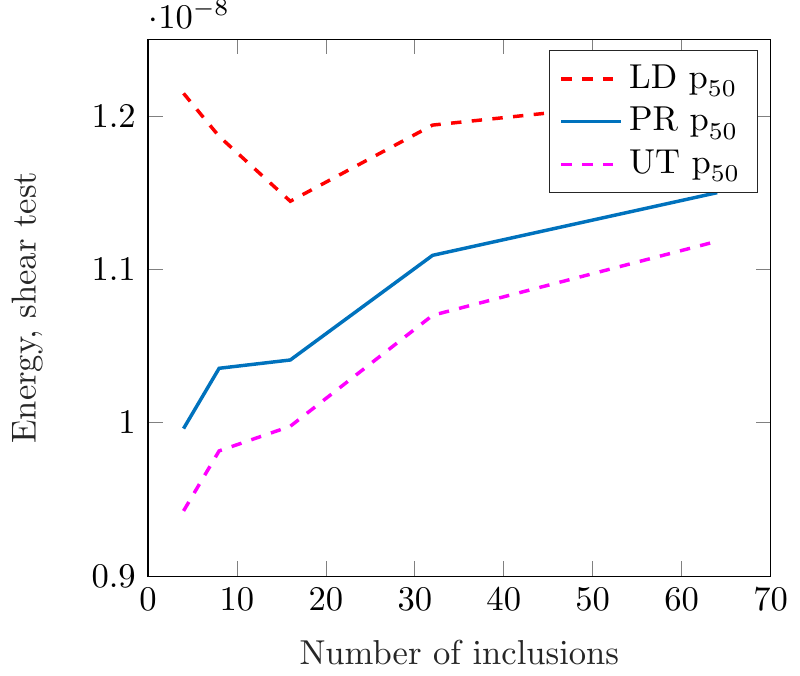}
\includegraphics{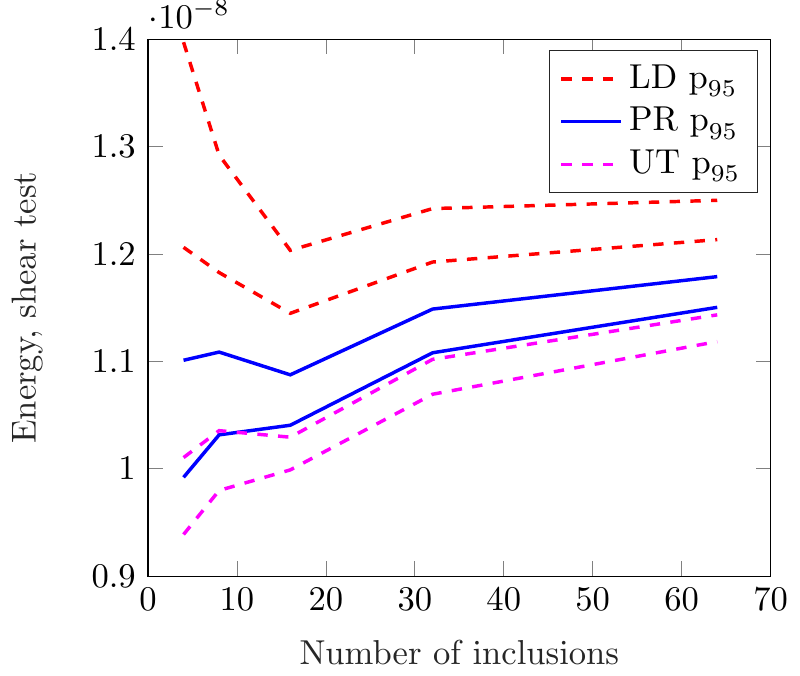}
\caption{Energy wrt boundary conditions and number of inclusions}
\label{fig_shear_diff_loading}
\end{center}
\end{figure}

Besides previous analysis, the impact of boundary conditions on the upscaled quantities 
is another important factor to study. In Fig.~\ref{fig_compr_diff_loading} and Fig.~\ref{fig_shear_diff_loading}
are depicted $50\%$ and $95\%$ quantiles  of energy for linear discplacement (LD), periodic (PR) and unifrom 
tension (UT) boundary conditions. According to these results, linear displacement defines the upper bound on the estimated energy, whereas 
uniform tension gives its lower limit. On the other hand, variations of energies are similar for all three types of boundary 
conditions, and are inverse proportional to the number of inclusions.

\begin{figure}
\begin{center}
\includegraphics{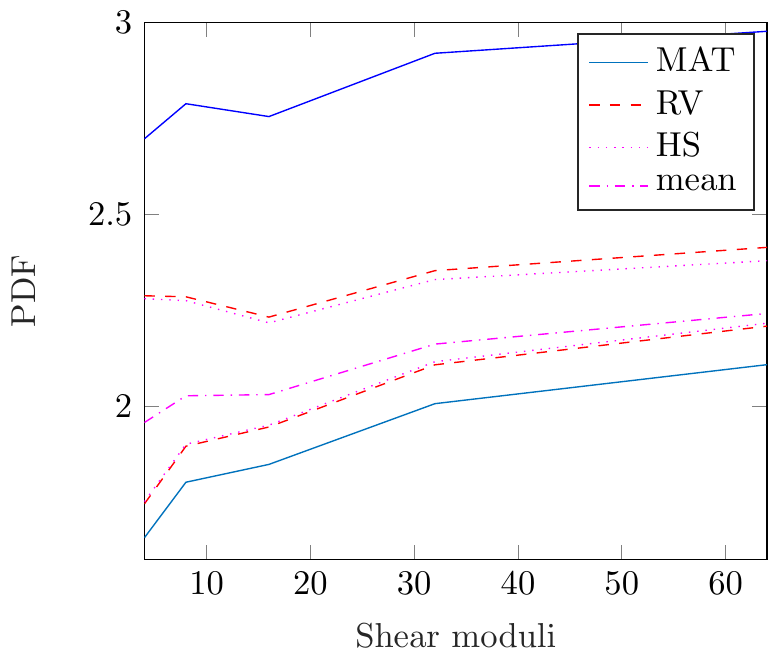}
\includegraphics{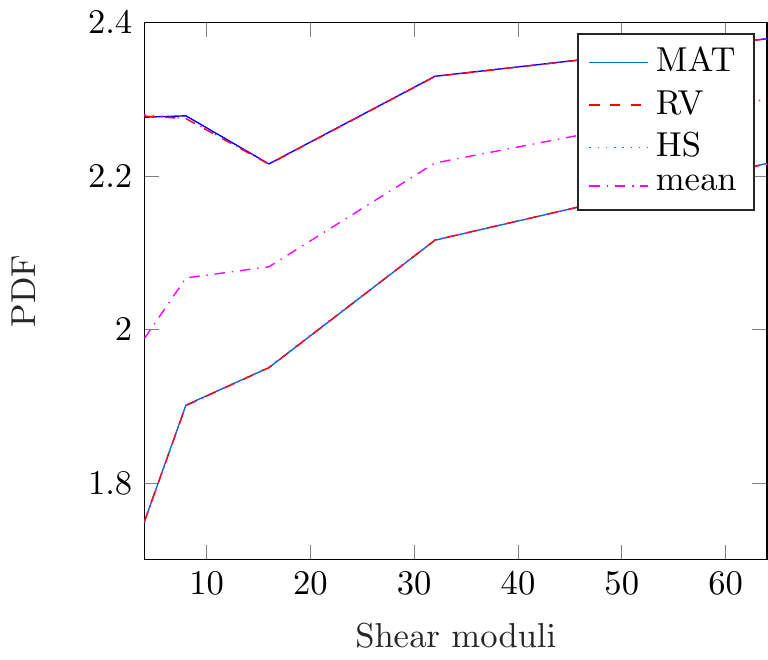}
\caption{Estimated shear moduli w.r.t. prior choice}
\label{Fig:shear_moduli_prior_choice}
\end{center}
\end{figure}

\begin{figure}
\begin{center}
\includegraphics{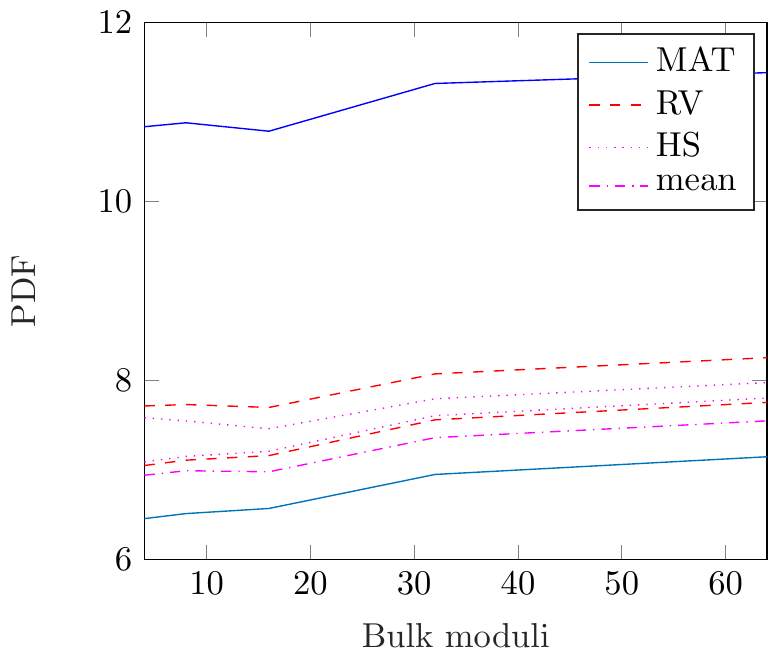}
\includegraphics{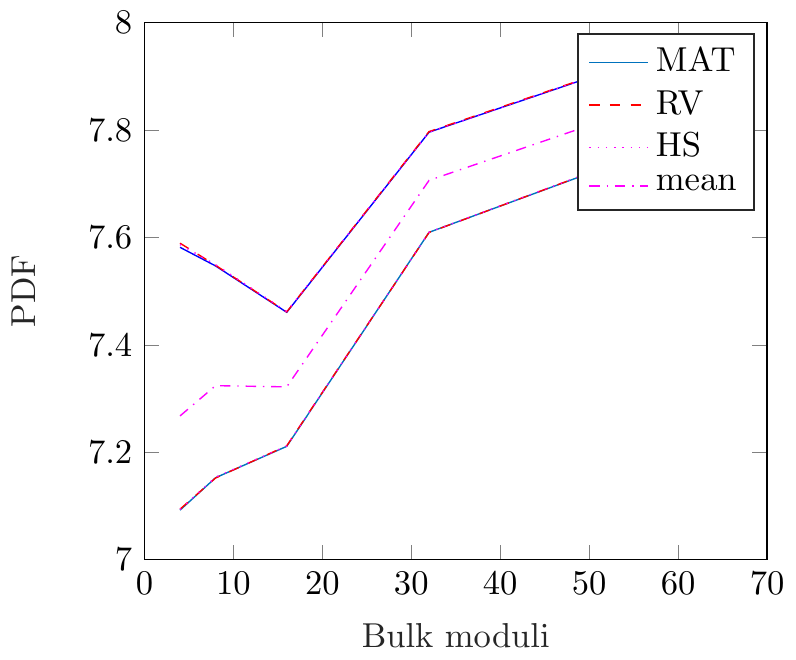}
\caption{Estimated bulk moduli w.r.t. prior choice}
\label{Fig:bulk_diff_prior}
\end{center}
\end{figure}

Once the measurement energy is identified, in the second step we use the proxy of $y_m$ to identify the 
elastic macro-scale material characteristics by using the generalised Kalman filter given in \refeq{hit_eq}.
When using this type of upscaling one is biased to prior knowledge on the material characteristics on the coarse scale. 
In a multiscale analysis, however, it is not an easy task to define the prior knowledge, or better to say
the limits of the prior 
distribution. Therefore, in Fig.~\ref{Fig:shear_moduli_prior_choice} is investigated the posterior change of shear moduli
w.r.t. prior knowledge. 
The prior distributons are chosen such that $95\%$ limits match interval
described by material property of the matrix phase 
and inclusion (in figure denoted by MAT),  Reuss-Voigt (RV)
or Hashin-Shtrikman (HS) bounds. Their corresponding $95\%$ posterior limits are depicted
in Fig.~\ref{Fig:shear_moduli_prior_choice}a). Interesting to note is that even though the posterior of the upscaled shear moduli 
changes w.r.t.~prior assumption, its mean remains the same, see  Fig.~\ref{Fig:shear_moduli_prior_choice}b).
Therefore, there is only one mean PDF in the plot. To better understand this point, the shear moduli update is obtained 
by using the linear Kalman formula
\begin{equation}
 \mu_a(\zeta,\theta)=\mu_f(\theta)+K(y_m(\zeta)-y_M(\theta))
\end{equation}
in which $K$ denotes the Kalman gain. 
Hence, in Fig.~\ref{Fig:shear_moduli_prior_choice}a are shown $95\%$ bounds 
of total posterior $\mu_a(\zeta,\theta)$, whereas 
in Fig.~\ref{Fig:shear_moduli_prior_choice}b are depicted $95\%$ bounds of $\mathbb{E}_\theta(\mu_a(\zeta,\theta))$ and hence only 
aleatory uncertainty. Therefore, all estimates match.  Similar holds for the bulk moduli, see Fig.~\ref{Fig:bulk_diff_prior}.

\begin{figure}
\begin{center}
\includegraphics{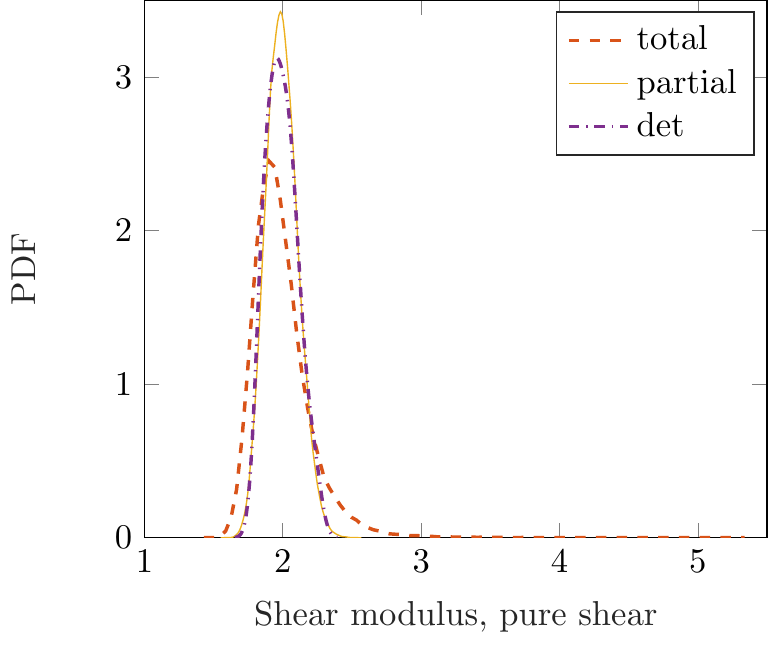}
\caption{Estimated shear moduli: total: $\mu_a(\zeta,\theta)$, partial: $\mathbb{E}_\theta(\mu_a)$, det: deterministic homogenisation}
\label{fig_valid_distrib}
\end{center}
\end{figure}

To validate our result, further in Fig.~\ref{fig_valid_distrib} we compare the mean posterior distribution
$\mathbb{E}_\theta(\mu_a(\zeta,\theta))$
with the posterior distribution obtained by repeating the deterministic homogenisation, see \cite{dethom}, on each of the 
meso-scale realisations. As one may notice, the distributon coming from the deterministic homogenisation and the one obtained by our approach 
are matching. They are further compared with the posterior distribution represented by $\mu_a(\zeta,\theta)$, i.e.
the total uncertainty that includes both aleatory and epistemic knowledge. 

\begin{figure}
\begin{center}
\includegraphics{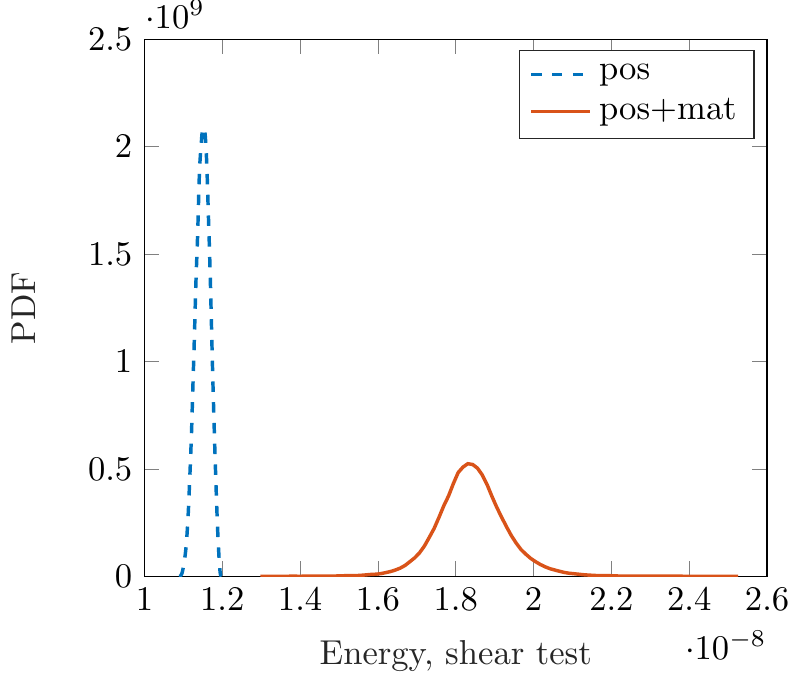}
\includegraphics{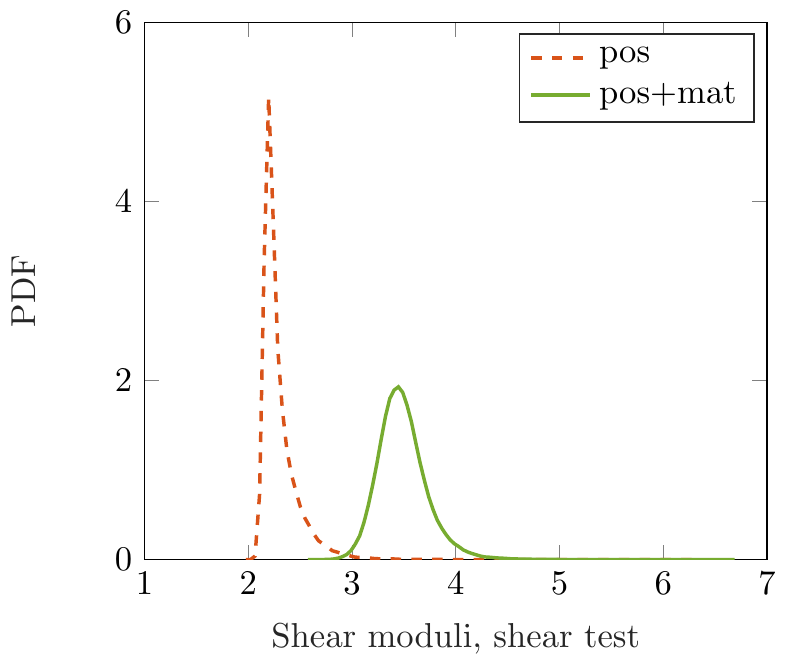}
\caption{Estimated shear energy and shear moduli w.r.t. randomness in position and material properties on the mesoscale}
\label{Fig:posplusmatshear}
\end{center}
\end{figure}

\begin{figure}
\begin{center}
\includegraphics{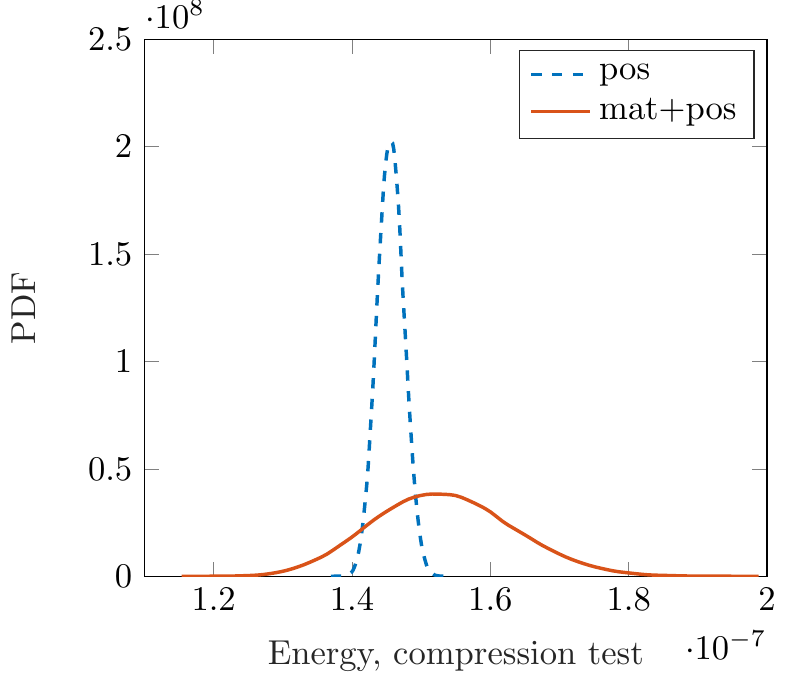}
\includegraphics{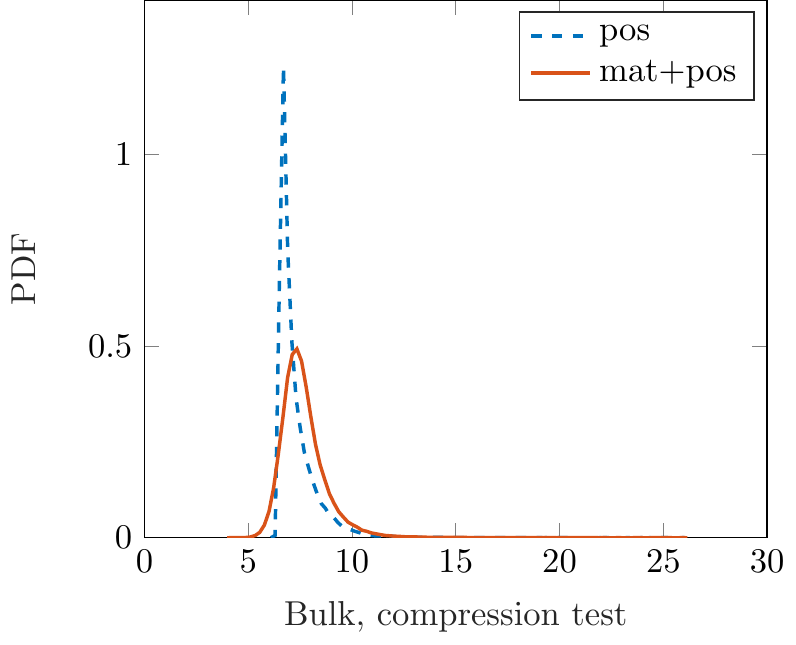}
\caption{Estimated compression energy and bulk moduli w.r.t. randomness in position and material properties on the mesoscale}
\label{Fig:posall}
\end{center}
\end{figure}

\begin{figure}
\begin{center}
\includegraphics{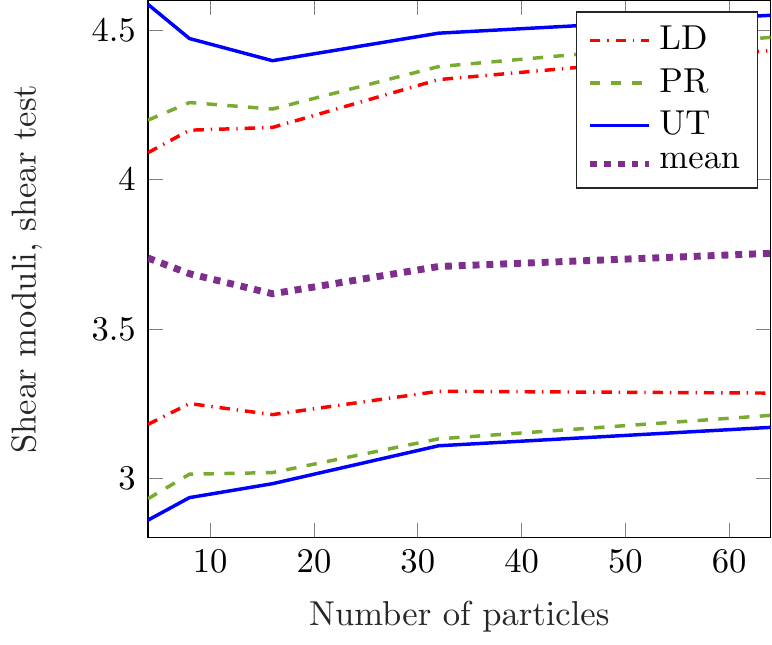}
\includegraphics{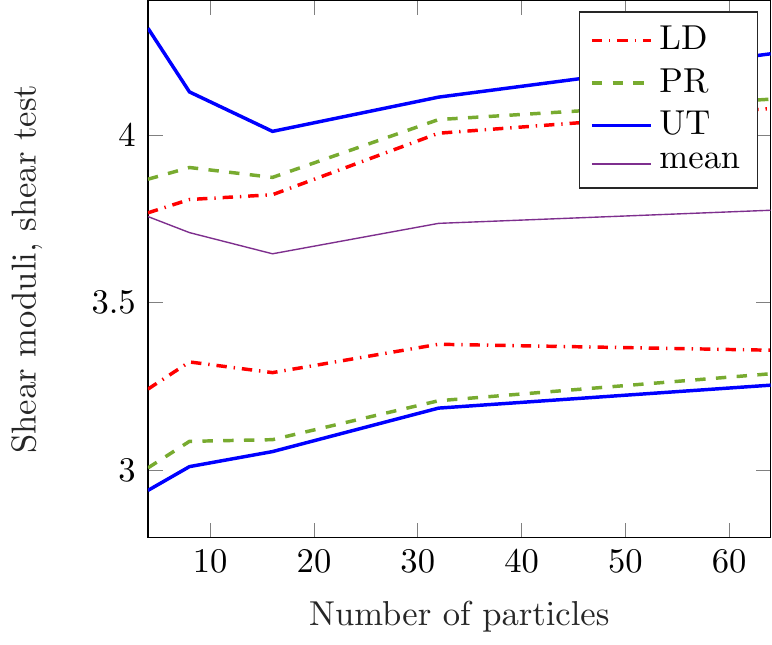}
\caption{Estimated upper and lower 95\% quantiles of shear moduli w.r.t. type of boundary condition. Left  is the total a posteriori uncertainty (aleatory+epistemic), whereas right is only aleatory one. }
\label{Fig:bulk_random_all}
\end{center}
\end{figure}

\begin{figure}
\begin{center}
\includegraphics{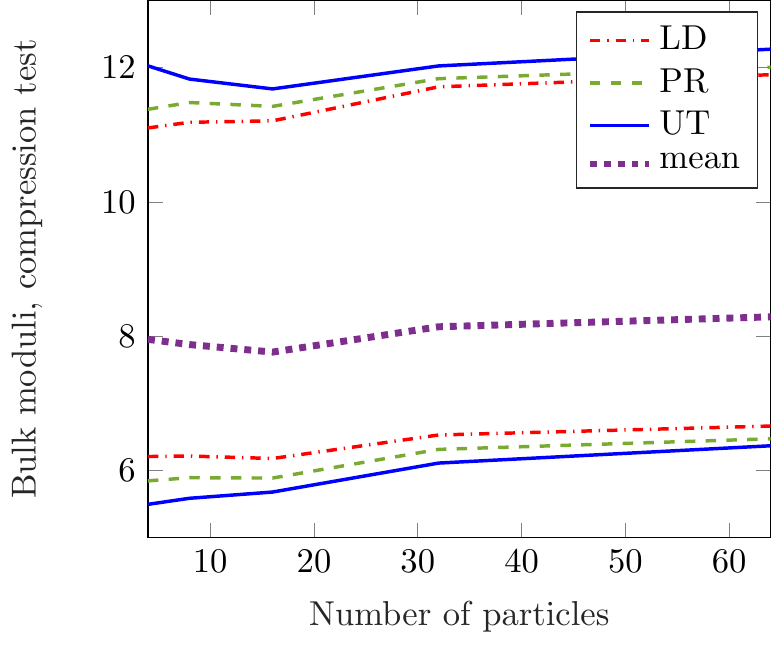}
\includegraphics{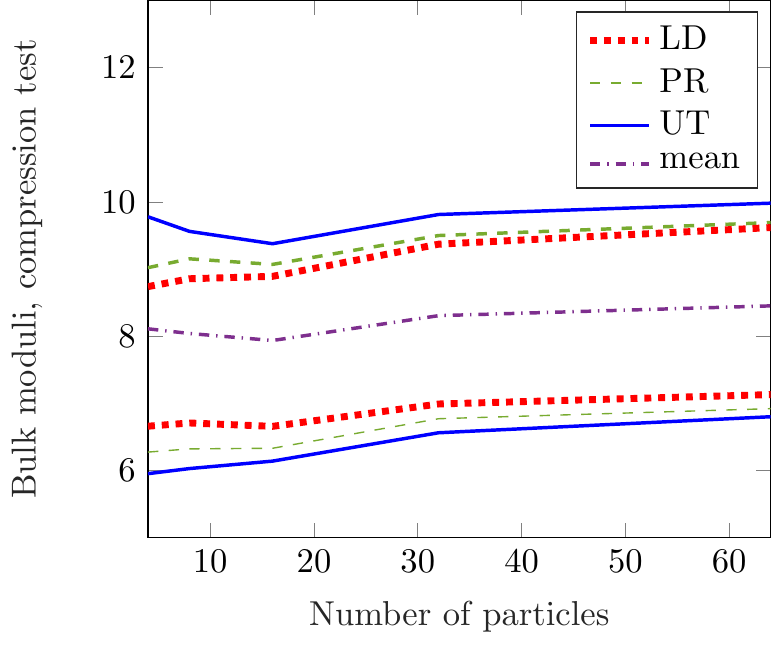}
\caption{Estimated upper and lower 95\% quantiles of bulk moduli w.r.t. type of boundary condition. Left  is the total a posteriori uncertainty (aleatory+epistemic), whereas right is only aleatory one. }
\label{Fig:part_total_random_all}
\end{center}
\end{figure}

So far the fine-scale randomness has been considered 
for the positioning of particles in the matrix phase. In the ensuing discussion, next to the random position
the uncertainty associated with the material properties is also considered. 
For the numerical experiment bi-axial compression or pure shear periodic boundary conditions are chosen. 
The number of particles in the sample are 64 constituting $40\%$ of the total volume.
In Fig.~\ref{Fig:posplusmatshear} is depicted change in the energy PDF, equipped with randomness 
in: only position (pos) and both position and material (pos+mat), respectively. To material properties 
on fine scale are assigned lognormal distributions with the same mean as in the previous experiment and $10\%$
of the variation. 
As one would suspect, the  PDF for the random position and material case 
has bigger spread than the one for only random position. Similar behaviour is observed 
in the PDF of the  upscaled bulk modulus - as shown  in Fig.~\ref{Fig:posall}. This spread is 
also observed when using 
different types of boundary conditions as shown in  Fig.~\ref{Fig:part_total_random_all} 
in which the left figure stands for the total uncertainty upper
$95\%$
limits (aleatory+epistemic), whereas the right figure stands for the $95\%$ limits of aleatory randomness only. Similar 
is depicted for the shear moduli in Fig.~\ref{Fig:bulk_random_all}.

%% file: example2.tex
\subsection{Upscaling of damage pheonmena}

In this subsection, the proposed approach is applied to another interesting problem. For this purpose a phenomenological elasto-damage model is considered as described in the beginning of this section. The goal is to compute a homogenized description of random material parameters on coarse-scale given fine-scale measurements. Fine-scale is assumed to follow same constitutive model as coarse-scale, however the material parameters are being random. For validation purpose they are also assumed to be homogeneous, and then in the second case scenario also spatially varying.  In both experiments we only simulate the displacement controlled biaxial compression of 2D block with unit length ( similar to the experiment in the previous example). The deformation tensor - applied to the boundary nodes of the specimen - is specified as:

\begin{equation}
\mathcal{F} = \begin{bmatrix}
-p & 0 \\
0 & -p 
\end{bmatrix}
\end{equation}
\
\\
where $p$ is the propotional factor for applied displacement, which varies with time increments $t_i$  according to the values in set $(t_i,p_i):\lbrace (0,0),(3,0.00025),(10,0.00035)\rbrace $  
\subsubsection{Validation}

\begin{figure}
\begin{center}
\includegraphics[width=0.45\textwidth]{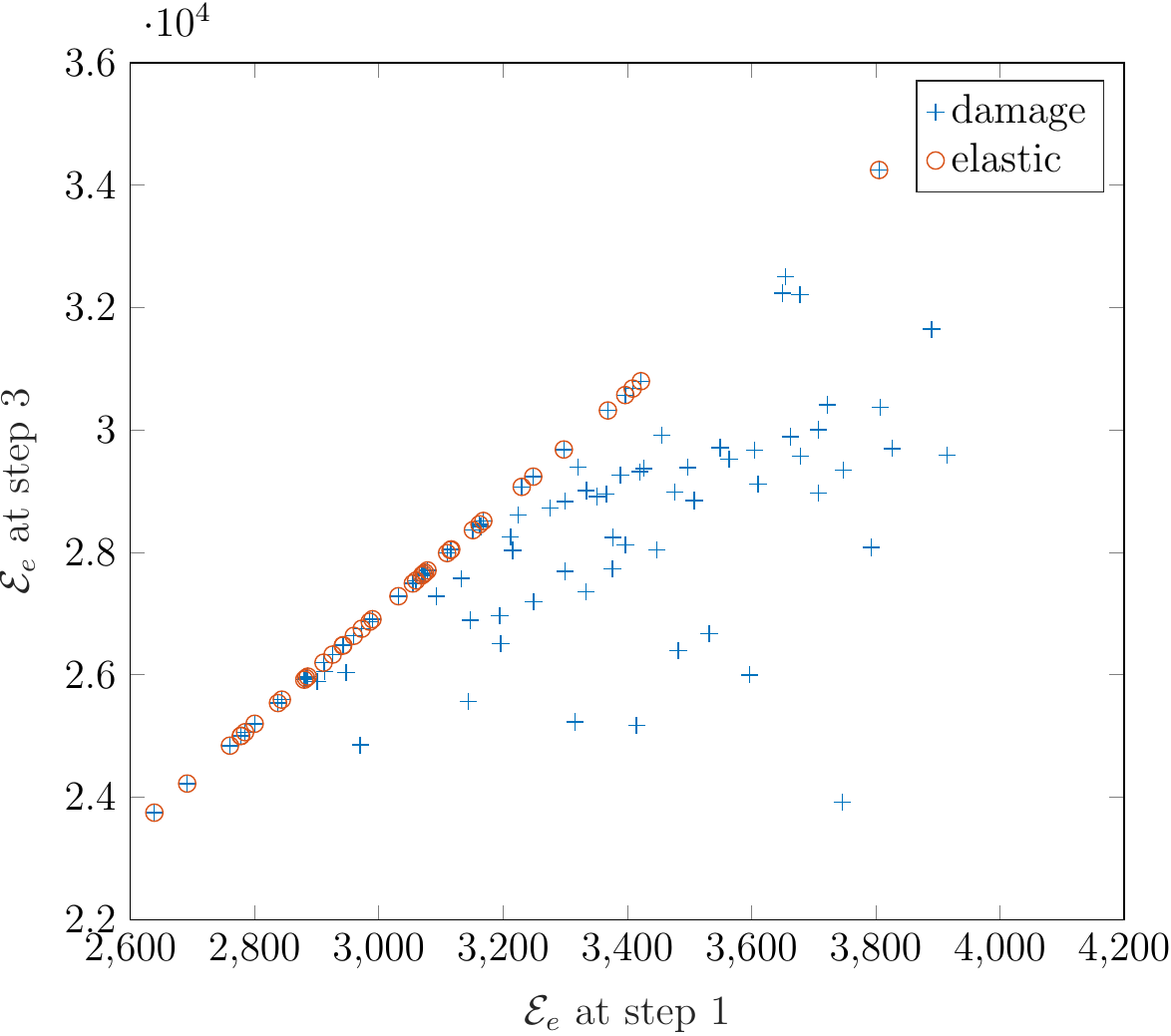}
\includegraphics[width=0.46\textwidth]{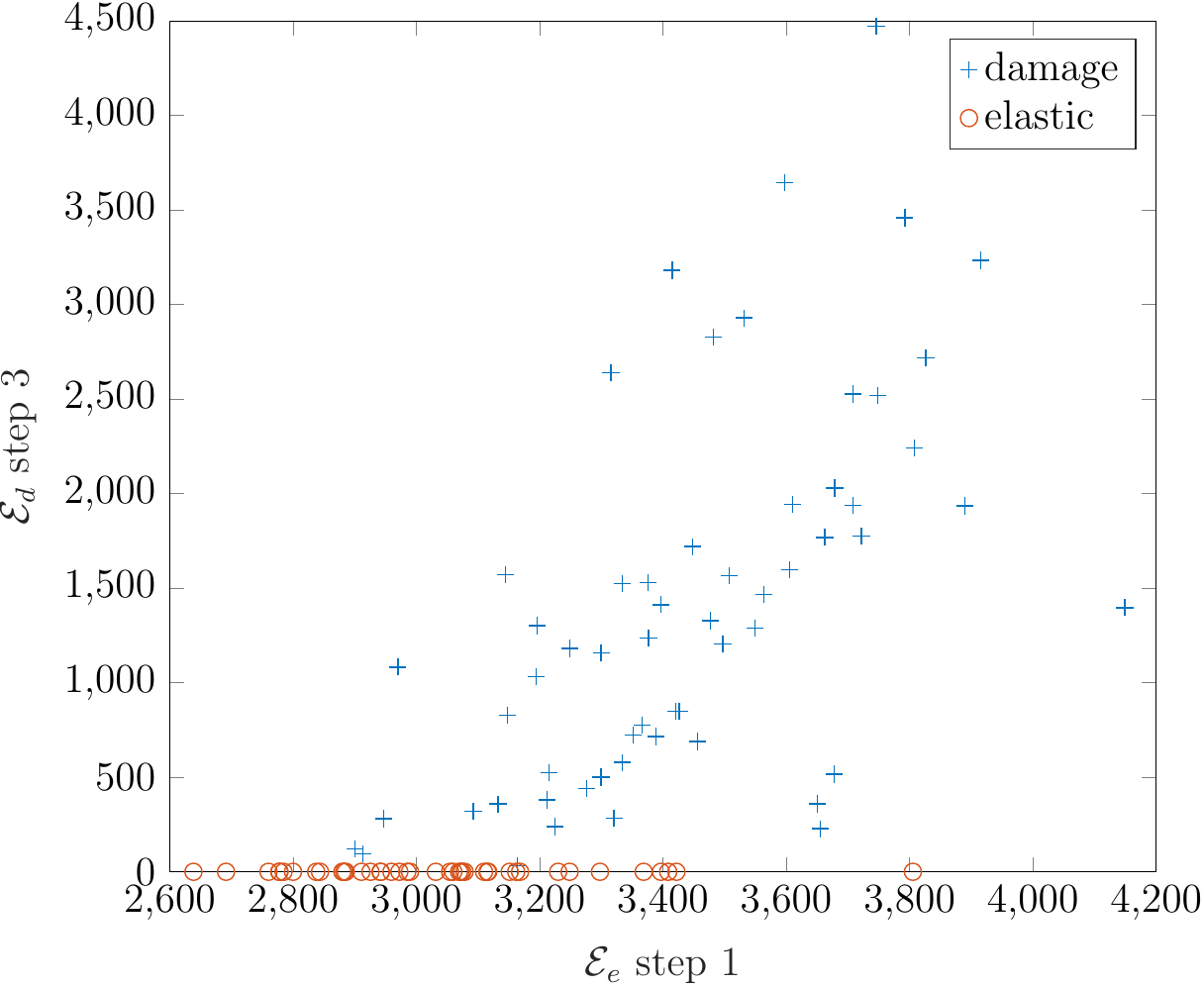}
\caption{Dependence of elastic ($\mathcal{E}_e$) and dissipative  ($\mathcal{E}_d$) energies between the linear elastic step 1 and nonlinear step 3}
\label{figvalidpic1pic2}
\end{center}
\end{figure}

\begin{figure}
\includegraphics[width=0.45\textwidth]{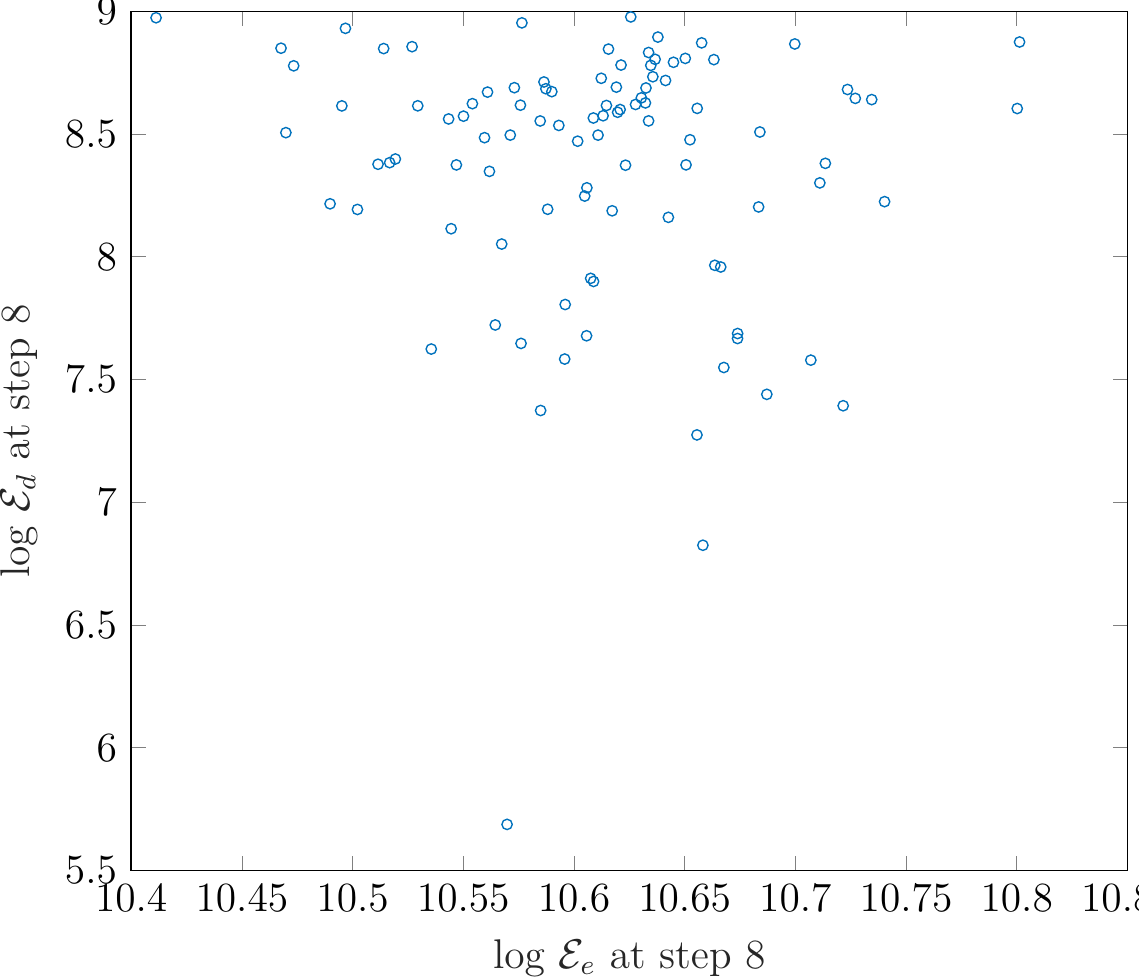}
\includegraphics[width=0.45\textwidth]{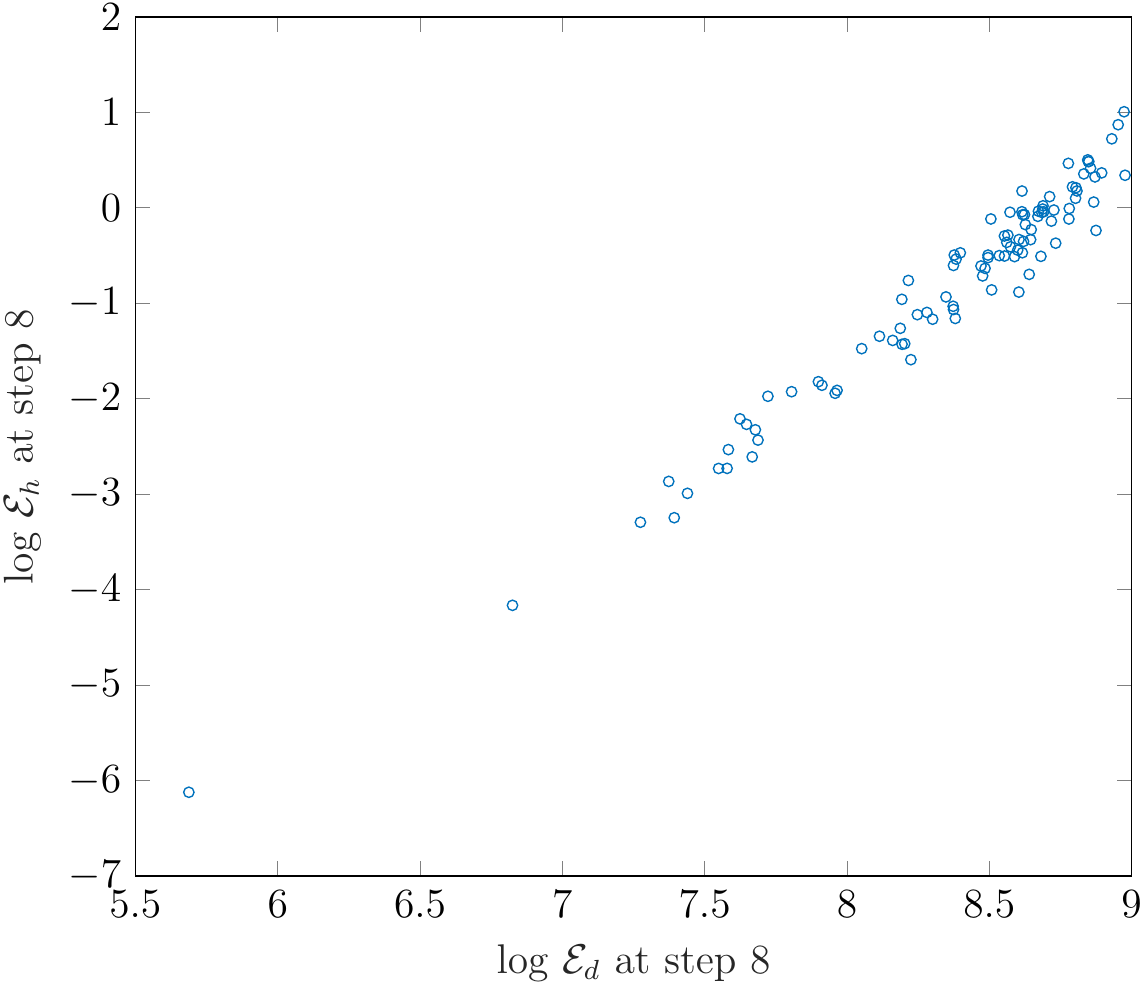}
\caption{Dependence of elastic ($\mathcal{E}_e$) and dissipative  ($\mathcal{E}_d$) energies between at the full damage step 8}
\label{fig_valid_pic3_pic4}
\end{figure}

\begin{figure}
\includegraphics[width=0.45\textwidth]{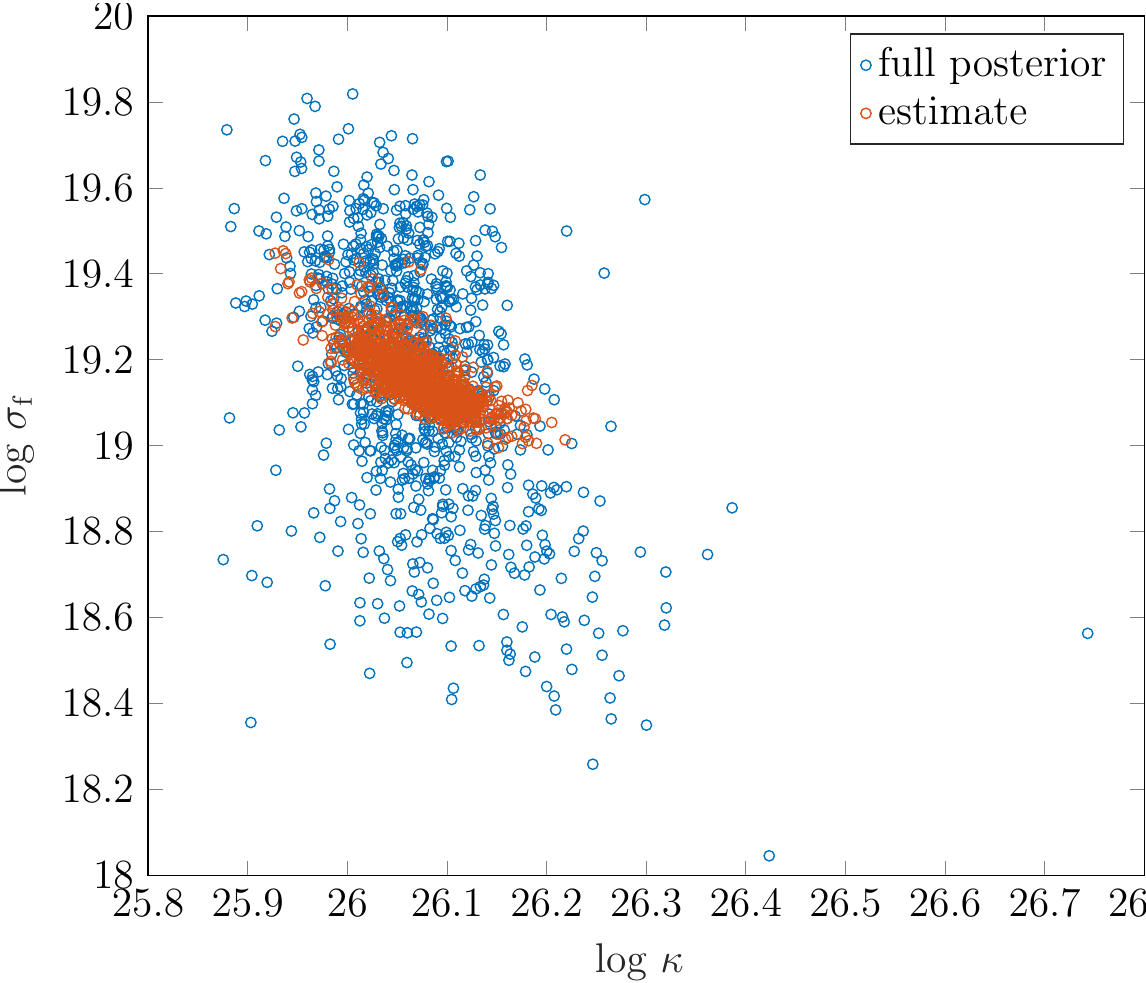}
\includegraphics[width=0.45\textwidth]{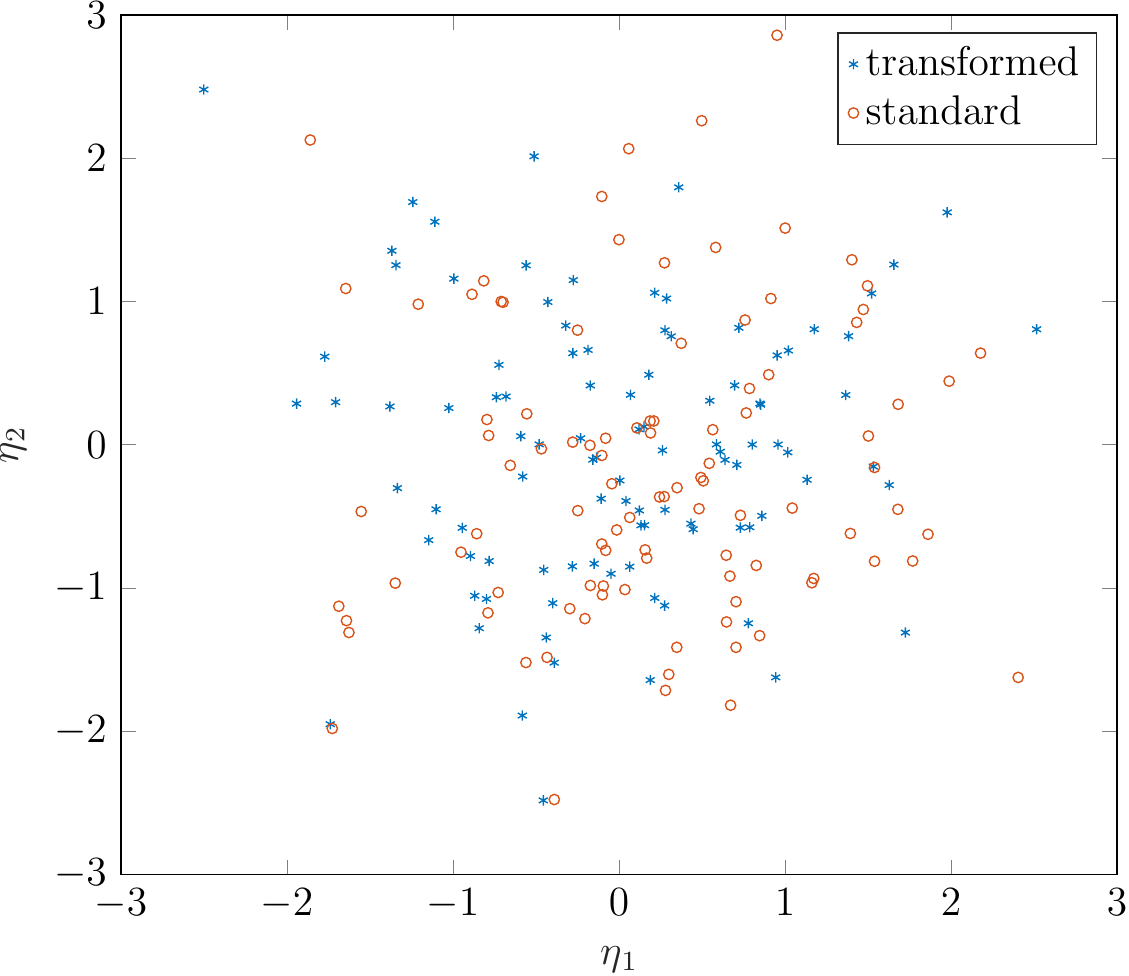}
\caption{Left: The relationship between the macro-scale estimated properties w.r.t. to full posterior measure (aleatory+epistemic uncertainty) and its epistemic mean (only aleatory uncertainty). Right: Comparison of 100 samples of mapped Gaussian random variables from the estimated macro-parameters and independent standard Gaussians. }
\label{fig_valid_pic7}
\end{figure}

\begin{figure}
\includegraphics[width=0.45\textwidth]{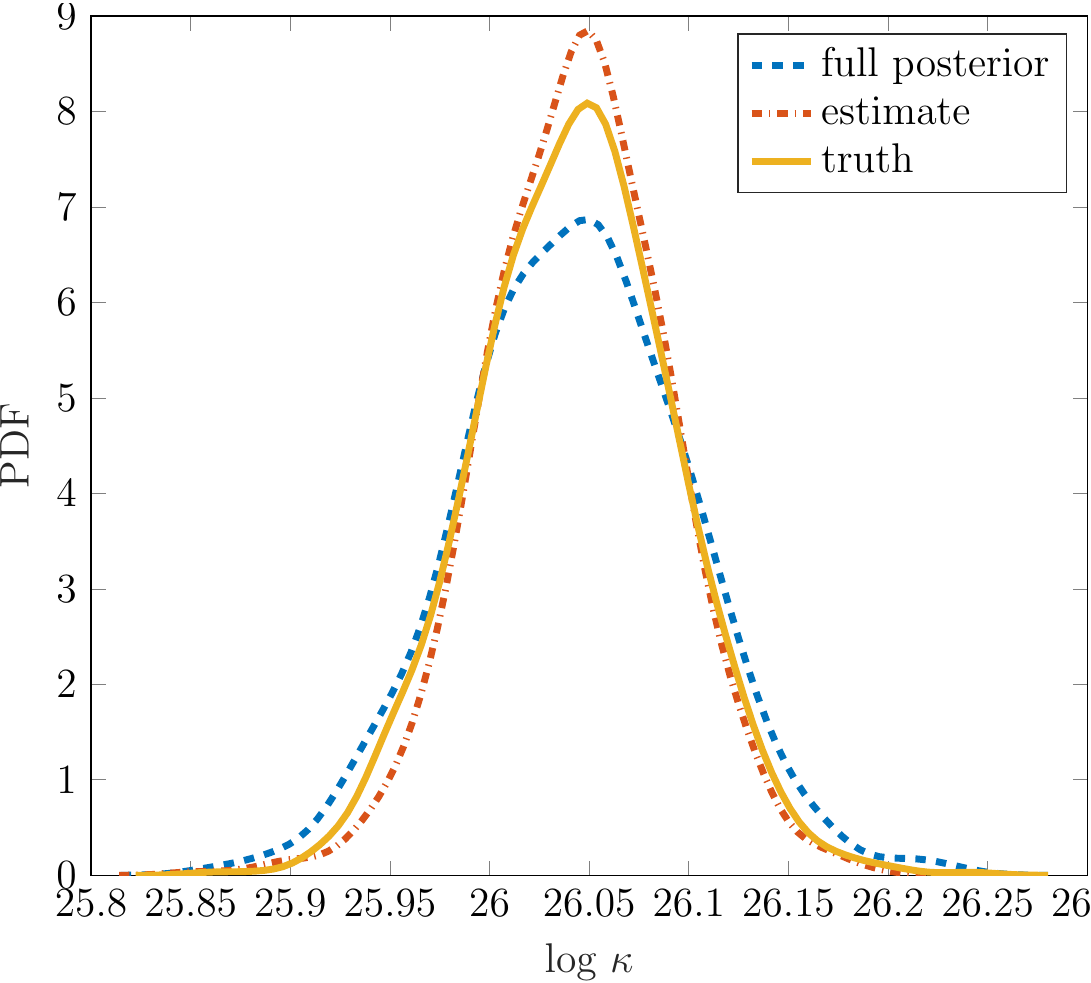}
\includegraphics[width=0.45\textwidth]{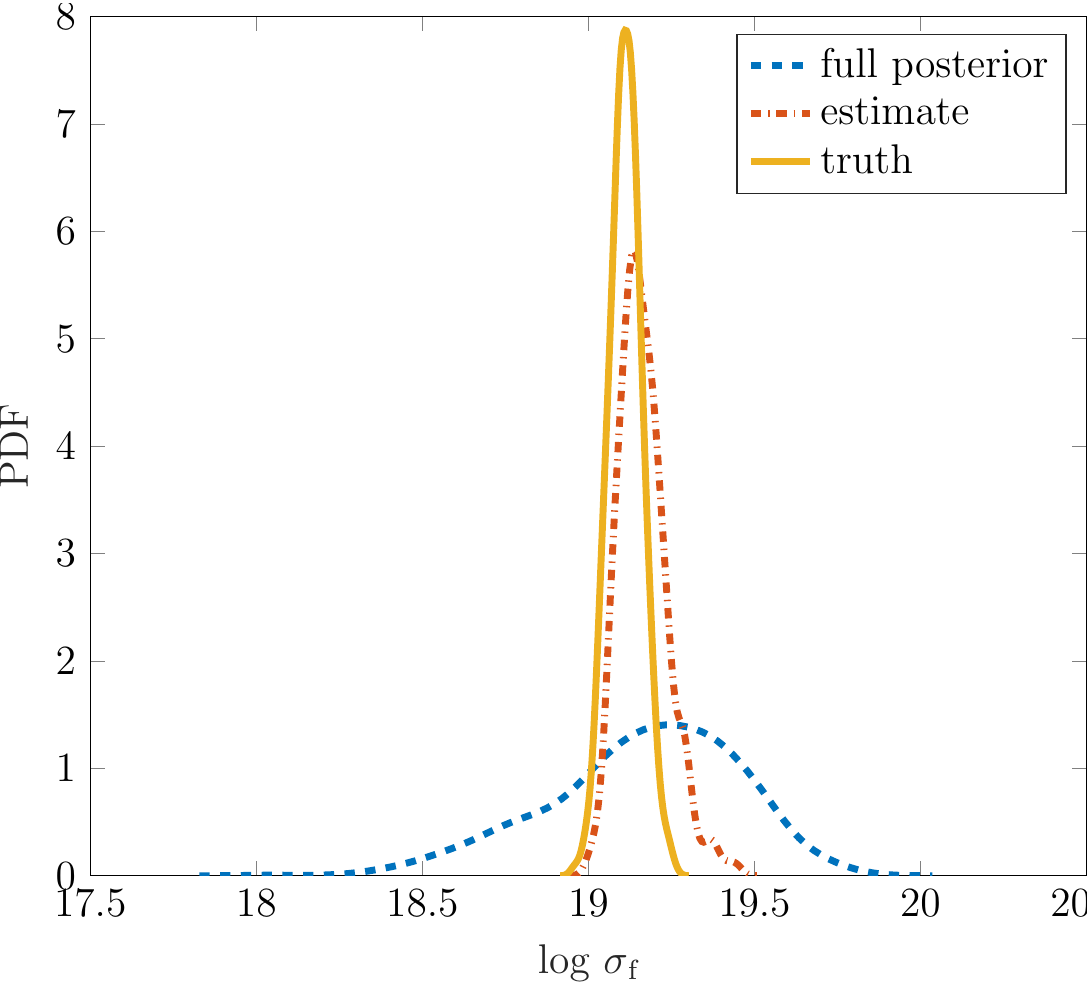}
\caption{Estimated posterior macroparameters w.r.t. to their true value. Full posterior represents both aleatory and epistemic uncertainty, whereas estimate is only aleaotry one.}
\label{fig_valid_pic5}
\end{figure}

For validation purposes the material 
parameters $\kappa_m(\omega_\kappa)$ on the fine-scale are modelled as lognormal random variables with the mean $\mu_\kappa:=\{2.0444,0.92,0.0031,0.0047)\cdot 10^{11}$ and the coefficient of correlation $5\%$. 
After propagating the variables through the elaso-damage model, the corresponding elastic, damage and hardening energies are estimated. We assume that the polynomial chaos approximation of the energy measurement is not given, but only a set of 100 samples. Therefore, 
the log of energies are modelled as copula Gaussian mixtures with the unknown number of components. The simulation is run in 8 equidistant time steps, the first two being elastic. In the third to sixth steps the behaviour is combination of elasticity and damage, whereas in the last step is predominated by the damage component. In Fig.~\ref{figvalidpic1pic2} are shown scatter plots of energies in the first and the third step, both depicting two states in the response: elastic and damage. 
The linear components are linearly related as expected, whereas the correlation between the damage and the elastic step is random. To uncouple 
the measurement data, we observe energies at the last time step (full damage state) as depicted in Fig.~\ref{fig_valid_pic3_pic4}.  Clearly, the hardening and damage energies are almost linearly related in the log-space, whereas 
this doesnt hold for the damage and elastic part of energy.  Both copula and inividual components are estimated as discussed in Section 4. In order to decouple the complete set, the inverse transform is applied in order to obtain uncorrelated Gaussian samples. These are further used to generate the polynomial chaos expansions for all time steps by Bayes's rule as presented in Section 3. In this manner one may obtain prediction of energies at the meso-scale for all time steps. 
These are further adoped as approximation of the measurement $y_m$. Furthermore, one assumes that the material parameters used for prediction of coarse-scale simulation are not known, and are assumed to be uncertain. Due to positive definitness requirement they are also modelled as lognromal random variables with the mean $20\%$ bigger than in the fine-scale case and coefficient of correlation $20\%$.
The resulting updated coarse-scale properties are shown in Fig.~\ref{fig_valid_pic5}. The bulk moduli and the limit stress that initiates the damage are both updated and match the 
true distribution, whereas their correlation and the mapping to the normal space is shown in Fig.~\ref{fig_valid_pic7}. On the other hand, due to chosen experiment both shear and hardening moduli stay unindentified as they are not observable. 

\begin{figure}
\includegraphics[width=0.45\textwidth]{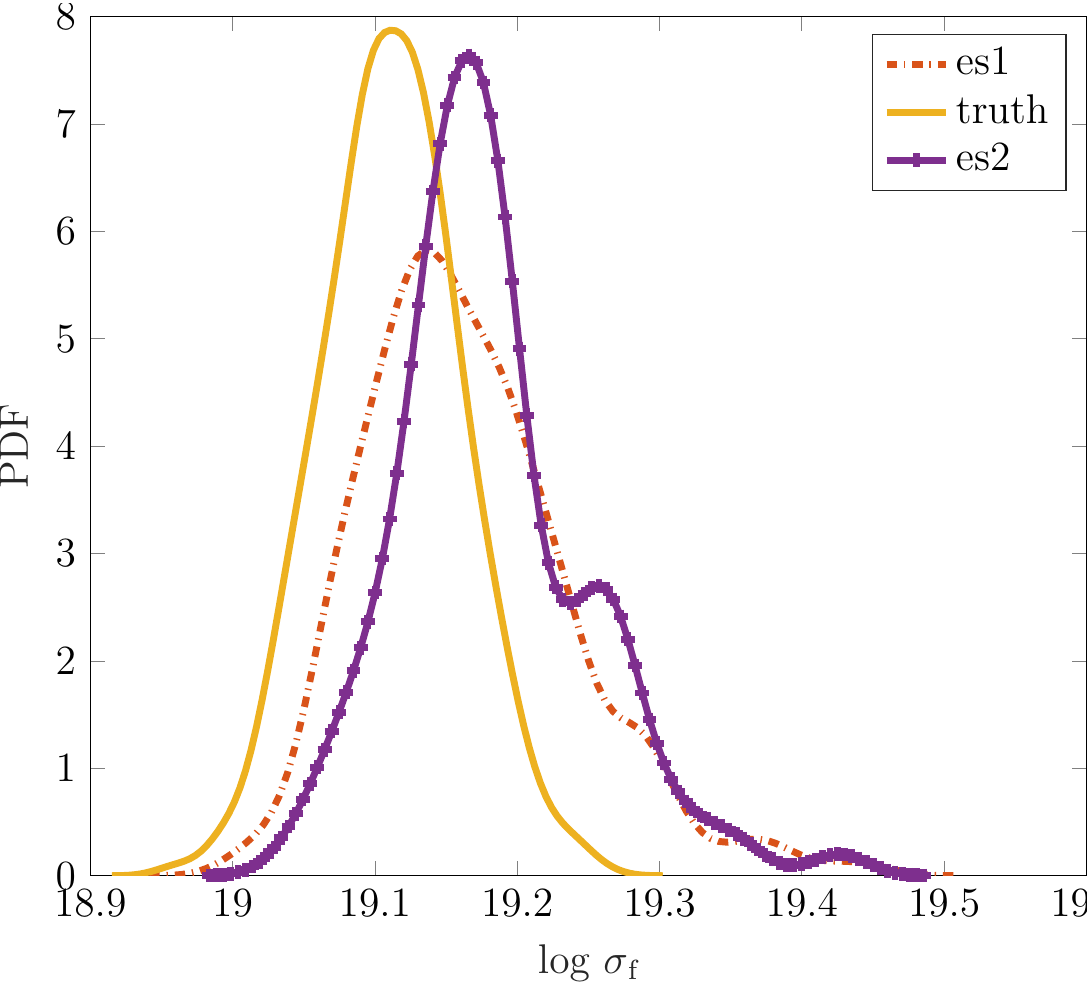}
 \includegraphics[width=0.45\textwidth]{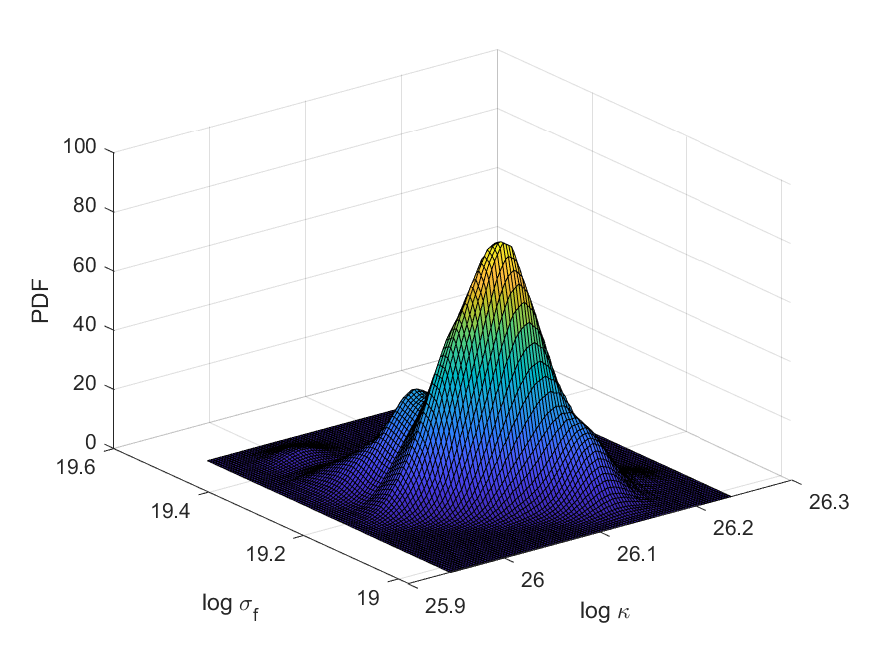}
\caption{Left: comparison of two upscaling strategies, right: correlation between the macro-scale parameters.}
\label{fig_valid_pic8}
\end{figure}

In the previously described experiment the relationships between the measurement data and their approximations are too complex in order to be properly modelled. Therefore, the experiment is repeated in same setting, only this time the measurement is not functionally aproximated. Instead, the inverse problem is solved for each individual sample of measurement (each RVE), and then the updated parameters are collected into the set of parameter samples. 
This calculation is expected to be simpler than the previous one as the relationships between the material parameters are easier to model. By collecting mean value of posterior distributions,  one may model the set of parameters as copula Gaussian mixture similarly to the previous case. After mapping to the Gaussian space the 
corresponding coarse-scale estimates can be functionally approximated by a polynomial chaos expanion obtained by Bayes' rule. In Fig.~\ref{fig_valid_pic8} is depicted the difference between this approach (est1) and the previous one (est2), as well as  joint distribution between the bulk moduli and $\sigma_f$. 

\subsubsection{Upscaling of heterogeneous medium}

\begin{figure}[!h]
\begin{center}
\includegraphics[width=0.6\textwidth]{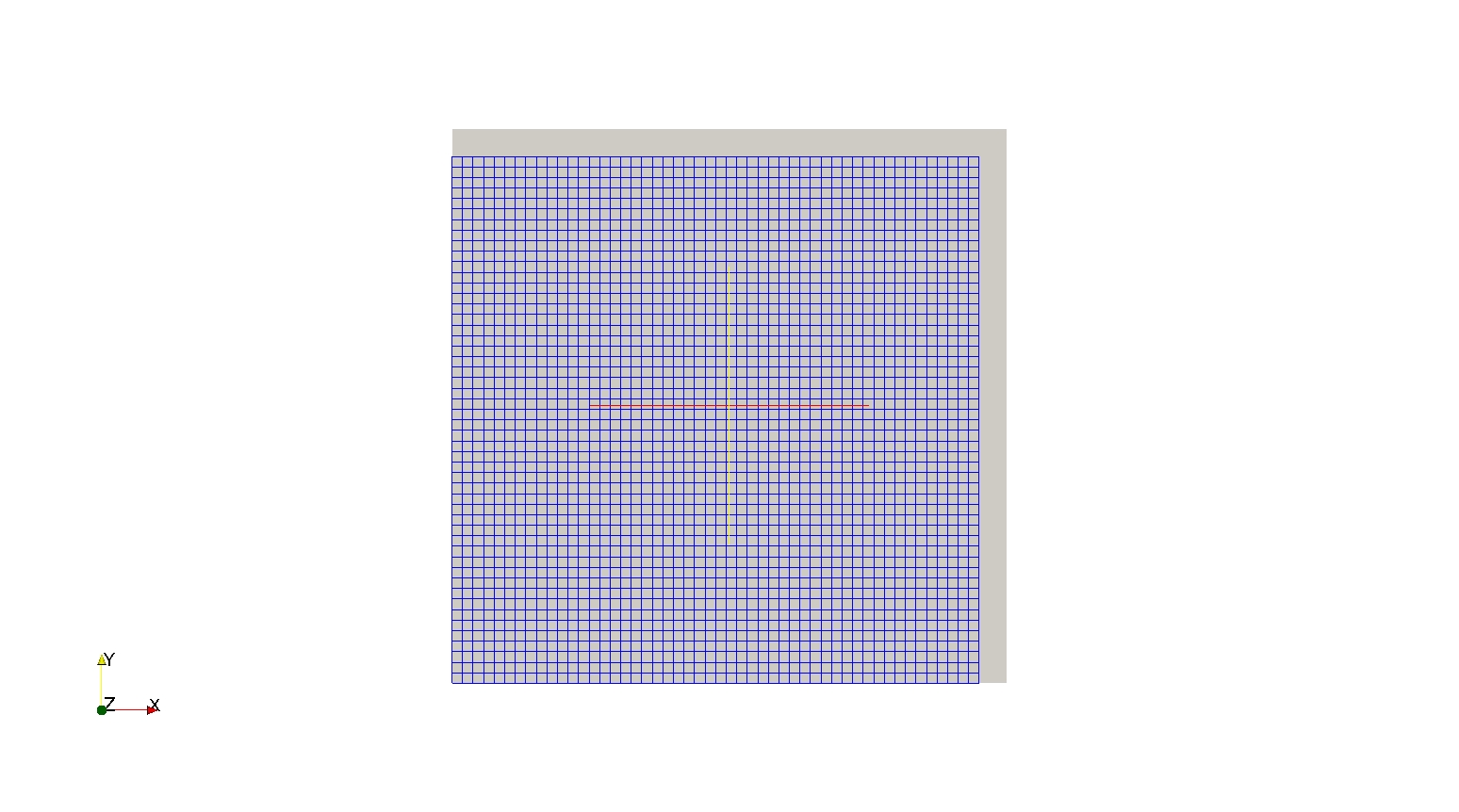}
        \caption{The deformation of the 2d block by bi-axial compression}
   \label{fig:update_homo_Plas}    
\end{center} 
\end{figure}

 The model is implemented in a finite element framework
( see \cite{ibrahimbegovic,sadiq} for more details) and is used to demonstrate the proposed up-scaling strategy. For this purpose, a 2D square
block of unit length is considered. To emulate the coarse and fine-scale descriptions in the finite element setting, the 2D block
is considered as one element on the coarse-scale, whereas the fine-scale comprises of 2500 elements. The block is deformed by discplacement controlled bi-axial compression as 
shown in Fig.~\ref{fig:update_homo_Plas}.

\begin{figure}[!h]
\subfloat[$\ell_c=5\ell_e, c_{var}=5\%$]{
      \includegraphics[width=0.35\textwidth]{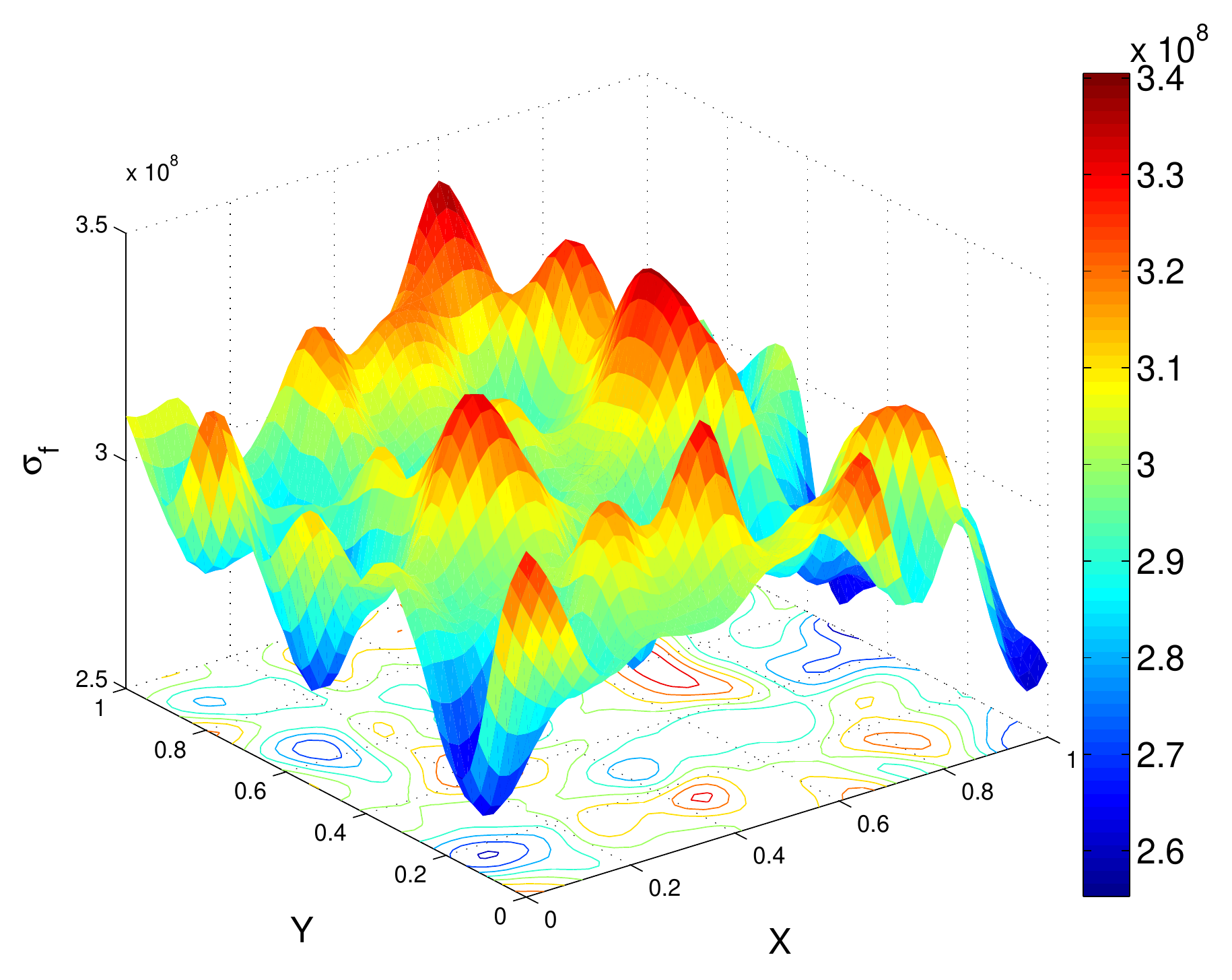}
     }
     \subfloat[$\ell_c=10\ell_e, c_{var}=5\%$]{
       \includegraphics[width=0.35\textwidth]{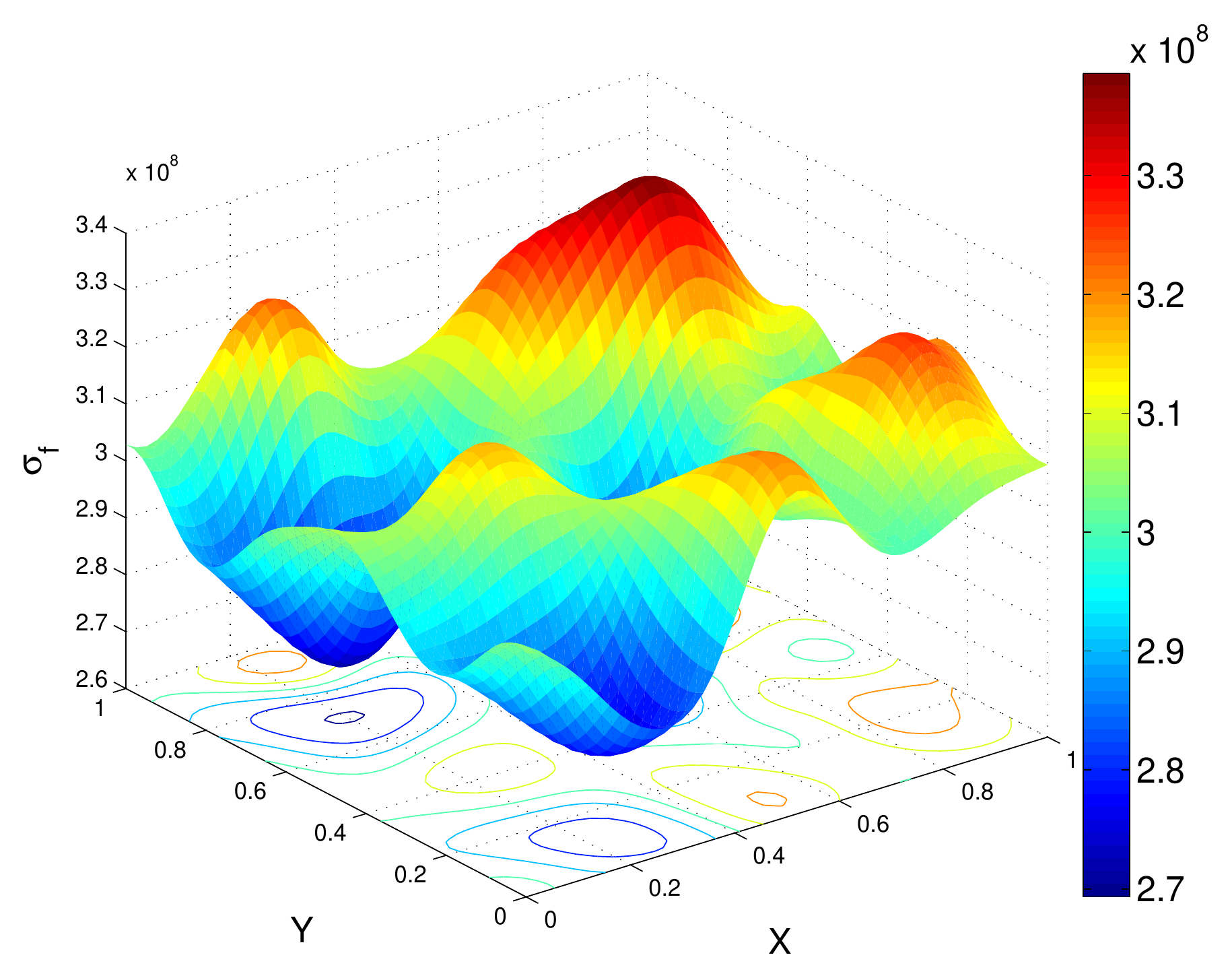}
     }
 \subfloat[$\ell_c=25\ell_e, c_{var}=5\%$]{
       \includegraphics[width=0.35\textwidth]{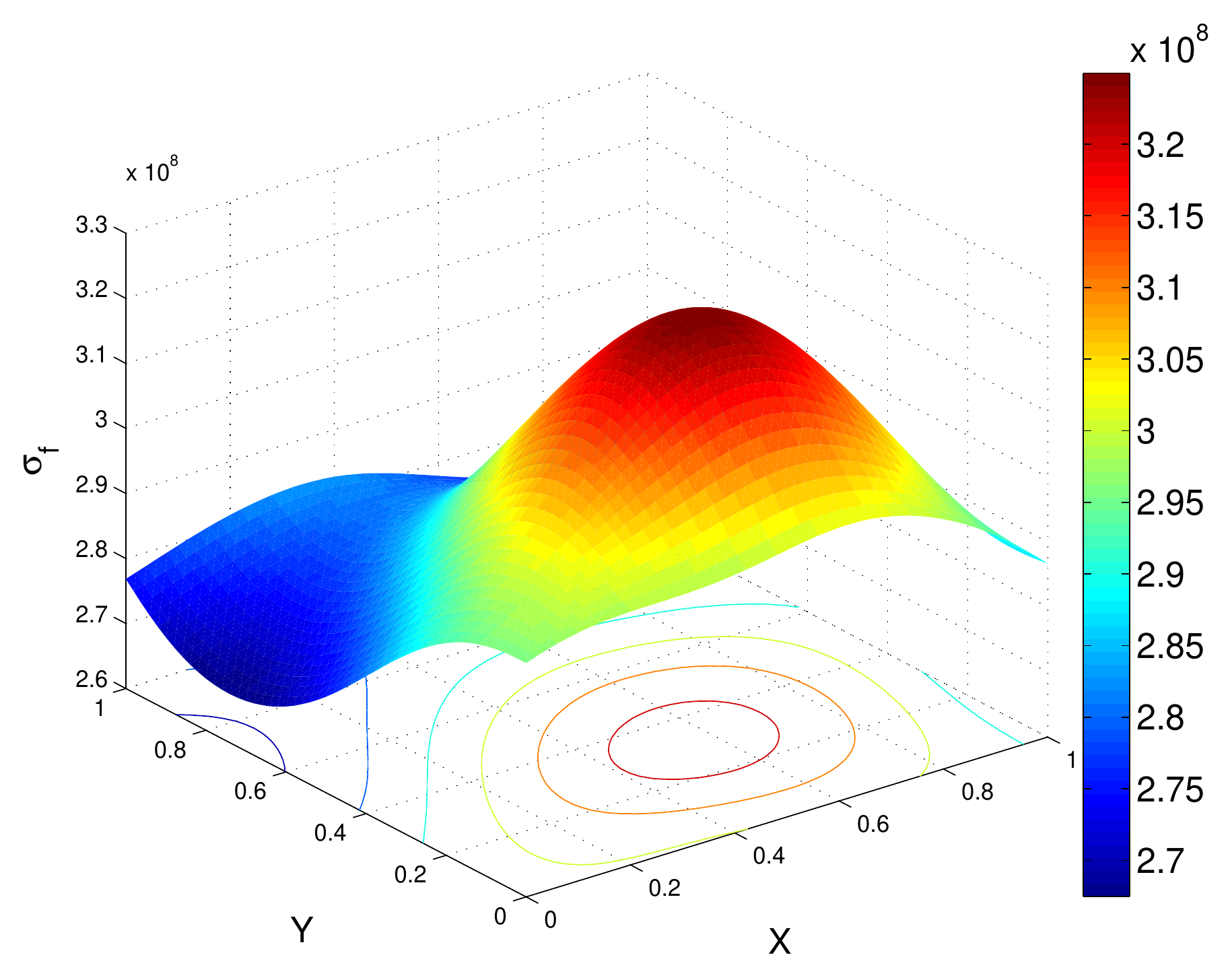}
     }
     \caption{The damage stress realisations using different values of the correlation length $\ell_c$}
     \label{fig:update_homo_Plas1}     
\end{figure}

As far as material description is concerned, the material properties on the heterogeneous fine-scale are a priori assumed to be realizations 
of log-normal random fields with the statistics depicted in \reftab{tab:damageprior}, and Gaussian covariance functions. These are simulated using different values of
correlation length $l_c \in \lbrace 5l_e,10l_e,25l_e \rbrace$ ($l_e$ is  the element length on the fine-scale)  and
 coefficients of variation $c_{var} \in \lbrace 5\%,10\%\rbrace$.  In Fig.~\ref{fig:update_homo_Plas1} is shown an example 
of the meso-scale random field realisations given different correlation lengths. The realisation is becoming smoother when  
the correlation length increases. This means that the material becomes homogeneous in the limit $\ell_c=\infty$. On the other hand, 
the coarse-scale material properties are taken a priori as a lognormal random variables, with the same mean and standard deviation as their meso-scale counterpart.

\begin{figure}
\begin{center}
\includegraphics[width=0.45\textwidth]{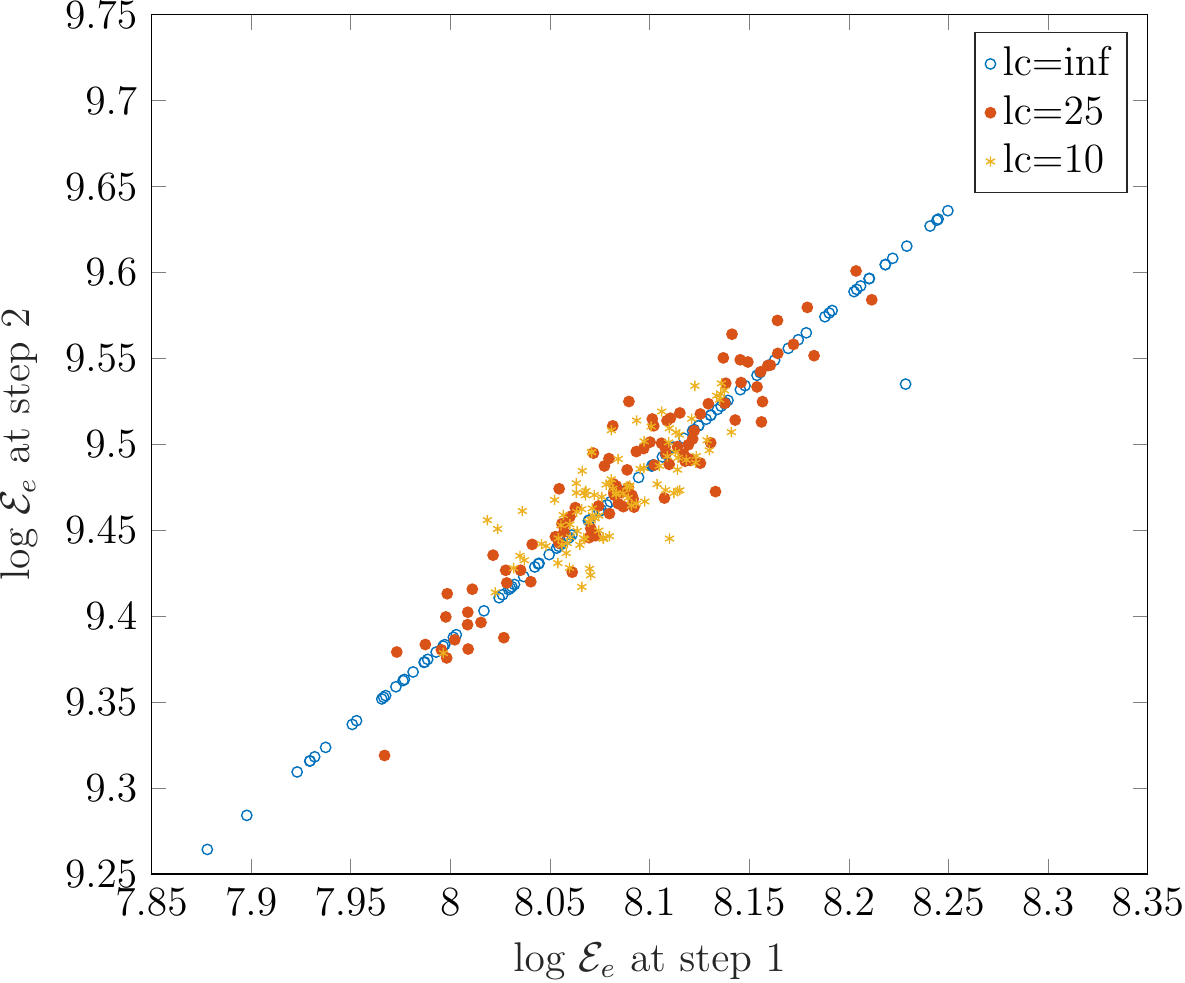}
\includegraphics[width=0.45\textwidth]{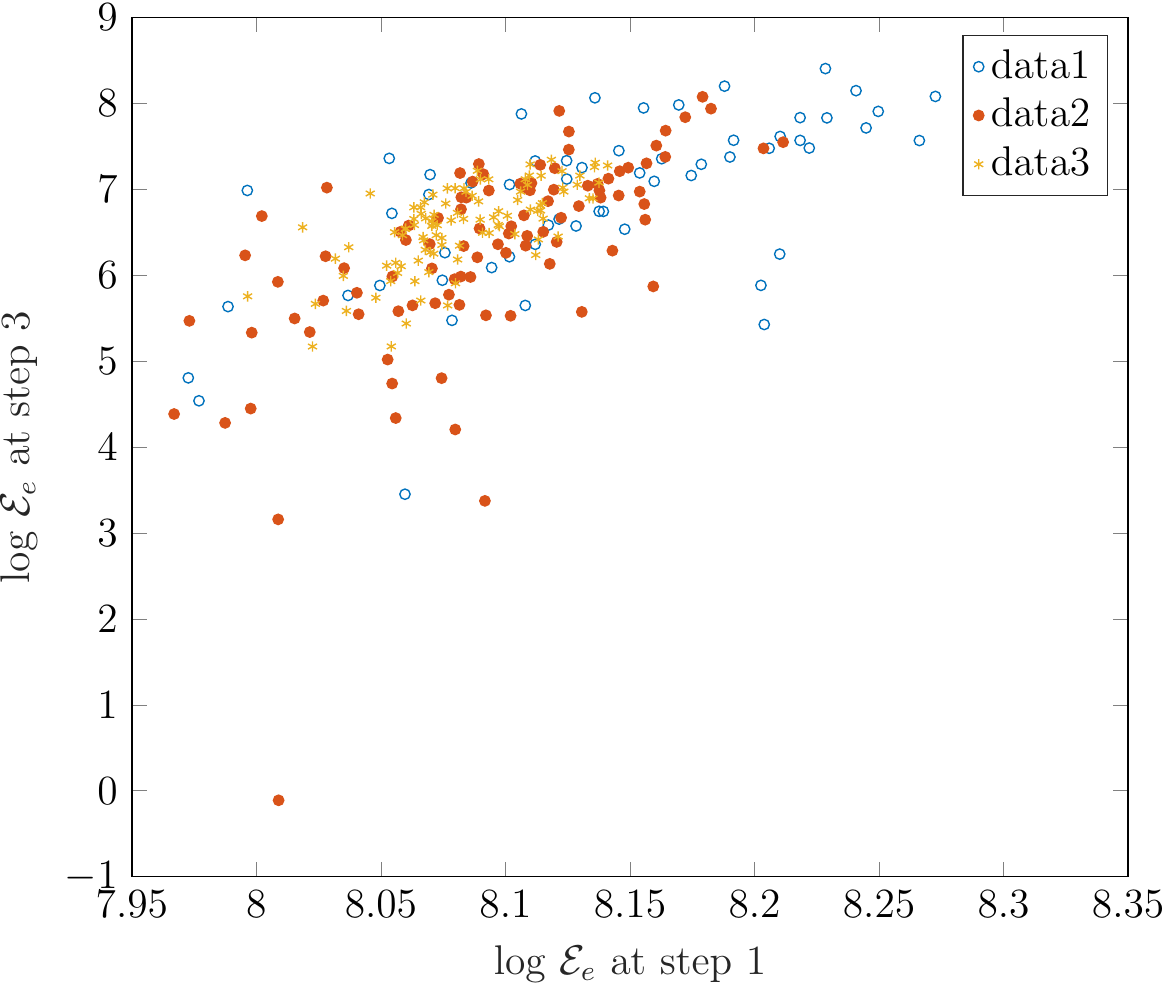}

\vspace*{20 mm}

\includegraphics[width=0.47\textwidth]{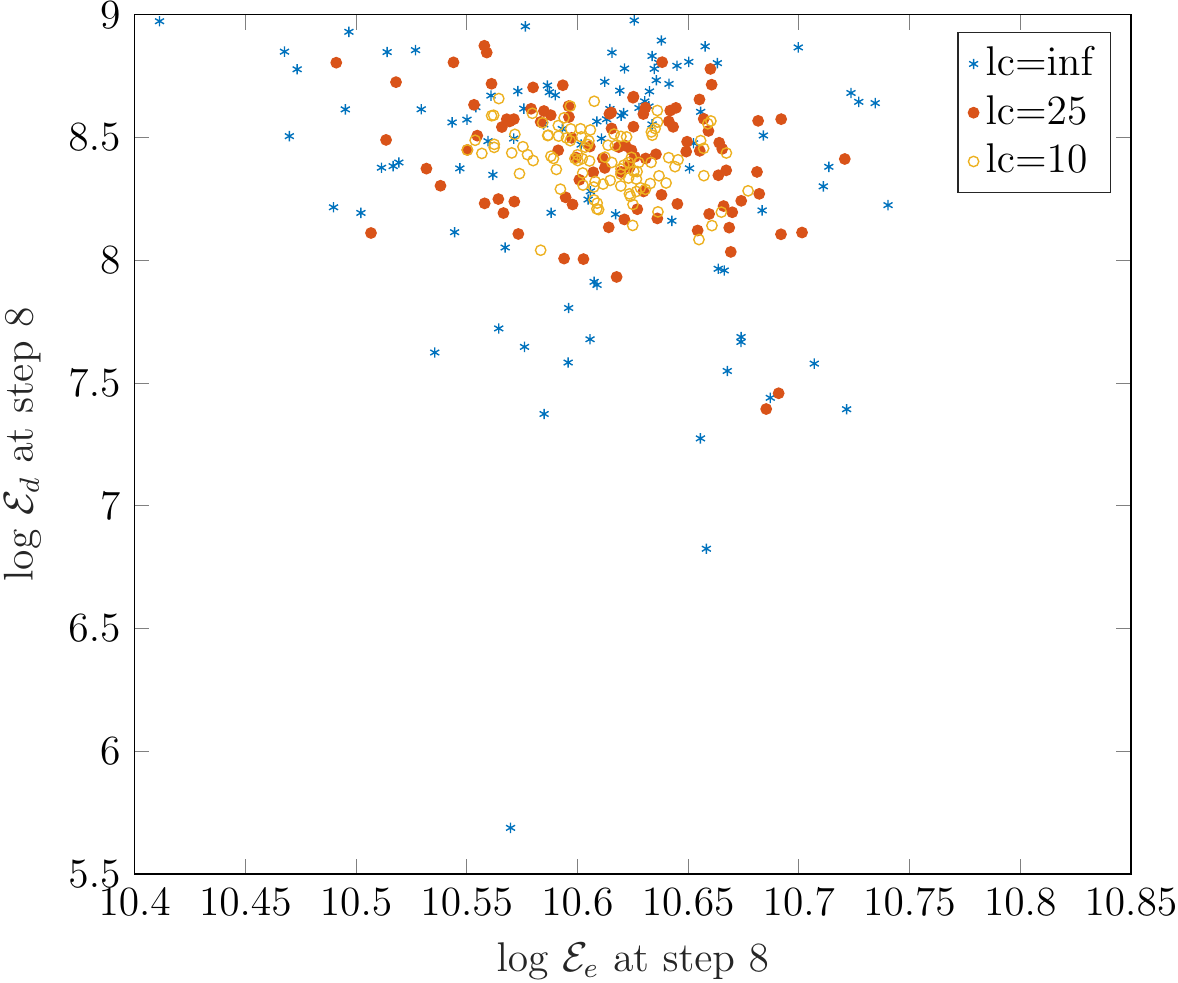}
\includegraphics[width=0.45\textwidth]{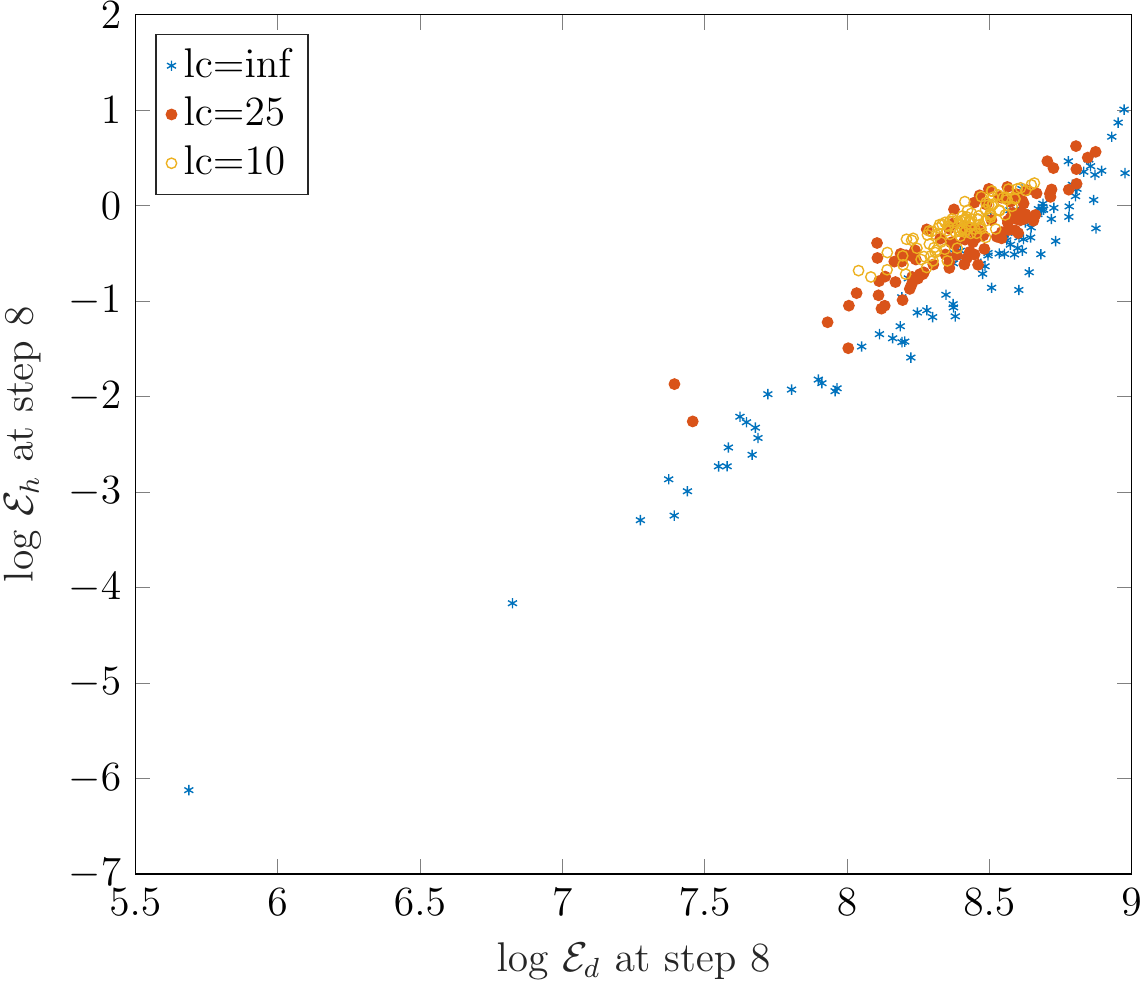}
\end{center}
\caption{Energies dependence w.r.t to different time steps}
\label{figvalidpic8a}
\end{figure}

The measurement data are made of three type of avaraged energies: the stored elastic energy $\mathcal{E}_e$ (see \refeq{measur2}), the damage dissipation $\mathcal{E}_d$ (i.e.~the first term in \refeq{measur2} ), and the hardening part $\mathcal{E}_h$ of the dissipation energy (i.e.~the second term in \refeq{measur2}). Their logarithms are simulated using mixture models and vine copulas, and further identified using variational Bayes rule. Similarly to 
the experiment in the validation section, in the first two simulation steps one may observe only elastic energy as the dissipation effects do not appear. 
Therefore, we start the simulation with the last step, and approximate the corresponding energies by  mixture models and vine copulas, see Fig.~\ref{figvalidpic8a}. The complete simulation 
steps from the previous section are repeated. 

\begin{figure}
\begin{center}
\includegraphics[width=0.4\textwidth]{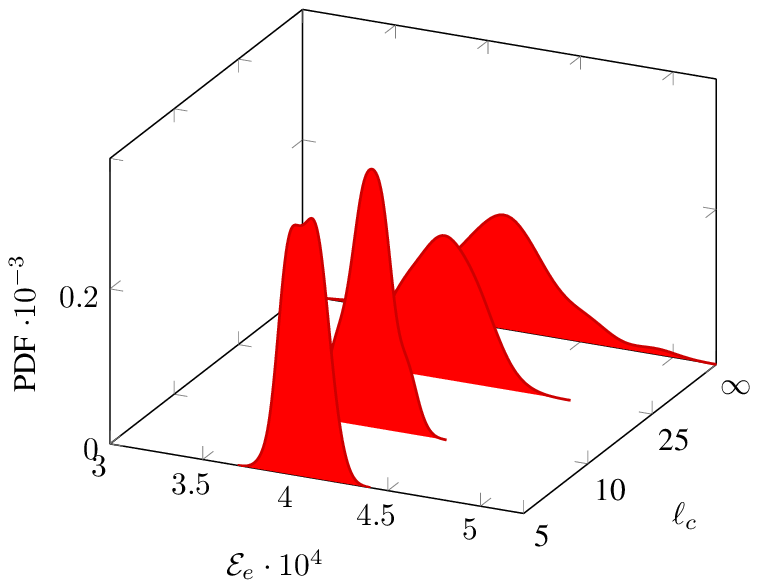}
\includegraphics[width=0.4\textwidth]{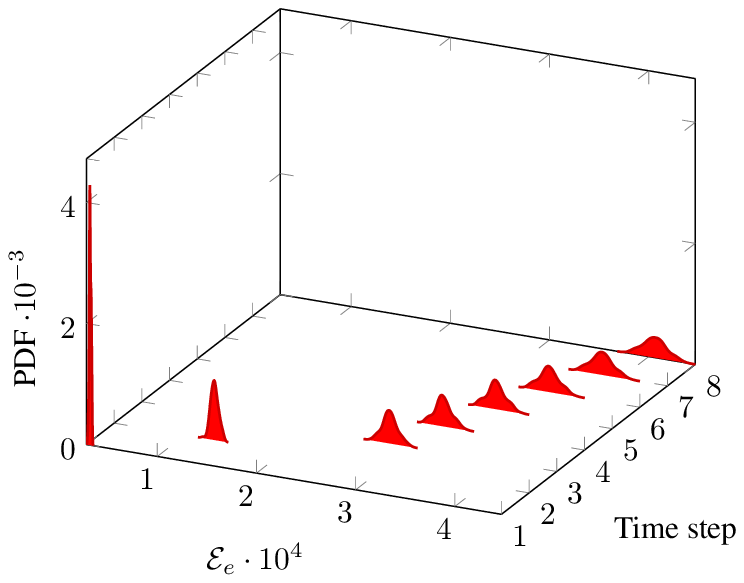}
\includegraphics[width=0.4\textwidth]{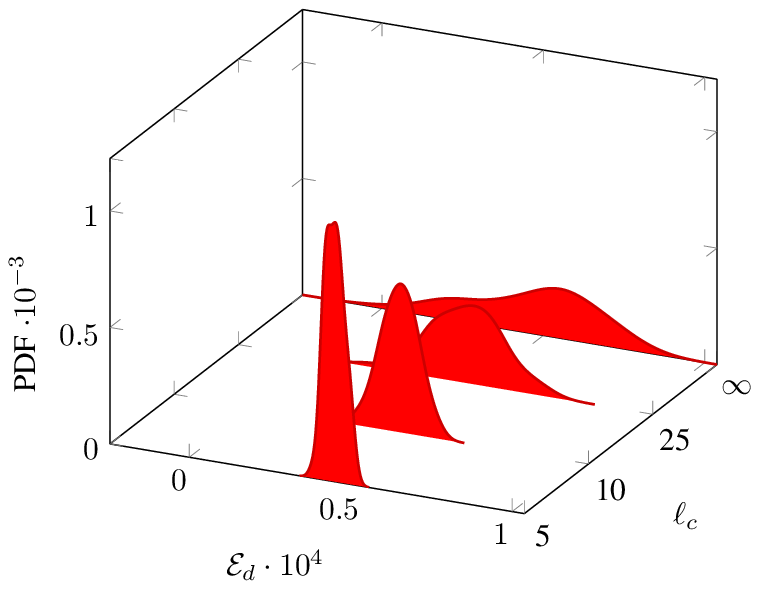}
\includegraphics[width=0.4\textwidth]{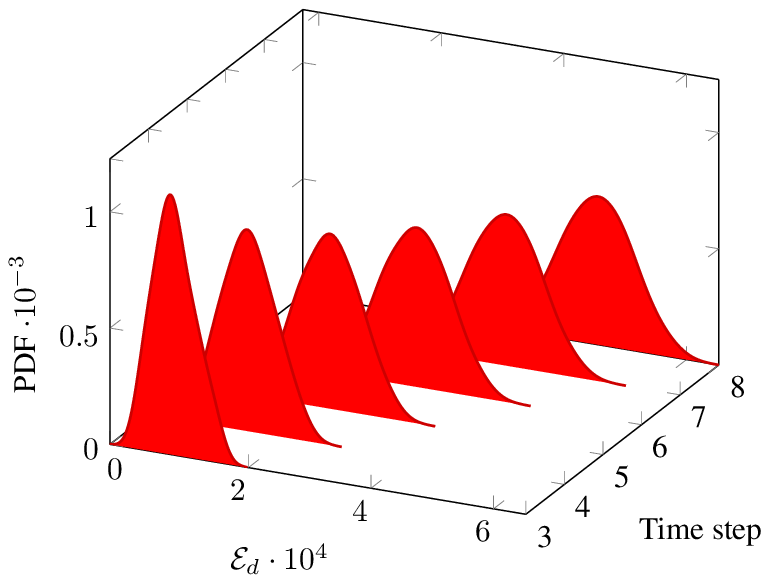}
\includegraphics[width=0.4\textwidth]{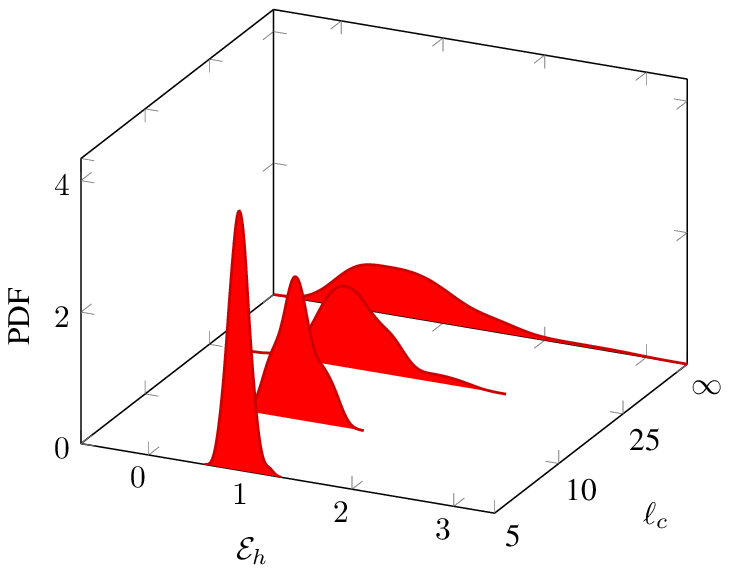}
\includegraphics[width=0.4\textwidth]{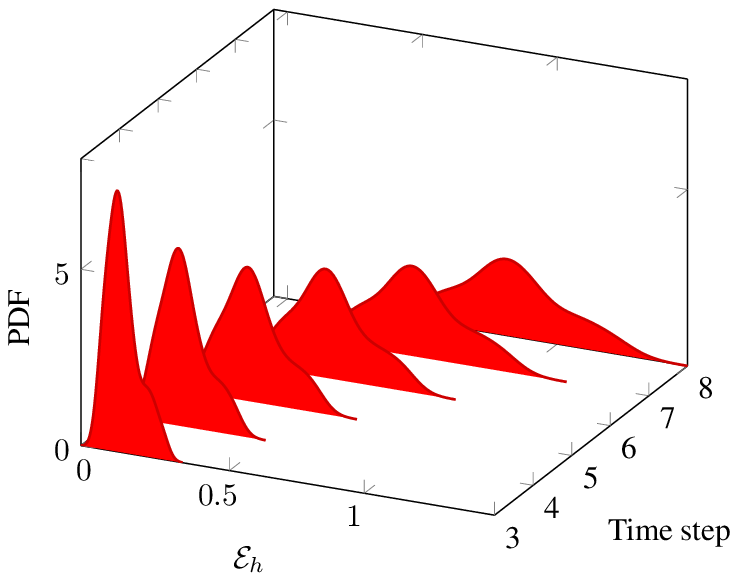}
\caption{The energy PDF w.r.t. the correlation length (left) $\ell_c$ and time steps (right)}
\label{fig_pure_shear_ran_pos_PR_ll}
\end{center}
\end{figure}

\begin{figure}
\begin{center}
\includegraphics[width=0.3\textwidth]{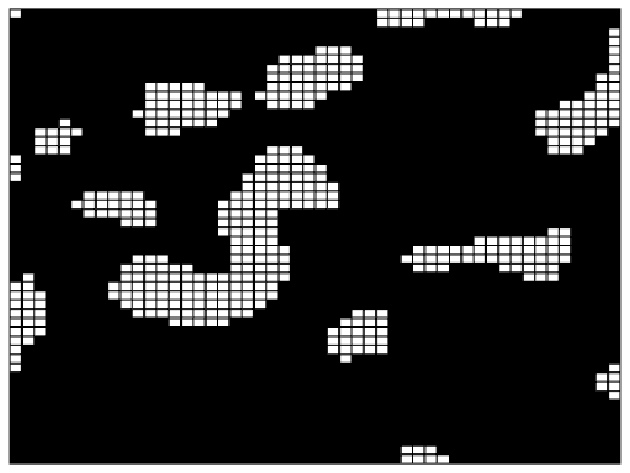}
\includegraphics[width=0.3\textwidth]{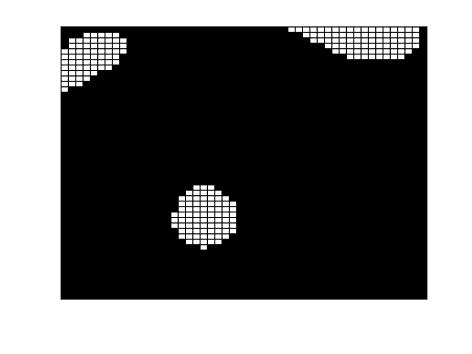}
\includegraphics[width=0.3\textwidth]{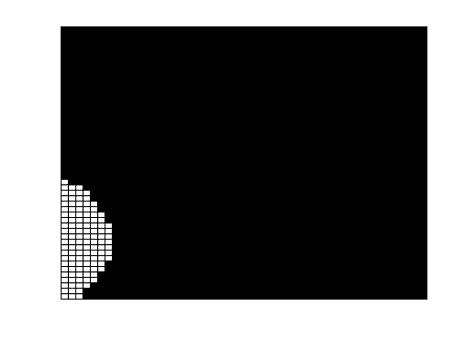}
\caption{The presence of damage on the meso-scale for $c_{var}=0.1$ and $\ell_c=5\ell_e$,  $\ell_c=10\ell_e$, $\ell_c=25\ell_e$ }
\label{damagfig1}
\end{center}
\end{figure}

As shown in Fig.~\ref{fig_pure_shear_ran_pos_PR_ll} the variation of energies increases with the correlation length size for the case when the $c_{var}$ of the meso-scale random field is taken to be $10\%$. The reason for 
this is that energy realisations are less fluctuating with increasing correlation length, but their avarage value is more pronaunced as prospective fluctuations do not cancel out, as similar can be concluded when observing Fig.~\ref{fig:update_homo_Plas1} on the left hand side. This holds for all energy estimates, and can be explained by Fig.~(\ref{damagfig1}) in which the presence of damage on one of the responses is shown. With increase of the correlation length the damage is more pronaunced and hence one expected higher variations. 
In addition, one may also conclude that the corresponding PDFs are becoming more skewed when the material model approaches homogeneous case.
The skewness in terms of long tailes is not completely caused by variations of the random meso-scale, but also by inaccuracy of the variational method used for the PCE estimation of the measurement due to overestimation. 

\begin{figure}
\begin{center}
\includegraphics[width=0.3\textwidth]{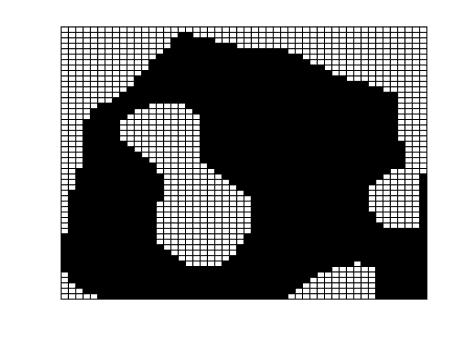}
\includegraphics[width=0.3\textwidth]{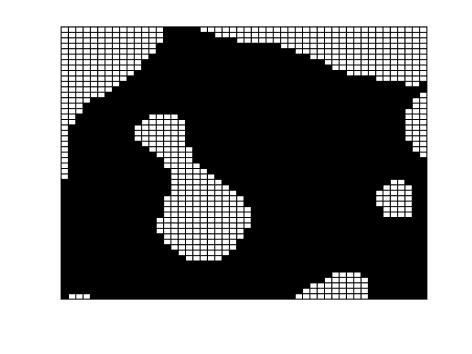}
\includegraphics[width=0.3\textwidth]{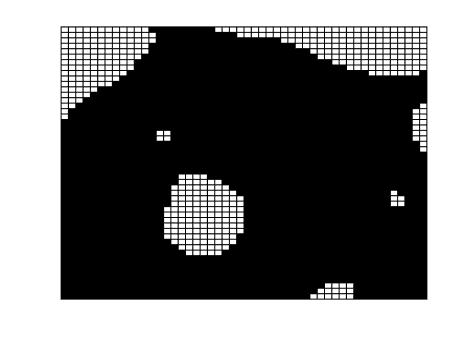}
\includegraphics[width=0.3\textwidth]{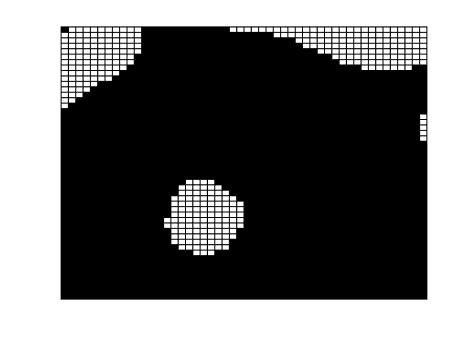}
\includegraphics[width=0.3\textwidth]{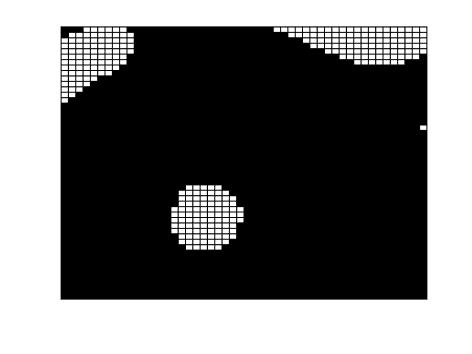}
\includegraphics[width=0.3\textwidth]{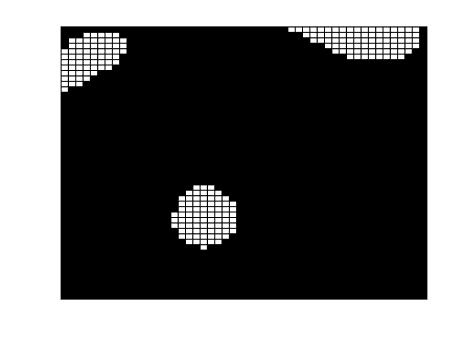}
\caption{The presence of damage on the meso-scale w.r.t. time step for $c_{var}=0.1$ and $\ell_c=10\ell_e$}
\label{damage_time}
\end{center}
\end{figure}

 On the right side of Fig.~\ref{fig_pure_shear_ran_pos_PR_ll} one may 
observe the energy evolution w.r.t. to time. Here, $c_{var}$ of the corresponding meso-scale random field is chosen to be $10\%$ and the correlation length is $\ell_c=10\ell_e$ . The top figure depicts the elastic energy. As expected, the energy variation grows in time. On the contrary, the damage energy seen in the middle does not alter much the PDF form. The damage initialises in the thrid step, and mostly shifts towards higher average value due to increased presence of damage as shown in Fig.~\ref{damage_time}.
 Finally, the hardening energy increases, but also changes the PDF form significantly in time. 

\begin{figure}
\begin{center}
\includegraphics[width=0.4\textwidth]{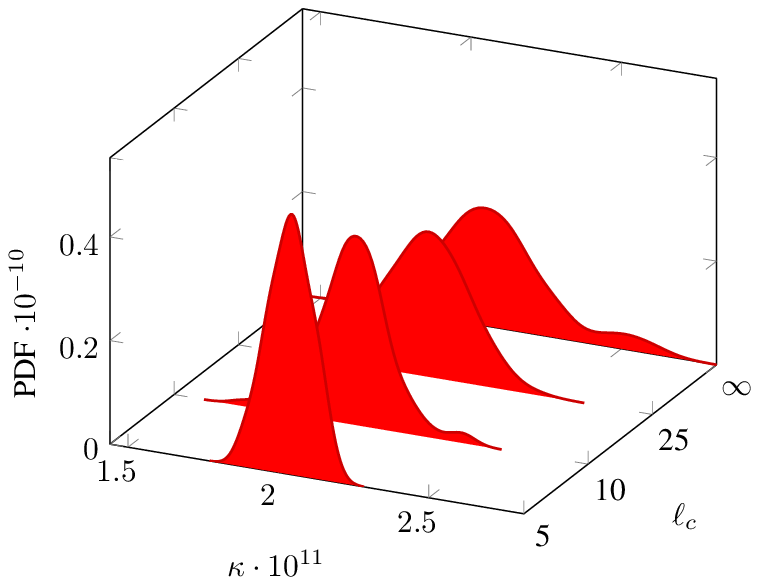}
\includegraphics[width=0.4\textwidth]{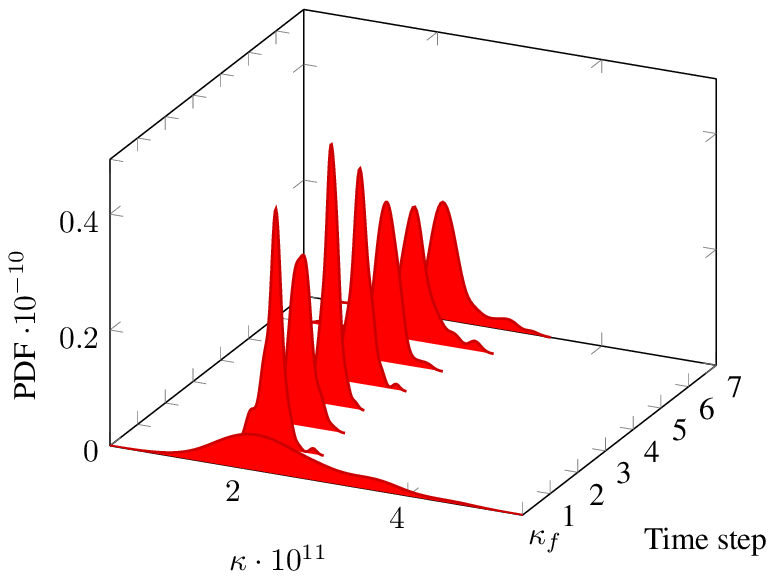}
\includegraphics[width=0.4\textwidth]{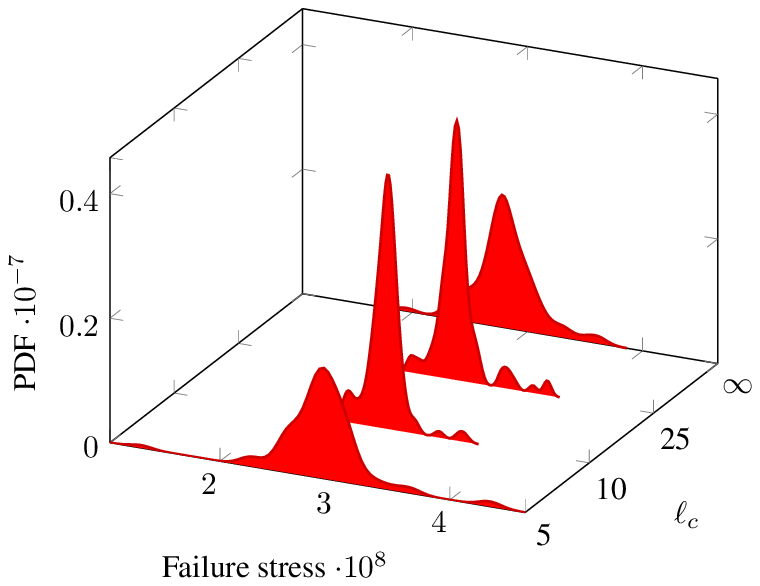}
\includegraphics[width=0.4\textwidth]{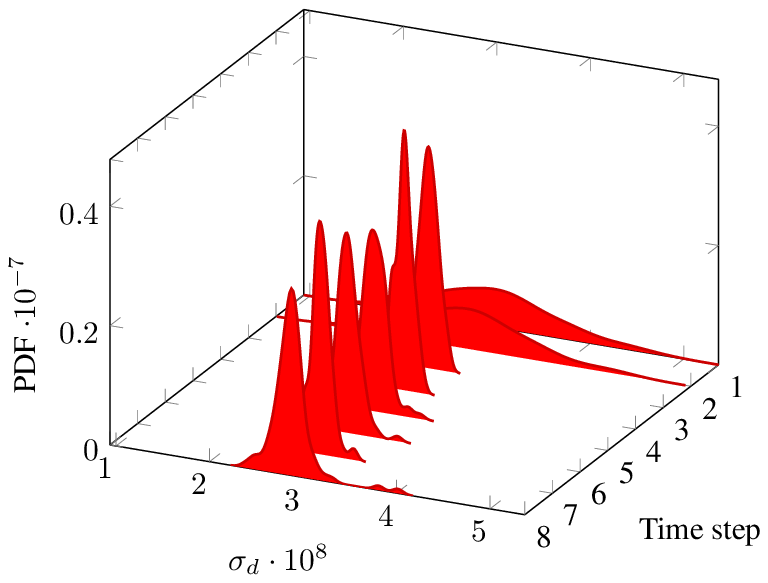}
\includegraphics[width=0.4\textwidth]{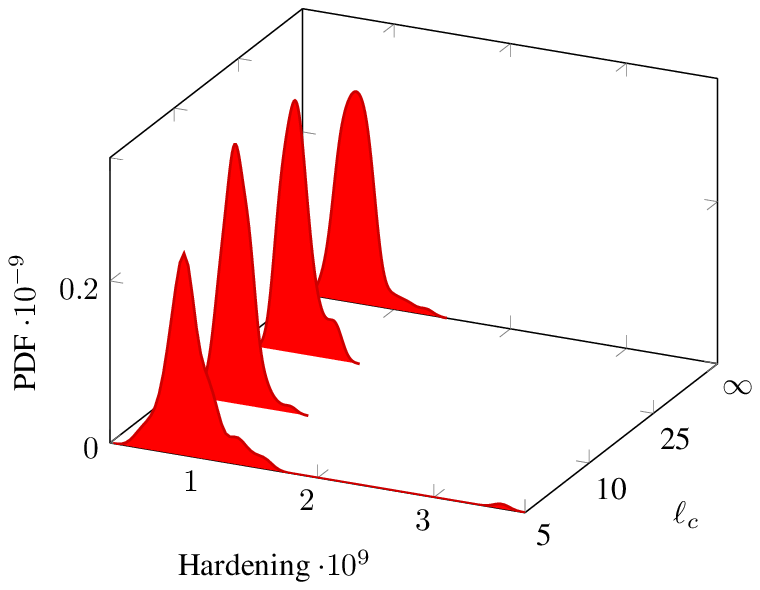}
\includegraphics[width=0.4\textwidth]{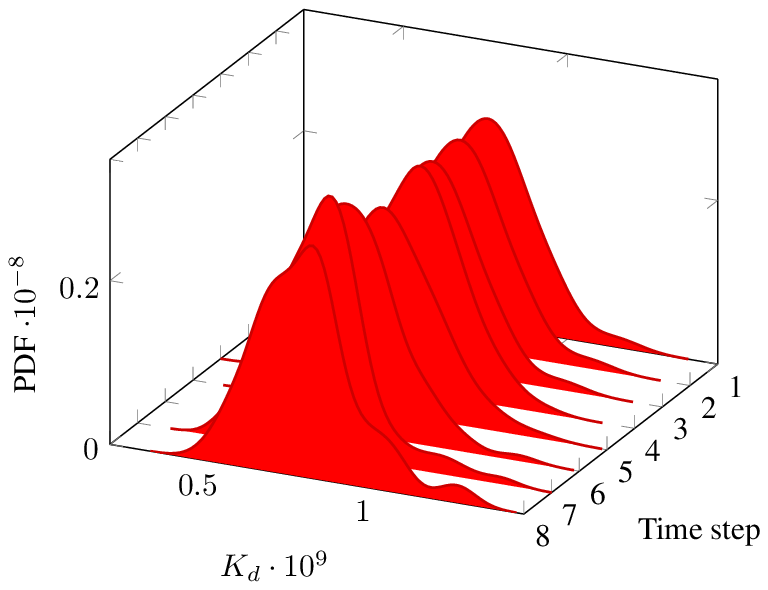}
\caption{The PDF of parameters $\kappa$, Failure stress $\sigma_d$ and Hardening modulus $K_d$ w.r.t. the correlation length (left) $\ell_c$ and time steps (right)}
\label{fig_pure_shear_ran_pos_PR1}
\end{center}
\end{figure}

The up-scaled parameter estimates behave similarly to the energy estimates as shown in Fig.~\ref{fig_pure_shear_ran_pos_PR1}. The hardening parameter doesnt get updated, and stays constant over time. Similar holds for the correlation length as the hardening energy doesnt change with the correlation length.

\begin{table}[h!]
  \centering  
  \begin{tabular}{cccc}
    \toprule
    Property&$\kappa$ & $\sigma_{f}$& $K_{d}$\\
    \midrule
    $\mu$    &204440 &300&450  \\
    $\sigma$ &10222 &15 &22.5 \\
    \bottomrule
  \end{tabular}
  \caption{Fine-scale statistics for elasto-damage constitutive model (in MPa)}
   \label{tab:damageprior}
\end{table}

%% file: conclusion.tex
\section{Conclusion}

The stochastic multiscale analysis as previously presented is one praticular kind of inverse problem in which the coarse-scale 
parameters are to be estimated given the fine-scale information. In this paper we employed the energy conservation principle in order to 
estimate the macro-scale parameters given meso-scale energy descriptions. Such an approach is then allowing derivation of new constitutive laws 
on the macroscale counterpart, the ones that are optimally matching the energy information. 
Furthermore, we show that 
in a case when the fine scale energy information is of deterministic kind, i.e.
describes the particular RVE, the process of estimation can be easily done by employing the nonlinear version of the Kalman filter. The filter then 
represents the map between the observation and the quantity of interest, i.e. the macro-scale model parameters, or the model itself. In addition, we have shown that this kind of mapping can be also used in a more general situation in which the fine scale information is described by uncertainty. The only requirement to achieve this is to fully specify the random variable representing the data, i.e. to describe its probability distribution. For this purpose we 
employ the Bayes variational inference in combination with the copula theory. Computationally, the measurement probability distribution is then represented by a functional approximation in terms of the polynomial chaos expansion obtained by mapping the measurement data to the Gaussian space, applying an inverse transform and using an additional sparse Bayes variational inference for the purpose of 
estimation of the expansion coefficients.  As the inervse map from the energy space to the Gaussian one is not easy to approximate, we recommend to first 
discretise the energy space (i.e.~sample), and then to map each sample to the macroscale model parameters.  As shown on both linear elastic and elasto-damage examples, the latter ones can be more accurate approximated. Note that in this paper we have only observed the elasto-damage models on two scales under one specific loading condition. However, in practice this will not suffice to achieve good macro-scale representation. Therefore, the next step to be considered is to add different loading conditions into estimation.

\section*{Acknowledgements} 
We acknowledge the support provided by the German Science Foundation (Deutsche Forschungsgemeinschaft, DFG) as part of priority programs SPP 1886 and SPP 1748.